\documentclass[aos,preprint]{imsart}
\RequirePackage{natbib}

\usepackage{wrapfig}

\RequirePackage[OT1]{fontenc}
\RequirePackage{amsthm,amsmath,amsfonts}
\RequirePackage[colorlinks,citecolor=blue,urlcolor=blue]{hyperref}
\usepackage{graphicx}
\usepackage{xr}
\usepackage{xcite}
\usepackage{scrextend} 

\externalcitedocument{supp_final}
\startlocaldefs

\newcommand{\R}{\mathbb{R}}

\renewcommand{\hat}[1]{\widehat{#1}}

\newcommand{\D}{\mathcal{D}}

\newcommand{\F}{\mathbb{F}}

\newcommand{\X}{\mathcal{X}}

\newcommand{\tsum}{\textstyle\sum}

\renewcommand{\P}{\mathbb{P}}
\newcommand{\E}{\mathbb{E}}

\renewcommand{\L}{\mathcal{L}}

\newcommand{\e}{\epsilon}
\newcommand{\ve}{\varepsilon}

\newcommand{\var}{\operatorname{var}}

\newcommand{\vol}{\operatorname{vol}}
\newcommand{\ts}{\textstyle}

\newcommand{\ttop}{^{\top}}

\newcommand{\Err}{\textsc{Err}}

\newcommand{\mnorm}[1]{\left\vert\kern-1.5pt\left\vert\kern-1.5pt\left\vert #1\right\vert\kern-1.5pt\right\vert\kern-1.5pt\right\vert}

\newcommand{\id}{\text{id}}
\renewcommand{\div}{\text{div}}

\newcommand{\err}{\text{err}}
\numberwithin{equation}{section}
\theoremstyle{plain}

\newtheorem{assumption}{Assumption}
\newtheorem{definition}{Definition}
\usepackage{natbib}

\endlocaldefs

\numberwithin{equation}{section}
\theoremstyle{plain}
\newtheorem{thm}{Theorem}[section]
\newtheorem{proposition}{Proposition}[section]

\newtheorem{lemma}{Lemma}[section]

\begin{document}

\begin{frontmatter}
\title{Estimating the Algorithmic Variance of Randomized Ensembles via the Bootstrap}
\runtitle{Algorithmic Variance of Randomized Ensembles}

\begin{aug}
\author{\fnms{Miles E.} \snm{Lopes}\thanksref{t1}\ead[label=e1]{melopes@ucdavis.edu}}

\thankstext{t1}{This research was partially supported by NSF grant DMS 1613218.}
\runauthor{M. E. Lopes}

\affiliation{University of California, Davis}

\address{Miles E. Lopes\\
Department of Statistics\\
 Mathematical Sciences Building 4118\\
  399 Crocker Lane\\
University of California\\
   Davis One Shields Avenue\\
Davis, CA 95616\\
\printead{e1}\\
\phantom{E-mail:\ }\printead*{}}

\end{aug}

\begin{abstract}
 Although the methods of  bagging and random forests are some of the most widely used prediction methods, relatively little is known about their algorithmic convergence. 
 In particular, there are not many theoretical guarantees for deciding when an ensemble  is ``large enough'' 
 --- so that its accuracy is close to that of an ideal infinite ensemble.
  Due to the fact that bagging and random forests are randomized algorithms, the choice of  ensemble size is closely related to the notion of ``algorithmic variance'' (i.e.~the variance of prediction error due only to the training algorithm). In the present work, we propose a bootstrap method to estimate this variance for bagging, random forests, and related methods in the context of classification.
To be specific, suppose the training dataset is fixed, and let the random variable $\Err_t$ denote the prediction error of a randomized ensemble of size $t$.
Working under a ``first-order model'' for randomized ensembles, we prove that the centered law of $\Err_t$ can be consistently approximated via the proposed method as $t\to\infty$. Meanwhile, the computational cost of the method is quite modest, by virtue of an extrapolation technique.  As a consequence, the method offers a practical guideline for deciding when the algorithmic fluctuations of $\Err_t$ are negligible.
 
\end{abstract}



\end{frontmatter}

\section{Introduction}\label{sec:intro}

Random forests and bagging are some of the most widely used prediction methods~\citep{breiman1996, breiman2001}, and over the course of the past fifteen years, much progress has been made in analyzing their statistical performance~\citep{buhlmannyu,hallsamworth,biau2008,biau2012,scornetconsistency}. However, from a computational perspective, relatively little is understood about the \emph{algorithmic convergence} of these methods, and in practice, ad hoc criteria are generally used to assess this convergence.

 To clarify the idea of algorithmic convergence, recall that when bagging and random forests are used for classification, a large collection of $t$ randomized classifiers is trained, and then new predictions are made by taking the plurality vote of the classifiers. If such a method is run several times on the same training data $\D$, the prediction error $\Err_t$ of the ensemble will vary with each run, due to the randomized training algorithm. As the ensemble size increases $(t\to\infty)$ with $\D$ held fixed, the random variable $\Err_t$ typically decreases and eventually stabilizes at a limiting value $\err_{\infty}=\err_{\infty}(\D)$. In this way, an ensemble reaches algorithmic convergence when its prediction error nearly  matches that of an infinite ensemble trained on the same data.

Meanwhile, with regard to computational cost, larger ensembles are more expensive to train, to store in memory, and to evaluate on unlabeled points. For this reason, it is desirable to have  a quantitative guarantee that an ensemble of a given size will perform nearly as well as an infinite one.
This type of guarantee also prevents wasted computation, and assures the user that extra classifiers are unlikely to yield much improvement in accuracy.

\subsection{Contributions and related work}\label{sec:contrib} 
To measure algorithmic convergence, we propose a new bootstrap method for approximating the distribution $\mathcal{L}(\sqrt{t}(\Err_t-\err_{\infty})|\D)$ as $t\to\infty$. Such an approximation allows the user to decide when the algorithmic fluctuations of $\Err_t$ around $\err_{\infty}$ are negligible. If particular, if we refer to the \emph{algorithmic variance} 
$$\sigma_t^2:=\var(\Err_t|\D)=\E[\Err_t^2|\D]-(\E[\Err_t|\D])^2,$$
 as the variance of $\Err_t$ due only the training algorithm, then the parameter $\sigma_t$ is a concrete measure of convergence that can be estimated via the bootstrap. 
 In addition, the computational  cost of the method turns out to be quite modest, by virtue of an extrapolation technique, as described in Section~\ref{sec:practical}.

Although the bootstrap is an established approach to distributional approximation and variance estimation, our work applies the bootstrap in a relatively novel way. Namely, the method is based on ``bootstrapping an algorithm'', rather than ``bootstrapping data'' ---  and in essence, we are applying an inferential method in order to serve a computational purpose. The opportunities for applying this perspective to other randomized algorithms can also be seen in the papers \citet{nocedal,lopes_matrix,lopes_LS}, which deal with stochastic gradient methods, as well as randomized versions of matrix multiplication and least-squares.

\paragraph{Bootstrap consistency} From a theoretical standpoint, our main result (Theorem~\ref{THM:CLT}) shows that the proposed method consistently approximates the distribution \smash{$\mathcal{L}(\sqrt{t}(\Err_t-\err_{\infty})|\D)$} as $t\to\infty$ under a ``first-order model'' for randomized ensembles. The proof also offers a couple of theoretical contributions related to Hadamard differentiability and the functional delta method~\citep[Chapter 3.9]{vaartWellner}.
 The first ingredient is  a lifting operator $\boldsymbol L$, which transforms a univariate empirical c.d.f. \smash{$\F_t:[0,1]\to [0,1]$} into a multivariate analogue \smash{$\boldsymbol L(\F_t):\Delta\to\Delta$}, where $\Delta$ is a simplex.
 In addition to having interesting properties in its own right, the lifting operator will allow us to represent $\Err_t$ as a functional of an empirical process. 
The second ingredient is the calculation of this functional's Hadamard derivative, which leads to
 a surprising connection with the classical \emph{first variation formula} for smooth manifolds~\citep{simon,white}.\footnote{Further examples of problems where geometric analysis plays a role in understanding the performance of numerical algorithms may be found in the book~\cite{cucker}.}

To briefly comment on the role of this formula in our analysis, consider the following informal statement of it. Let $\mathcal{M}$ be a smooth $d$-dimensional manifold contained in $\R^d$, and let $\{f_{\delta}\}_{\delta \in(-1,1)}$ be a one-parameter family of diffeomorphisms \smash{$f_{\delta} : \mathcal{M}\to f_{\delta}(\mathcal{M})\subset \R^d$},  satisfying $f_{\delta}\to\id_{\mathcal{M}}$ as $\delta\to 0$, where $\id_{\mathcal{M}}$ denotes the identity map on $\mathcal{M}$. Then,
\begin{equation}\label{firstvarm}
\ts\frac{d}{d\delta} \vol(f_{\delta}(\mathcal{M}))\Big|_{\delta=0} = \displaystyle\int_{\mathcal{M}} \text{div}(Z)(\theta) d\theta,
\end{equation}
where $\vol(\cdot)$ is a volume measure, the symbol $\text{div}(Z)$ denotes the divergence of the vector field $Z(\theta):=\frac{\partial}{\partial \delta} f_{\delta}(\theta) \big|_{\delta=0}$, and the symbol $d\theta$ is a volume element on $\mathcal{M}$.
In our analysis, it is necessary to adapt this result to a situation where the maps $f_{\delta}$ are non-smooth, the manifold $\mathcal{M}$ is a non-smooth subset of Euclidean space, and the vector field $Z(\cdot)$ is a non-smooth Gaussian process. Furthermore, applying a version of Stokes' theorem to the right side of equation~\eqref{firstvarm} leads to a particular linear functional of $Z(\cdot)$, which turns out to be the Hadamard derivative relevant to understanding $\Err_t$. A more detailed explanation of this connection is given below equation~\eqref{maindiff} in Appendix~\ref{sec:highlevel}.

\paragraph{Related work} In the setting of binary classification, a few papers analyze the bias $\E[\Err_t-\err_{\infty}|\D]$, and show that it converges at the fast rate of $1/t$ under various conditions~\citep{ngjordan,lopes2016,cannings2015}.  A couple of other works study alternative measures of convergence. For instance, the paper~\citet{lam1997application} considers the probability that the majority vote commits an error at a fixed test point, and the paper~\citet{hernandez} provides an informal analysis of the probability that an ensemble of size $t$ disagrees with an infinite ensemble
 at a random test point, but these approaches do not directly control $\Err_t$. In addition, some empirical studies of algorithmic convergence may be found in the papers~\citet{oshiro,latinne}.

Among the references just mentioned, the ones that are most closely related to the current paper are~\citet{lopes2016} and~\citet{cannings2015}. These works derive theoretical upper bounds on $\var(\Err_t|\D)$ or $\var(\Err_{t,l}|\D)$, where $\Err_{t,l}$ is the error rate on a particular class $l$ (cf. Section~\ref{sec:practical}). The paper~\citet{lopes2016} also proposes a method to  estimate the unknown parameters in such bounds. In relation to these works, the current paper differs in two  significant ways. First, we offer an approximation to the full distribution $\mathcal{L}(\Err_t-\err_{\infty}|\D)$, and hence provide a \emph{direct estimate} of algorithmic variance, rather than a bound. Second, the method proposed here is relevant to a wider range of problems, since it can handle any number of classes, whereas the analyses in~\citet{lopes2016} and~\citet{cannings2015} are specialized to the binary setting. Moreover, the theoretical analysis of the bootstrap approach is entirely different from the previous techniques used in deriving variance bounds.

Outside of the setting of randomized ensemble classifiers, the papers \citet{sexton2009standard,arlot2014,wager2014confidence,mentch2016,scornet2016asymptotics} look at the algorithmic fluctuations of ensemble regression functions at a fixed test point.

\subsection{Background and setup}\label{sec:background}
We consider the general setting of a classification problem with $k\geq 2$ classes. The set of training data is denoted $\mathcal{D}:=\{(X_1,Y_1),\dots,(X_{n},Y_{n})\}$, which is contained in a sample space $\X\times \mathcal{Y}$. The feature space $\mathcal{X}$ is arbitrary, and the space of labels $\mathcal{Y}$ has cardinality $k$.
An ensemble of $t$ classifiers is denoted by $Q_i:\mathcal{X}\to\mathcal{Y}$, with $i=1,\dots,t$.

\noindent \paragraph{Randomized ensembles}  
The key issue in studying the algorithmic convergence of bagging and random forests is randomization. In the method of bagging, randomization is introduced by generating random sets $\D_1^*, \dots,\D_t^*$, each of size $n$, via  sampling with replacement from $\D$. For each $i=1,\dots,t$, a classifier $Q_i$ is trained on $\D_i^*$, with the same classification method being used each time. When each $Q_i$ is trained with a decision tree method (such as CART~\citep{CART}), the random forests procedure extends bagging by adding a randomized feature selection rule~\citep{breiman2001}.

It is helpful to note that the classifiers in bagging and random forests can be represented in a common way. Namely, there is a deterministic function, say $g$, such that for any fixed $x\in\mathcal{X}$, each classifier $Q_i$ can be written as
\begin{equation}\label{rfrep}
Q_i(x)=g(x, \D, \xi_i),
\end{equation}
 where $\xi_1,\xi_2,\dots$, is an i.i.d. sequence of random objects, independent of $\D$, that specify the ``randomizing parameters'' of the classifiers (cf.\,\citet[Definition 1.1]{breiman2001}). For instance, in the case of bagging, the object $\xi_i$ specifies the randomly chosen points in $\D_i^*$.
 
 Beyond bagging and random forests, our proposed method will be generally applicable to ensembles that can be represented in the form~\eqref{rfrep}, such as those in~\citet{tinkamho1998,dietterich2000,buhlmannyu}. This representation should be viewed abstractly, and it is not necessary for the function $g$ or the objects $\xi_1,\xi_2,\dots$ to be explicitly constructed in practice. Some further examples include  a recent ensemble method based on random projections~\citep{cannings2015}, as well as the voting Gibbs classifier~\citep{ngjordan}, which is a Bayesian ensemble method based on posterior sampling.
 More generally, if the functions $Q_1,Q_2,\dots$ are i.i.d.~conditionally on $\D$, then the ensemble can be represented in the form~\eqref{rfrep}, as long as the classifiers lie in a standard Borel space~\cite[Lemma 3.22]{Kallenberg}. 
 Lastly, it is important to note that the representation~\eqref{rfrep} generally does not hold for classifiers generated by boosting methods~\citep{schapire2012boosting}, for which the analysis of algorithmic convergence is quite different.

\noindent \paragraph{Plurality vote} For any $x\in\mathcal{X}$, we define the ensemble's \emph{plurality vote} as the label receiving the largest number of votes among $Q_1(x),\dots,Q_t(x)$. In the exceptional case of a tie, it will simplify technical matters to define the plurality vote as a symbol not contained in $\mathcal{Y}$, so that a tie always counts as an error. We also use the labeling scheme, $\mathcal{Y}:=\{\boldsymbol e_0,\dots,\boldsymbol e_{k-1}\}\subset\R^{k-1}$, where $\boldsymbol e_0:=\boldsymbol 0$,  and $\boldsymbol e_l$ is the $l$th standard basis vector for $l\geq 1$. One benefit of this scheme is that the plurality vote  is determined by the average of the labels, $\bar{Q}_t(x):=\ts\frac{1}{t}\sum_{i=1}^t Q_i(x)$. For this reason, we denote the plurality vote as ${\tt{V}}(\bar{Q}_t(x))$.

\noindent \paragraph{Error rate} Let $\nu=\mathcal{L}(X,Y)$ denote the distribution of a test point $(X,Y)$ in $\mathcal{X}\times \mathcal{Y}$, drawn independently of $\D$ and $Q_1,\dots,Q_t$.
Then, for a particular realization of the classifiers $Q_1,\dots,Q_t$, trained with the given set $\D$, the prediction error rate is defined as
\begin{equation}\label{errdef}
\Err_{t}:=\int_{\mathcal{X}\times \mathcal{Y}} 1\{{\tt{V}}(\bar{Q}_t(x))\neq y\}d\nu(x,y)  \ = \ \P\Big({\tt{V}}(\bar{Q}_t(X)) \neq Y \,  \big| \,   \D, \boldsymbol \xi_t\Big),
\end{equation}
where $\boldsymbol \xi_t:=(\xi_1,\dots,\xi_t)$. (Class-wise error rates $\Err_{t,l}$, with $l=0,\dots,k-1$ will also be addressed in Section~\ref{sec:classwise}.)
Here, it is  crucial to note that $\Err_{t}$ is a random variable, since $\bar{Q}_t$ is a random function. Indeed, the integral above shows that $\Err_t$ is a functional of $\bar{Q}_t$. Moreover, there are two sources of randomness to consider: the algorithmic randomness arising from $\boldsymbol \xi_t$, and the randomness arising from the training set $\mathcal{D}$. Going forward, we will focus on the algorithmic fluctuations of $\Err_{t}$ due to $\boldsymbol \xi_t$, and our analysis will always be conditional on $\D$.

\noindent \paragraph{Algorithmic variance}  At first sight, it might not be obvious how to interpret the algorithmic fluctuations of $\Err_t$ when $\D$ is held fixed. These fluctuations are illustrated below in Figure~\ref{fig:intro}.
 The left panel shows how $\Err_{t}$ changes as decision trees are added incrementally during a single run of the random forests algorithm.  For the purposes of illustration, if we run the algorithm repeatedly on $\D$ to train many ensembles, we obtain a large number of sample paths of $\Err_{t}$ as a function of $t$, shown in the right panel.
Averaging the sample paths at each value of $t$ produces the red curve, representing $\E[\Err_{t}\big| \D]$ with $\boldsymbol \xi_t$ averaged out. Furthermore, the
 blue envelope curves for the sample paths are obtained by plotting  \mbox{$\E[\Err_{t}| \D]\pm3\sqrt{\var(\Err_t|\D)}$} as a function of $t$. 

 \begin{figure}[h!]
\centering
{\includegraphics[angle=0,
  width=.44\linewidth,  height=.44\linewidth]{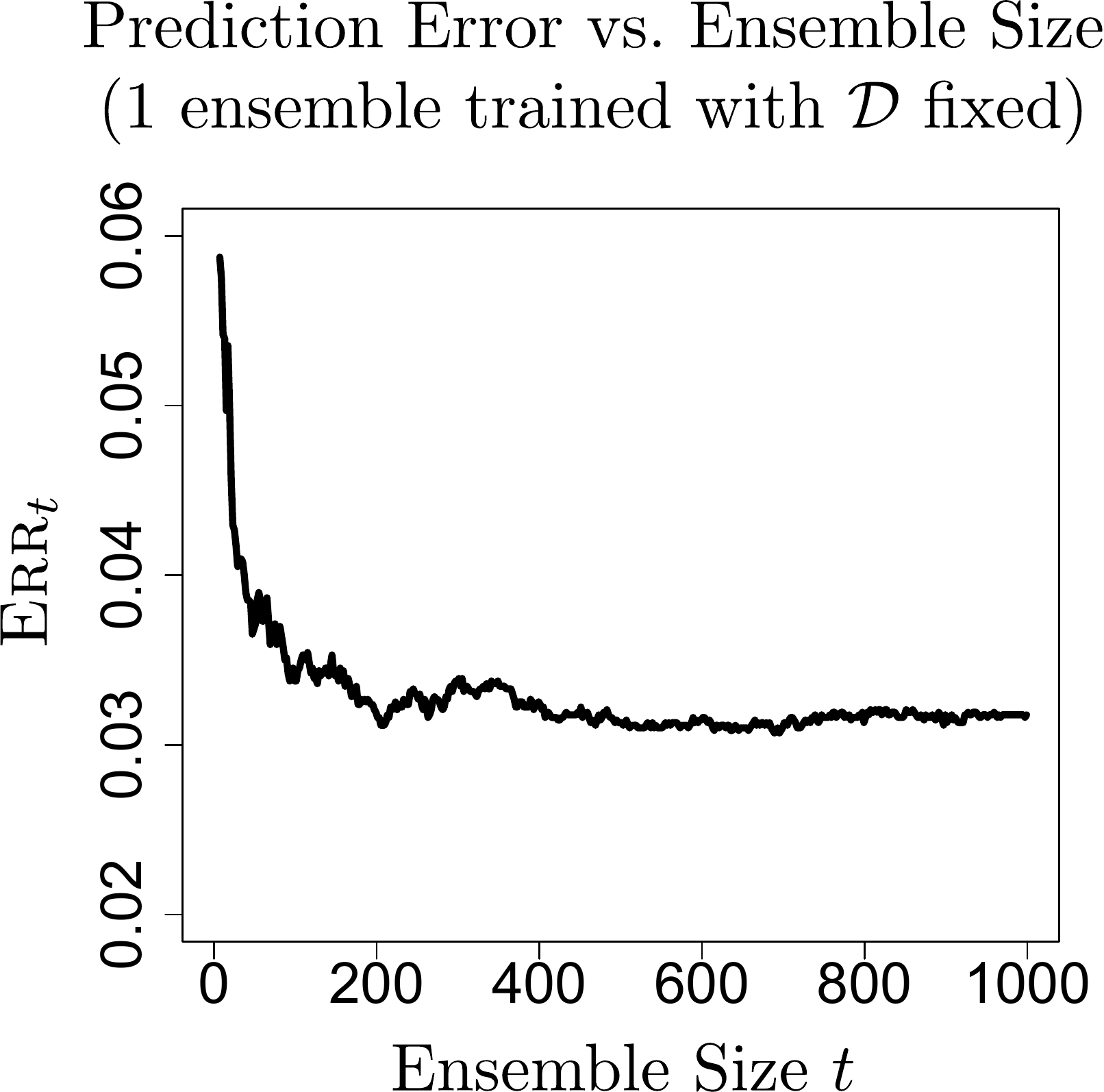}} \ \ \ \
{\includegraphics[angle=0,
  width=.44\linewidth,height=.44\linewidth]{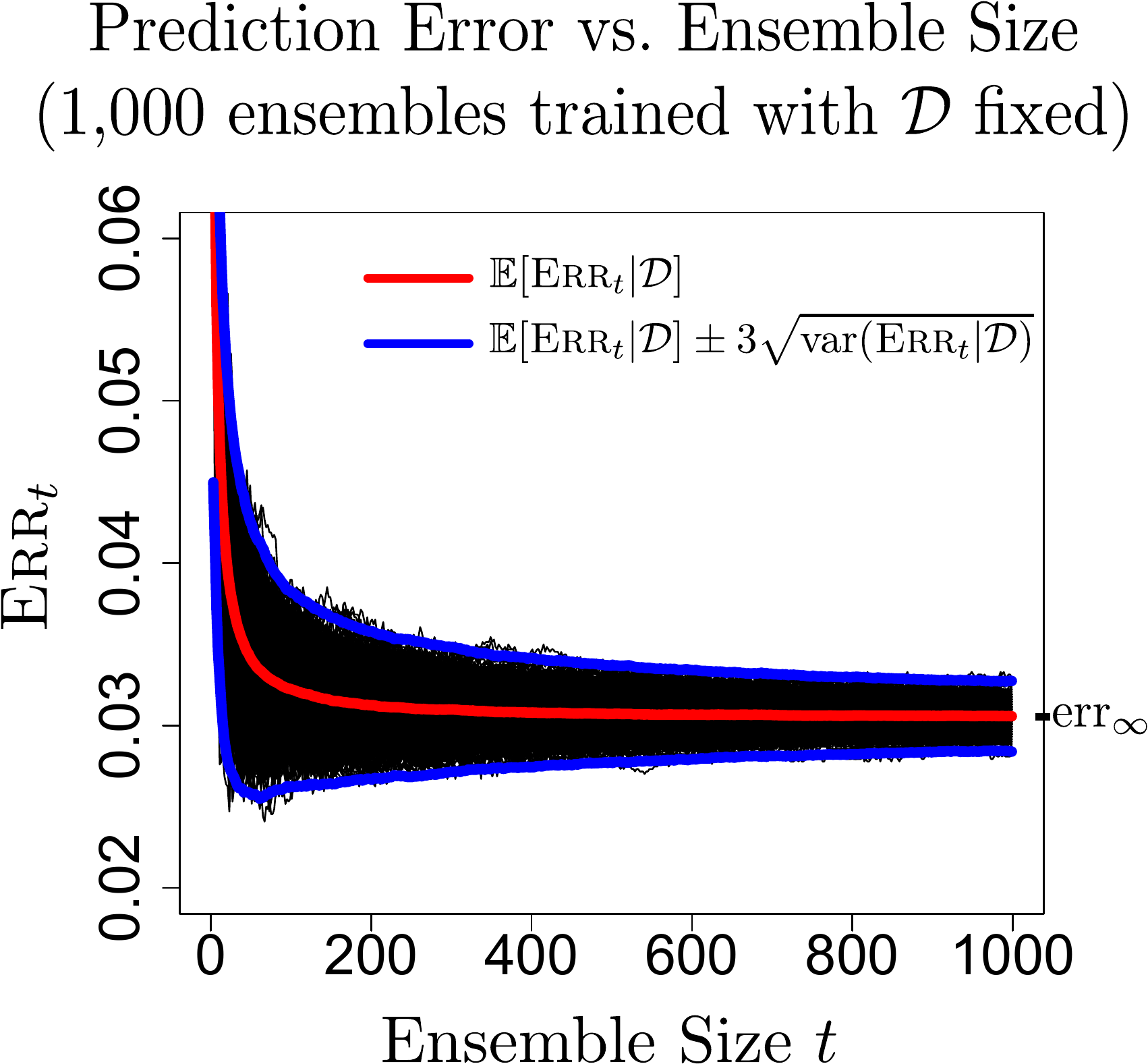}}
 \vspace{-0.1cm}
\caption{ Left panel: The fluctuations of $\Err_{t}$ for a single run of random forests on the ``nursery data'' (cf. Section~\ref{sec:expt}). Right panel: The fluctuations of $\Err_{t}$ for 1,000 runs of random forests on the same data (i.e. $1,\!000$ different realizations of the ensemble $Q_1,\dots,Q_t$). For each ensemble, the value $\Err_t$ was approximated by testing on the set of points denoted $\D_{\textup{ground}}$, as described in Section~\ref{sec:expt}.
}
\label{fig:intro}
\end{figure}

\noindent \paragraph{Problem formulation}\label{probformulation}
Recall that the value $\err_{\infty}=\err_{\infty}(\D)$ represents the ideal prediction error of an infinite ensemble trained on $\D$.  Hence, a natural way of defining algorithmic convergence is to say that it occurs when $t$ is large enough so that the condition \smash{$|\Err_t-\err_{\infty}|\leq \e$} holds with high probability, conditionally on $\D$, for some user-specified tolerance $\e$. However, the immediate problem we face is that it is not obvious how to check such a condition in practice.

From the right panel of Figure~\ref{fig:intro}, we see that for most $t$, the inequality \smash{$|\Err_t-\err_{\infty}|\leq 3\sigma_t$} is highly likely to hold --- and this observation can  be formalized using Theorem~\ref{THM:CLT} later on. For this reason, we propose to estimate $\sigma_t$ as a route to measuring algorithmic convergence. It is also important to note that estimating the quantiles of $\L(\Err_t-\err_{\infty}|\D)$ would serve the same purpose, but for the sake of simplicity, we will focus on $\sigma_t$.
In particular, there are at least two ways that an estimate $\hat{\sigma}_t$ can be used in practice:

\label{probformulation}

\begin{enumerate}
\item \emph{Checking convergence for a given ensemble.} \label{probformulation}
If an ensemble of a given size $t_0$ has been trained, then convergence can be checked by asking whether or not the observable condition $3\hat{\sigma}_{t_0}\leq \e$ holds.  Additional comments on possible choices for $\e$ will be given shortly.

\item \emph{Selecting $t$ dynamically.} In order to make the training process as computationally efficient as possible, it is desirable to select the smallest $t$ needed so that $|\Err_t-\err_{\infty}|\leq \e$ is likely to hold. It turns out that this can be accomplished using an extrapolation technique, due to the fact that $\sigma_t$ tends to scale like $1/\sqrt{t}$ (cf. Theorem~\ref{THM:CLT}). More specifically, if the user trains a small initial ensemble of size $t_0$ and computes an estimate $\hat{\sigma}_{t_0}$, then ``future'' values of $\sigma_t$ for $t\gg t_0$ can be estimated at no additional cost with the re-scaled estimate $\sqrt{t_0/t}\,\hat{\sigma}_{t_0}$.  In other words, it is possible to look ahead and predict how many additional classifiers are needed to achieve $3\sigma_t\leq \e$. Additional details are given in Section~\ref{sec:extrap}.

\end{enumerate}
\paragraph{Sources of difficulty in estimating $\sigma_t$} Having described the basic formulation of the problem,  it is important to identify what challenges are involved in estimating $\sigma_t$.
First, we must keep in mind that the parameter $\sigma_t$ describes how  $\Err_t$ fluctuates over \emph{repeated ensembles} generated from $\D$ --- and so it is not obvious that it is possible to estimate $\sigma_t$  from the output of a \emph{single ensemble}. Second, 
the computational cost to estimate $\sigma_t$ should not outweigh the cost of training the ensemble, and consequently, the proposed method should be computationally efficient. 
These two obstacles will be described in Sections~\ref{sec:consist} and~\ref{sec:extrap} respectively.

\paragraph{Remarks on the choice of error rate}
Instead of analyzing the random variable $\Err_t$, a natural inclination would be to consider the ``unconditional'' error rate \smash{$\E[\Err_t]=\P({\tt{V}}(\bar{Q}_t(X))\neq Y))$}, which is likely more familiar, and  reflects averaging over both $\D$ and the randomized algorithm. Nevertheless, there are a few reasons why the ``conditional'' error $\Err_t$ may be more suitable for analyzing algorithmic convergence. First, the notion of algorithmic convergence typically describes how much computation should be applied to a given input  --- and in our context, the given input for the training algorithm is $\D$. Second, from an operational standpoint, once a user has trained an ensemble on a given dataset, their \emph{actual} probability of misclassifying a future test point is $\Err_t$, rather than $\E[\Err_t]$.

There is also a sense in which the convergence of $\E[\Err_t]$ may be misleading.  Existing theoretical results suggest that  $\E[\Err_t-\err_{\infty}]$ converges at the fast rate of $1/t$, and in the binary case, \smash{$k=2$}, it can be shown that \smash{$\E[\Err_t-\err_{\infty}|\D]=\ts\frac{1}{t}c(\D)+o(\ts\frac{1}{t})$,}  for some number $c(\D)$, under certain conditions~\citep{lopes2016,cannings2015}. However, in Theorem 1 of Section~\ref{sec:main}, we show that conditionally on $\D$, the difference $\Err_t-\err_{\infty}$ has fluctuations of order $1/\sqrt{t}$. 
 In this sense, if the choice of $t$ is guided only by $\E[\Err_t]$ (rather than the fluctuations of $\Err_t$), then the user may be misled into thinking that algorithmic convergence occurs much faster than it really does --- and this distinction is apparent from the red and blue curves in Figure~\ref{fig:intro}.

\paragraph{Outline}Our proposed bootstrap method is described in Section~\ref{sec:method}, and our main consistency result is given in Section~\ref{sec:main}. Practical considerations are discussed in Section~\ref{sec:practical}, numerical experiments are given in Section~\ref{sec:expt}, and conclusions are stated in Section~\ref{sec:conc}. The essential ideas of the proofs are explained in Appendices~\ref{app:mainthm} and~\ref{sec:highlevel}, while the technical arguments are given Appendices~\ref{app:variation}-\ref{lemmas}. Lastly, in Appendix~\ref{app:assumption}, we provide additional assessment of technical assumptions. All appendices are in the supplementary material.

\section{Method}\label{sec:method}
Based on the definition of $\Err_t$ in equation~\eqref{errdef}, we may view $\Err_t$ as a functional of $\bar{Q}_t$, denoted
 \begin{equation*}
 \Err_t = \varphi(\bar{Q}_t).
 \end{equation*}
 From a statistical standpoint, the importance of this expression is that $\Err_t$ is a functional of a sample mean, which makes it plausible that $\sigma_t$ is amenable to bootstrapping, provided that $\varphi$ is sufficiently smooth. 

To describe the bootstrap method, let $(Q_1^*,\dots,Q_t^*)$ denote a random sample with replacement from the trained ensemble $(Q_1,\dots,Q_t)$, and put
$\bar{Q}_t^*(\cdot):=\ts\frac{1}{t}\tsum_{i=1}^t Q_i^*(\cdot).$
%
In turn, it would be natural to regard the quantity
\begin{equation}\label{errtstar}
\Err_t^*:= \varphi(\bar{Q}_t^*)
\end{equation}
as a bootstrap sample of $\Err_t$, but strictly speaking, this is an  ``idealized'' bootstrap sample, because the functional $\varphi$ depends on the unknown test point distribution $\nu=\mathcal{L}(X,Y)$. Likewise, in Section~\ref{hoalg} below, we explain how each value $\varphi(\bar{Q}_t^*)$ can be estimated. So, in other words, if $\hat{\varphi}$ denotes an estimate of $\varphi$, then an estimate of $\Err_t$ would be written as
\begin{equation*}
\hat{\Err}_t:=\hat{\varphi}(\bar{Q}_t),
\end{equation*}
and the corresponding bootstrap sample is
\begin{equation*}
\hat{\Err}_t^*:=\hat{\varphi}(\bar{Q}_t^*).
\end{equation*}
Altogether, a basic version of the proposed bootstrap algorithm is summarized as follows. 

\paragraph{Algorithm 1 (Bootstrap estimation of $\sigma_t$)}
~\\
\vspace{-0.2cm}
\hrule
~\\[0.1cm]
\noindent \textbf{For} $b=1,\dots, B \,  \textbf{:}$

\begin{itemize}
\item Sample $t$ classifiers $(Q_1^*,\dots,Q_t^*)$ with replacement from $(Q_1,\dots,Q_t)$.
\item Compute $z_b:= \hat{\varphi}(\bar{Q}_t^*)$.
\end{itemize}
\textbf{Return:} the sample standard deviation of $z_1,\dots,z_B$, denoted $\hat{\sigma}_t$.\\[-0.1cm]
\hrule
~\\
\vspace{-0.3cm}
\noindent 
\paragraph{Remark} While the above algorithm is conceptually simple, it suppresses most of the implementation details, and these are explained below. Also note that in order to approximate quantiles of \smash{$\mathcal{L}(\Err_t-\err_{\infty}|\D)$}, rather than $\sigma_t$, it is only necessary to modify the last step, by returning the desired quantile of the centered values $z_1-\bar{z},\dots,z_B-\bar{z}$, with $\bar{z}=\frac{1}{B}\sum_{b=1}^Bz_b$. 

\subsection{Resampling algorithm with hold-out or ``out-of-bag'' points}\label{hoalg}
Here, we consider a version of Algorithm 1 where $\varphi$ is estimated implicitly with hold-out points, and then later on, we will explain how hold-out points can be avoided using ``out-of-bag'' (\textsc{oob}) points. To begin, suppose we have a hold-out set of size $m$, denoted \smash{$\D_{\text{hold}}:=\{(\tilde{X}_1,\tilde{Y}_1),\dots,(\tilde{X}_m,\tilde{Y}_m)\}$.} Next, consider an array $\tilde{A}$ of size $t\times m$, whose $i$th row $\tilde{\mathbf{a}}_i$ is given by the predicted labels of $Q_i$ on the hold-out points. That is,
\begin{equation}\label{atildedef}
\tilde{\mathbf{a}}_i:= {\Large{\boldsymbol [}} \, Q_i(\tilde{X}_1)\, , \dots, \, Q_i(\tilde{X}_m) \,{\Large{\boldsymbol ]} },
\end{equation}
and
\begin{eqnarray} 
\tilde{A}:= \left[
            \begin{array}{l}
            - \tilde{\mathbf{a}}_1 - \\
            \ \ \  \vdots\\
            - \tilde{\mathbf{a}}_t -
            \end{array}
            \right].
\end{eqnarray}
The estimated error rate is easily computed as a function of this array, i.e.~$\hat{\Err}_t = \hat{\Err}_t(\tilde{A})$. To spell out the details, let $\tilde{A}_j=(Q_1(\tilde{X}_j),\dots,Q_t(\tilde{X}_j)$) denote the $j$th column of $\tilde{A}$, and with a slight abuse of our earlier notation, let ${\tt{V}}(\tilde{A}_j)$ denote the plurality vote of the labels in $\tilde{A}_j$.  Then, the estimated error rate is obtained from simple column-wise operations on $\tilde{A}$,
\begin{equation}\label{errAdef}
 \hat{\Err}_t(\tilde{A}) := \frac{1}{m}\displaystyle\sum_{j=1}^m 1\{{\tt{V}}(\tilde{A}_j)\neq \tilde{Y}_j\}.
\end{equation}
In other words, $\hat{\Err}_t(\tilde{A})$ is just the proportion of columns of $\tilde{A}$ for which plurality vote is incorrect. (Note that equation~\eqref{errAdef} is where $\varphi$ is implicitly estimated.) Finally, since there is a one-to-one correspondence between the rows $\tilde{\mathbf{a}}_i$ and the classifiers $Q_i$, the proposed method is equivalent to resampling the rows $\tilde{\mathbf{a}}_i$, as given below.

\paragraph{Algorithm 2 (Bootstrap estimation of $\sigma_t$ with hold-out points)}
~\\
\vspace{-0.2cm}
\hrule
~\\[0.1cm]
\noindent \textbf{For} $b=1,\dots, B \,  \textbf{:}$

\begin{itemize}
\item Draw a $t\times m$ array $\tilde{A}^*$ whose rows $(\tilde{\mathbf{a}}_1^*,\dots,\tilde{\mathbf{a}}_t^*)$ are sampled with replacement from $(\tilde{\mathbf{a}}_1,\dots,\tilde{\mathbf{a}}_t)$.

\item Compute $z_b:= \hat{\Err}(\tilde{A}^*)$.
\end{itemize}
\textbf{Return:} the sample standard deviation of $z_1,\dots,z_B$, denoted $\hat{\sigma}_t$.\\[0.0cm]

\hrule
~\\[-0.5cm]

\paragraph{Extension to OOB points}  Since the use of a hold-out set is often undesirable in practice, we instead consider \textsc{oob} points --- which are a special feature of bagging and random forests. To briefly review this notion, recall that each classifier $Q_i$ is trained on a set of $n$ points $\mathcal{D}_i^*$ obtained by sampling with replacement from $\mathcal{D}$. Consequently, each set $\D_i^*$ excludes approximately $(1-\ts\frac{1}{n})^n\approx $ 37\% of the points in $\mathcal{D}$, and these excluded points may be used as test points for the particular classifier $Q_i$. If a point $X_j$ is excluded from $\mathcal{D}_i^*$, then we say ``the point $X_j$ is \textsc{oob} for the classifer $Q_i$'', and we write $i \in \textsc{\textsc{oob}}(X_j)$, where the set $\textsc{\textsc{oob}}(X_j)\subset\{1,\dots,t\}$ indexes the classifiers for which $X_j$ is \textsc{oob}.

In this notation, the error estimate $\hat{\Err}_t(\tilde{A})$ in Algorithm 2 can be given an analogous definition in terms of \textsc{oob} points. Define a \emph{new}  $t\times n$ array $A$ whose $i$th row is given by 
\begin{equation*}
\mathbf{a}_i:=\big[ Q_i(X_1),\dots,Q_i(X_n)\big].
\end{equation*}
Next, letting $A_j$ be the $j$th column of $A$, define ${{\tt{V}}}_{\textsc{o}}(A_j)$ to be the plurality vote of the set of labels $\{Q_i(X_j) \, |\, i \in \textsc{\textsc{oob}}(X_j)\}$. (If this set of labels is empty, then we treat this case as a tie, but this is unimportant, since it occurs with probability $[1-(1-\ts\frac{1}{n})^{n}]^t \approx (0.63)^t$.)
So, by analogy with $\hat{\Err}_{t}(\tilde{A})$, we define
\begin{equation}\label{errOOB}
 \hat{\Err}_{t,\textsc{o}}(A) :=\frac{1}{n}\sum_{j=1}^n 1\{{\tt{V}}_{\textsc{o}}(A_j)\neq Y_j\}.
\end{equation}
Hence, the \textsc{oob} version of Algorithm 2 may be implemented by simply interchanging  $\tilde A$ and $A$, as well as $\hat{\Err}_{t}(\tilde{A}^*)$ and  $\hat{\Err}_{t,\textsc{o}}(A^*)$. The essential point to notice is that the sum in equation~\eqref{errOOB} is now over the training points in $\mathcal{D}$, rather than over the hold-out set $\D_{\text{hold}}$, as in equation~\eqref{errAdef}. 

\section{Main result}\label{sec:main}

Our main theoretical goal is to prove that the bootstrap yields a consistent approximation of $\mathcal{L}(\sqrt{t}(\Err_t-\err_{\infty})|\D)$ as $t$ becomes large. Toward this goal, we will rely on two simplifications that are customary in analyses of bootstrap and ensemble methods. First, we will exclude the Monte-Carlo error arising from the finite number of $B$ bootstrap replicates, as well as the error arising from the estimation of $\Err_t$. For this reason, our results do not formally require the training or hold-out points to be i.i.d. copies of the test point $(X,Y)$ --- but from a practical standpoint, it is natural to expect that this type of condition should hold in order for Algorithm 2 (or its \textsc{oob} version) to work well.

Second, we will analyze a simplified type of ensemble, which we will refer to as a \emph{first-order model}. This type of approach has been useful in gaining theoretical insights into the behavior of complex ensemble methods in a variety of previous works~\citet{biau2008,biau2012,arlot2014,linjeon,genuer,scornet2016asymptotics,scornetIEEE}. In our context, the value of this simplification is that it neatly packages the complexity of the base classifiers, and clarifies the relationship between $t$ and quality of the bootstrap approximation. Also, even with such simplifications, the theoretical problem of proving bootstrap consistency still leads to considerable technical challenges. Lastly, it is important to clarify that the first-order model is introduced only for theoretical analysis, and our proposed method does not rely on this model.

\subsection{A first-order model for randomized ensembles}  Any randomized classifier  $Q_1:\mathcal{X}\to  \{\boldsymbol e_0,\dots,\boldsymbol e_{k-1}\}$ may be viewed as a stochastic process indexed by $\mathcal{X}$. From this viewpoint, we say that another randomized classifier \smash{$T_1:\mathcal{X}\to \{\boldsymbol e_0,\dots,\boldsymbol e_{k-1}\}$} is a \emph{first-order model} for $Q_1$ if it has the same marginal distributions as $Q_1$, conditionally on $\D$, which means
\begin{equation}\label{firsto}
\mathcal{L}(Q_1(x)|\mathcal{D}) = \ \mathcal{L}(T_1(x)|\D) \ \textup{ for all } x\in\mathcal{X}.
\end{equation}
Since $Q_1(x)$ takes values in the finite set of binary vectors $\{\boldsymbol e_0,\dots,\boldsymbol e_{k-1}\}$, the condition~\eqref{firsto} is equivalent to
\begin{equation}\label{firste}
\E[Q_1(x)|\D] \ = \ \E[T_1(x)|\D] \ \ \textup{ for all } x\in\mathcal{X},
\end{equation}
where the expectation is only over the algorithmic randomness in $Q_1$ and $T_1$.
A notable consequence of this matching condition is that the ensembles associated with $Q_1$ and $T_1$
have the \emph{same} error rates on average.
Indeed, if we let  $\Err_{t}'$ be the error rate associated with an ensemble of $t$ independent copies of $T_1$, then it turns out that
\begin{equation}\label{sameerr}
\begin{split}
\E[\Err_{t}|\D]= \E[\Err_{t}'|\D],
\end{split}
\end{equation}
for all $t\geq 1$, where $\Err_t$ is the error rate for $Q_1,\dots,Q_t$, as before. (A short proof is given in Appendix~\ref{lemmas}.)
In this sense, a first-order model $T_1$ is a meaningful proxy for $Q_1$ with regard to statistical performance --- even though the internal mechanisms of $T_1$ may be simpler.

\subsubsection{Constructing a first-order model}

Having stated some basic properties that are satisfied by any first-order model, we now construct a particular version that is amenable to analysis. Interestingly, it is possible to start with an \emph{arbitrary} random classifier $Q_1:\mathcal{X}\to\{\boldsymbol e_0,\dots, \boldsymbol e_{k-1}\}$, and construct an associated $T_1$ in a relatively explicit way. 

To do this, let $x\in\mathcal{X}$ be fixed, and consider the function
\begin{equation}\label{varthetadef}
\vartheta(x):=\E[Q_1(x)|\D],
\end{equation}
which takes values in the ``full-dimensional'' simplex $\Delta\subset\R^{k-1}$, defined by

\begin{equation*}
\Delta :=\Big\{\theta\in [0,1]^{k-1}\Big| \ \theta_1+\dots+\theta_{k-1}\leq 1 \Big\}.
\end{equation*}
For any fixed $\theta\in\Delta$, there is an associated partition of the unit interval into sub-intervals  
\begin{equation*}
I_0(\theta)\cup \dots \cup I_{k-1}(\theta) = [0,1],
\end{equation*}
such that the width of interval $I_l(\theta)$ is equal to $\theta_l$ for $l\geq 1$. Namely, we put $I_1(\theta):=[0,\theta_1]$, and for $l=2,\dots,k-1$,
\begin{equation*}
I_l(\theta) := \Big((\theta_1+\cdots+\theta_{l-1})\, , \, (\theta_1+\cdots+\theta_{l})\Big].
\end{equation*}
Lastly, for $I_0$, we put
\begin{equation*}
I_0(\theta):=\big(\tsum_{l=1}^{k-1}\theta_l \, , \, 1\big].
\end{equation*}
Now, if we let $x\in\mathcal{X}$ be fixed, and let $U_1\sim $ Uniform$[0,1]$, then we define $T_1(x)\in \{\boldsymbol e_0,\dots,\boldsymbol e_{k-1}\}$ to have its $l$th coordinate equal to the following indicator variable
\begin{equation*}
[T_1(x)]_l:=1\big\{U_1\in I_l(\vartheta(x))\big\},
\end{equation*}
 where $l=1,\dots,k-1$. It is simple to check that the first-order matching condition~\eqref{firste} holds, and so $T_1$ is indeed a first-order model of $Q_1$.
Furthermore, given that $T_1$ is defined in terms of a single random variable $U_1\sim$ Uniform$[0,1]$, we obtain a corresponding ``first-order ensemble''  $T_1,\dots,T_t$ via an i.i.d.~sample of uniform variables $U_1,\dots,U_t$, which are independent of $\D$. (The $l$th coordinate of the $i$th classifier $T_i$ is given by $[T_i(x)]_l=1\{U_i\in I_l(\vartheta(x))\}$.) Hence, with regard to the representation $Q_i(x)=g(x,\D,\xi_i)$ in equation~\eqref{rfrep}, we may make the identification\label{thetaiui}
$$\xi_i=U_i,$$
when the first-order model holds with $Q_i=T_i$.

\paragraph{Remark} To mention a couple of clarifications, the variables \smash{$U_1,\dots,U_t$} are only used for the construction of a first-order model, and they play no role in the proposed method. Also, even though the ``randomizing parameters'' $U_1,\dots,U_t$ are independent of $\D$, the classifiers $T_1,\dots,T_t$ still depend on $\D$ through the function $\vartheta(x)=\E[Q_1(x)|\D]$.

\paragraph{Interpretation of first-order model}
 To understand the statistical meaning of the first-order model, it is instructive to consider the simplest case of binary classification, $k=2$. In this case, $T_1(x)$ is a Bernoulli random variable, where $T_1(x)=1\{U_1\leq \vartheta(x)\}$. Since $\bar{Q}_t(x)\to \vartheta(x)$ almost surely as \smash{$t\to\infty$} (conditionally on $\D$), the majority vote of an infinite ensemble has a similar form, i.e. $1\{\ts\frac{1}{2}\leq \vartheta(x)\}$. Hence, the classifiers $\{T_i\}$ can be viewed as ``random perturbations'' of the asymptotic majority vote arising from $\{Q_i\}$. Furthermore, if we view the number $\vartheta(x)$ as score to be compared with a threshold, then the variable $U_i$ plays the role of a random threshold whose expected value is $\frac{1}{2}$. Lastly, even though the formula $T_1(x)=1\{U_1\leq \vartheta(x)\}$ might seem to yield a simplistic classifier, the complexity of $T_1$ is actually wrapped up in the function $\vartheta$.
Indeed, the matching condition~\eqref{firste} allows for the function $\vartheta$ to be arbitrary.

\subsection{Bootstrap consistency}\label{sec:consist}
We now state our main result, which asserts that the bootstrap ``works'' under the first-order model.
To give meaning to bootstrap consistency, we first review the notion of conditional weak convergence.

\paragraph{Conditional weak convergence} 
Let $\lambda_0$ be a probability distribution on $\R$, and let $\{\lambda_{\boldsymbol \xi_t}\}_{t\geq 1}$ be a sequence of probability distributions on $\R$ that depend on the randomizing parameters \smash{$\boldsymbol \xi_t=(\xi_1,\dots,\xi_t)$}. Also, let $d_{\text{BL}}$ be the bounded Lipschitz metric for distributions on $\R$~\cite[Sec 1.12]{vaartWellner}, and let $\P_{\boldsymbol \xi}$ be the joint distribution of $(\xi_1,\xi_2,\dots)$. Then, as $t\to\infty$, we say that $\lambda_{\boldsymbol \xi_t}\xrightarrow{ \  w \  } \lambda_0$ in $\P_{\boldsymbol \xi}$-probability if the sequence $\{d_{\text{BL}}(\lambda_{\boldsymbol \xi_t},\lambda_0)\}_{t\geq 1}$  converges to 0 in  $\P_{\boldsymbol \xi}$-probability.

\paragraph{Remark} If a test point $X$ is drawn from class $Y=\boldsymbol e_l$, then  we denote the distribution of the random vector $\vartheta(X)$, conditionally on $\D$, as
\begin{equation*}
 \mu_l:=\mathcal{L}(\vartheta(X)|\D, Y=\boldsymbol e_l),
\end{equation*} 
which is a distribution on the simplex $\Delta\subset \R^{k-1}$. Since this distribution plays an important role in our analysis, it is worth noting that the properties of $\mu_l$ are \emph{not affected} by the assumption of a first-order model, since $\vartheta(x)=\E[T_1(x)|\D]=\E[Q_1(x)|\D]$ for all $x\in \mathcal{X}$. We will also assume that the measures $\mu_l$ satisfy the following extra regularity condition.

\begin{assumption}\label{ASSUMPTION}
For the given set $\D$, and each $l=0,\dots,k-1$, the distribution $\mu_l$ has a density $f_l:\Delta\to [0,\infty)$ with respect to Lebesgue measure on $\Delta$, and $f_l$ is continuous on $\Delta$. Also, if $\Delta^{\circ}$ denotes the interior of $\Delta$, then for each $l$,  the density $f_l$ is $C^1$ on $\Delta^{\circ}$, and $\|\nabla f_l\|_2$ is bounded on $\Delta^{\circ}$.
\end{assumption}

To interpret this assumption, consider a situation where the class-wise test point distributions $\mathcal{L}(X|Y=\boldsymbol e_l)$ have smooth densities on $\mathcal{X}\subset\R^p$ with respect to Lebesgue measure. In this case, the density $f_l$ will exist as long as $\vartheta$ is sufficiently smooth (cf. Appendix~\ref{app:assumption}, Proposition~\ref{prop:assumption}).  Still, Assumption~\ref{ASSUMPTION} might seem unrealistic in the context of random forests, because $\vartheta$ is obtained by averaging over all decision trees that can be generated from $\D$, and strictly speaking, this is a finite average of non-smooth functions. 
 However, due to the bagging mechanism in random forests, the space of trees that can be generated from $\D$ is very large, and consequently, the function $\vartheta$ represents a very fine-grained average. Indeed, the idea that bagging is actually a ``smoothing operation'' on non-smooth functions has received growing attention over the years~\citep{bujastuetzle2000,buhlmannyu,bujastuetzle2006,efron2014}, and the recent paper \citet{efron2014} states that bagging is ``also known as bootstrap smoothing''. In Appendix~\ref{app:assumption}, we provide additional assessment of Assumption~\ref{ASSUMPTION}, in terms of both theoretical and empirical examples.
\begin{thm}[Bootstrap consistency]\label{THM:CLT} \hspace{-.2cm}Suppose that the first-order model $Q_i=T_i$ holds for all $i\geq 1$, and that Assumption~\ref{ASSUMPTION} holds.
Then, for the given set $\D$, there are numbers $\textup{err}_{\infty}=\textup{err}_{\infty}(\D)$ and $\sigma=\sigma(\D)$ such that as $t\to\infty$,
\begin{equation}\label{basicclt}
\mathcal{L}\big(\sqrt{t}(\textup{Err}_t-\textup{err}_{\infty})\big|\D\big)
\xrightarrow{ \ \ w \ \ } N(0,\sigma^2),
\end{equation}
and furthermore, 
\begin{equation*}
\mathcal{L}\big(\sqrt{t}(\textup{Err}^*_t-\textup{Err}_t)\big| \D, \boldsymbol \xi_t)
\xrightarrow{ \  \ w \ \ }  N(0,\sigma^2) \text{  \   in  \ $\P_{\boldsymbol \xi}$-probability.}
\end{equation*}
\end{thm}
\paragraph{Remarks} 
In a nutshell, the proof of Theorem~\ref{THM:CLT} is composed of three pieces: showing that $\Err_t$ can be represented as a functional of an empirical process (Appendix~\ref{sec:lift}), establishing the smoothness of this functional (Appendix~\ref{sec:hadamard}), and employing the functional delta method (Appendix~\ref{sec:delta}).
With regard to theoretical techniques, there are two novel aspects of the proof. The problem of deriving this functional is solved by introducing a certain lifting operator, while the problem of showing smoothness is based on a non-smooth instance of the first-variation formula, as well as some special properties of Bernstein polynomials. Lastly, it is worth mentioning that the core technical result of the paper is Theorem~\ref{THM:HADAMARD}.

To mention some consequences of Theorem~\ref{THM:CLT}, the fact that the limiting distribution of $\mathcal{L}(\Err_t-\err_{\infty} \big| \D)$ has mean 0 shows that the fluctuations of $\Err_t$ have more influence on algorithmic convergence than the bias $\E[\Err_t-\err_{\infty}\big| \D]$  (as illustrated in Figure~\ref{fig:intro}). Second, the limiting distribution motivates a convergence criterion of the form $3\sigma_t\leq \e$, since asymptotic normality it indicates that the event \smash{$|\Err_t-\err_{\infty}|\leq 3\sigma_t$} should occur with high probability when $t$ is large, and again, this is apparent in Figure~\ref{fig:intro}. Lastly, the theorem implies that the quantiles of $\mathcal{L}(\Err_t-\err_{\infty}\big| \D)$ agree asymptotically with their bootstrap counterparts. This is of interest, because quantiles allow the user to specify a bound on  $\Err_t-\err_{\infty}$ that holds with a tunable probability.
Quantiles also provide an alternative route to estimating algorithmic variance, because if $r_t^*$ denotes the interquartile range of $\mathcal{L}(\Err_t^*-\Err_t\big| \D,\boldsymbol \xi_t)$, then the theorem implies
$ \ts\frac{\sqrt{t}}{c}r_t^* \xrightarrow{ \ \ \ } \sigma  \text{   in   $\P_{\boldsymbol \xi}$-probability},$
where $c=\Phi^{-1}(3/4)-\Phi^{-1}(1/4)$.

\section{Practical considerations}~\label{sec:practical}
In this section, we discuss some considerations that arise when the proposed method is used in practice, such as the choice of error rate, the computational cost, and the choice of a stopping criterion for algorithmic convergence.
\subsection{Extension to class-wise error rates}\label{sec:classwise} In some applications, class-wise  error rates may be of greater interest than the total error rate $\Err_t$.  For any $l=0,\dots,k-1$, let \smash{$\nu_l=\mathcal{L}(X|Y=\boldsymbol e_l)$} denote the distribution of the test point $X$ given that it is drawn from class $l$. Then, the error rate on class $l$ is defined as
\begin{equation}\label{classerrdef}
\Err_{t,l}:=\int_{\mathcal{X}} 1\{{\tt{V}}(\bar{Q}_t(x))\neq \boldsymbol e_l\}d\nu_l(x)  \ = \ \P\Big({\tt{V}}(\bar{Q}_t(X)) \neq \boldsymbol e_l \,  \big| \,   \D, \boldsymbol \xi_t, Y=\boldsymbol e_l\Big),
\end{equation}
and the corresponding algorithmic variance is
$$\sigma_{t,l}^2 :=\var(\Err_{t,l}|\D).$$
In order to estimate $\sigma_{t,l}$,  Algorithm~2 can be easily adapted using either hold-out or \textsc{oob} points from a particular class.
Our theoretical analysis also extends immediately to the estimation of $\sigma_{t,l}$ (cf. Section~\ref{sec:lift}).

\subsection{Computational cost and extrapolation}\label{sec:extrap}

A basic observation about Algorithm 2 is that it only relies on the array of predicted labels $\tilde{A}$ (or alternatively $A$). Consequently, the algorithm does not require any re-training of the classifiers. Also, with regard to computing the arrays $\tilde{A}$ or $A$, at least one of these is  \emph{typically computed anyway} when evaluating an ensemble's performance with hold-out or \textsc{oob} points --- and so the cost of obtaining $\tilde{A}$ or $A$ will typically not be an added expense of Algorithm 2.
Thirdly, the algorithm can be parallelized, since the bootstrap replicates can be computed independently.
Lastly,  the cost of the algorithm is \emph{dimension-free} with respect to the feature space $\mathcal{X}$, since all operations are on the arrays $\tilde A$ or $A$, whose sizes do not depend on the number of features.

To measure the cost of Algorithm 2 in terms of floating point operations, it is simple to check that at each iteration $b=1,\dots, B$, the cost of evaluating  $\widehat{\Err}_t(\tilde{A}^*)$ is of order $t\cdot m$, since $\tilde{A}$ has $m$ columns, and each evaluation of the plurality voting function has cost $\mathcal{O}(t)$. Hence, if the arrays $\tilde{A}$ or $A$ are viewed as given, and if $m=\mathcal{O}(n)$, then the cost of Algorithm 2 is $\mathcal{O}(B\cdot t \cdot n)$, for either the hold-out or \textsc{oob} versions. Below, we describe how this cost can be reduced using a basic form of extrapolation~\citep{bickelrichardson,sidi,brezinski}.

\paragraph{Saving on computation with extrapolation} 

To explain the technique of extrapolation, 
the first step produces an inexpensive estimate $\hat{\sigma}_{t_0}$ by applying Algorithm 2 to a small initial ensemble of size $t_0$. The second step then rescales $\hat{\sigma}_{t_0}$ so that it approximates $\sigma_t$ for $t\gg t_0$. 
This rescaling relies on Theorem~\ref{THM:CLT}, which leads to the approximation,
$\sigma_t \approx \ts\frac{\sigma}{\sqrt{t}}$.
Consequently, we define the extrapolated estimate of $\sigma_t$ as
\begin{equation}\label{extrapcondition}
\hat{\sigma}_{t,\text{extrap}}:=\ts\frac{\sqrt{t_0}}{\sqrt{t}}\cdot \hat{\sigma}_{t_0}.
\end{equation}
In turn, if the user desires $3\sigma_t\leq \e$ for some $\e\in(0,1)$, then $t$ should be chosen so that 
\begin{equation}\label{eq:epscriterion}
3\hat{\sigma}_{t,\text{extrap}}\leq \e,
\end{equation}
 which is equivalent to
$t\geq (\ts\frac{3\sqrt{t_0}}{\e}\cdot\hat{\sigma}_{t_0})^2$.

In addition to applying Algorithm 2 to a smaller ensemble, a second computational benefit is that extrapolation allows the user to ``look ahead'' and dynamically determine how much extra computation is needed so that $\sigma_t$ is within a desired range.  In Section~\ref{sec:expt}, some examples are given showing that $\sigma_{1,000}$ can be estimated well via extrapolation when $t_0=200$.

\paragraph{Comparison with the cost of training a random forest} Given that one of the main uses of Algorithm 2 is to control the size of a random forest, one would hope that the cost of Algorithm 2 is less than or similar to the cost of training a single ensemble.
In order to simplify this comparison, suppose that each tree in the ensemble is grown so that all nodes are split into exactly 2 child nodes (except for the leaves), and that all trees are grown to a common depth $d\geq 2$. Furthermore, suppose that $\mathcal{X}\subset \R^p$, and that random forests uses the default rule of randomly selecting from $\lceil \sqrt{p}\rceil$ features during each node split. Under these conditions, it is known that the cost of training a  random forest with $t$ trees via CART is at least of order $t\cdot \sqrt{p}\cdot d\cdot n$~\cite[p.166]{CART}. Additional background on the details of random forests and decision trees may be found in the book~\citet{elements}.

Based on the reasoning just given, the cost of running Algorithm 2 does not exceed the cost of training  $t$ trees, provided that 
$$B =\mathcal{O}(\ts\frac{t}{t_0}\cdot\sqrt{p} \cdot d),$$
where the factor $\frac{t}{t_0}$ arises from the extrapolation speedup described earlier. 
Moreover, with regard to the selection of $B$, our numerical examples in Section~\ref{sec:expt} show that the modest choice $B=50$ allows Algorithm 2 to perform well on a variety of datasets.

\paragraph{The choice of the threshold $\e$}  When using a criterion such as~\eqref{eq:epscriterion} in practice, the choice of the threshold $\e$ will typically be unique to the user's goals. For instance, if the user desires that $\Err_t$ is within 0.5\% of $\err_{\infty}$, then the choice $\e=.005$ would be appropriate. Another option is to choose $\e$ from a relative standpoint, depending on the scale of the error. If the error is high (say $\E[\Err_t|\D]$ is 40\%), then it may not be worth paying a large computational price to ensure that $3\sigma_t$ is less than 0.5\%. Conversely, if $\E[\Err_{t}|\D]$ is 2\%, then it may be worthwhile to train a very large ensemble so that $\sigma_t$ is a fraction of 2\%. In either of these cases, the size of $\sigma_t$ could be controlled in a relative sense by selecting $t$ when $\hat{\sigma}_t\leq \eta \hat{m}_t$, where $\hat{m}_t$ is an estimate of $\E[\Err_t|\D]$ obtained from a hold-out set, and $\eta\in (0,1)$ is a user-specified constant that measures the balance between computational cost and accuracy. But regardless of the user's preference for $\e$, the more basic point to keep in mind is that the proposed method makes it possible for the user to have direct control over the relationship between $t$ and $\sigma_t$, and this type of control has not previously been available.

\section{Numerical Experiments}\label{sec:expt}

To illustrate our proposed method, we describe experiments in which the random forests method is applied to natural and synthetic datasets (6 in total). More specifically, we consider the task of estimating the parameter \mbox{3$\sigma_t=3\sqrt{\var(\Err_t|\D)}$}, as well as 
$3\sigma_{t,l}=3\sqrt{\var(\Err_{t,l}|\D)}$.
Overall, the main purpose of the experiments is to show that the bootstrap can indeed produce accurate estimates of these parameters. A second purpose is to demonstrate the value of the extrapolation technique from Section~\ref{sec:extrap}.

\subsection{Design of experiments}\label{data:design}
 Each of the 6 datasets were partitioned in the following way. First, each dataset was evenly split into a training set $\mathcal{D}$ and a ``ground truth'' set $\mathcal{D}_{\text{ground}}$, with nearly matching class proportions in $\D$ and $\mathcal{D}_{\text{ground}}$. (The reason that a substantial portion of data was set aside for $\mathcal{D}_{\text{ground}}$ was to ensure that  ground truth values of $\sigma_t$ and $\sigma_{t,l}$ could be approximated using this set.) Next, a smaller set $\D_{\text{hold}}\subset\mathcal{D}_{\text{ground}}$ with cardinality satisfying $|\D_{\text{hold}}|/(|\D_{\text{hold}}|+|\mathcal{D}|)\approx 1/6$ was used as the hold-out set for implementing Algorithm 2. As before, the class proportions in $\D_{\text{hold}}$ and $\D$ were nearly matching. The smaller size of $\D_{\text{hold}}$ was chosen to illustrate the performance of the method when hold-out points are limited.

\paragraph{Ground truth values} After preparing $\mathcal{D}$, $\mathcal{D}_{\text{ground}}$, and $\D_{\text{hold}}$, a collection of 1,000 ensembles was trained on $\mathcal{D}$ by repeatedly running the random forests method. Each ensemble contained a total of 1,000 trees, trained under default settings from the package {\tt{randomForest}}~\citep{randomForestsCitation}. Also, we tested each ensemble on $\mathcal{D}_{\text{ground}}$ to approximate a corresponding sample path of $\Err_t$  (like the ones shown in Figure~\ref{fig:intro}).  Next, in order to obtain ``ground truth'' values for $\sigma_t$ with $t=1,\dots,1,\!000$, we used the sample standard deviation of the 1,000 sample paths at each $t$. (Ground truth values for each $\sigma_{t,l}$ were obtained analogously.)

\paragraph{Extrapolated estimates} With regard to our methodology, we applied the hold-out and \textsc{oob} versions of Algorithm 2 to each of the  ensembles --- yielding 1,000 realizations of each type of estimate of $\sigma_t$. In each case, the number of bootstrap replicates was set to $B=50$, and we applied the extrapolation rule, starting from $t_0=200$.  If we let $\hat{\sigma}_{200,\textsc{h}}$ and $\hat{\sigma}_{200,\textsc{o}}$  denote the initial hold-out and \textsc{oob} estimators, then the corresponding extrapolated estimators for $t\geq 200$ are given by
\begin{equation}\label{extrapests}
\hat{\sigma}_{t,\textsc{h},\text{extrap}}:=\ts\frac{\sqrt{200}}{\sqrt{t}} \hat{\sigma}_{200,\textsc{h}} \text{ \  \ \ and \ \ \ } \hat{\sigma}_{t,\textsc{o},\text{extrap}}:=\ts\frac{\sqrt{200}}{\sqrt{t}} \hat{\sigma}_{200,\textsc{o}}.
\end{equation}
Next, as a benchmark, we considered an enhanced version of the hold-out estimator, for which the entire ground truth set $\D_{\text{ground}}$ was
 used in place of $\D_{\text{hold}}$. In other words, this benchmark reflects a situation where a much larger hold-out set is available, and it is referred to as the ``ground estimate'' in the plots. Its value based on $t_0=200$ is
 denoted $\hat{\sigma}_{200,\textsc{g}}$, and for $t\geq 200$, we use
\begin{equation}\label{extrapground}
\hat{\sigma}_{t,\textsc{g},\text{extrap}}:=\ts\frac{\sqrt{200}}{\sqrt{t}} \hat{\sigma}_{200,\textsc{g}}
\end{equation}
to refer to its extrapolated version. Lastly, class-wise versions of all extrapolated estimators were computed in an analogous way.

\subsection{Description of datasets}\label{sec:dataprep}
The following datasets were each partitioned into $\D$, $\D_{\text{hold}}$ and $\D_{\text{ground}}$, as described above.

\paragraph{Census income data} A set of census records for 48,842 people were collected with 14 socioeconomic features (continuous and discrete)~\citep{uci}.  Each record was labeled as 0 or 1, corresponding to low or high income. The proportions of the classes are approximately (.76,.24). As a pre-processing step, we excluded three features corresponding to work-class, occupation, and native country, due to a high proportion of missing values.

\paragraph{Connect-4 data} The observations represent 67,557 board positions in the two-person game ``connect-4''~\citep{uci}. For each position, a list of 42 categorical features are available, and each position is labeled as a draw $l=0$, loss $l=1$, or win $l=2$ for the first player, with the class proportions being approximately $(.10, .25,.65)$.

\paragraph{Nursery data}
 This dataset was prepared from a set of 12,960 applications for admission to a nursery school~\citep{uci}. Each application was associated with a list of 8 (categorical) socioeconomic features.  Originally, each application was labeled as one of five classes, but in order to achieve reasonable label balance, the last three categories were combined. This led to approximate class proportions $(1/3,1/3,1/3)$.

\paragraph{Online news data} A collection of 39,797 news articles from the website \smash{{\tt{mashable.com}}} were associated with 60 features (continuous and discrete). Each article was labeled based on the number of times it was shared: fewer than 1000 shares $(l=0)$, between 1,000 and 5,000 shares ($l=1$), and greater than 5,000 shares ($l=2$), with approximate class proportions $(.28,.59, .13)$.

\paragraph{Synthetic continuous data} Two classes of data points in $\R^{100}$, each of size 10,000, were obtained by drawing samples from the multivariate normal distributions $N(\boldsymbol \mu_0,\Sigma)$ and $N(\boldsymbol \mu_1,\Sigma)$ with a common covariance matrix $\Sigma$. The first mean vector was chosen to be $\boldsymbol \mu_0 = \boldsymbol 0\in\R^{100}$, and the second mean vector was constructed to be a sparse vector in the following way. Specifically, we sampled 10 numbers $(i_1,\dots,i_{10})$ without replacement from $\{1,\dots,100\}$, and the coordinates of $\boldsymbol \mu_1$ indexed by $(i_1,\dots,i_{10})$ were set to the value $.05$ (with all other coordinates were set to 0).  Letting $U\Lambda U\ttop $ denote the spectral decomposition of $\Sigma$, we selected the matrix of eigenvectors $U$ by sampling from the uniform (Haar) distribution on $100\times 100$ orthogonal matrices. The eigenvalues were chosen as $\Lambda=\text{diag}(\frac{1}{1^2},\frac{1}{2^2},\dots, \frac{1}{100^2})$.

\paragraph{Synthetic discrete data} Two classes of data points in $\R^{100}$, each of size 10,000, were obtained by drawing samples from the discrete distributions Multinomial$(N_0,\boldsymbol p_0)$ and Multinomial$(N_1,\boldsymbol p_1)$, where $N_l$ refers to the number of balls in 100 cells, and the cell probabilities are specified by $\boldsymbol p_l\in\R^{100}$. Specifically, we set $N_0=N_1=100$, and  \smash{$\boldsymbol p_0=(\frac{1}{100},\dots,\frac{1}{100})$}. The vector $\boldsymbol p_1$ was obtained by perturbing $\boldsymbol p_0$ and then normalizing it. Namely, letting \smash{$\boldsymbol z\in\R^{100}$} be a vector of i.i.d. $N(0,1)$ variables, we defined the vector \smash{$\boldsymbol p_1= |\boldsymbol p_0+\ts\frac{1}{300} \boldsymbol z|/\|\boldsymbol p_0+\ts\frac{1}{300} \boldsymbol z\|_1$}, where $|\cdot|$ refers to coordinate-wise absolute value.

\subsection{Numerical results}\label{data:results}

\paragraph{Interpreting the plots}
For each dataset, we plot the ground truth value $3\sigma_t$ as a function of $t=1,\dots,1,\!000$, where  the y-axis is expressed in units of  \%, so that a value $3\sigma_t=.01$  is marked as 1\%.
Alongside each curve for $3\sigma_t$,
we plot the averages of $3\hat{\sigma}_{t,\textsc{o},\text{extrap}}$ (green) , $3\hat{\sigma}_{t,\textsc{h},\text{extrap}}$ (purple), and $3\hat{\sigma}_{t,\textsc{g},\text{extrap}}$ (orange)  over their 1,000 realizations, with error bars indicating the spread between the 10th and 90th percentiles of the estimates. Here, the error bars are only given to illustrate the variance of the estimates, conditionally on $\D$, and they are not proposed as confidence intervals for $\sigma_t$. (Indeed, our main focus is on the fluctuations of $\Err_t$, rather than the fluctuations of variance estimates.) Lastly, we plot results for the class-wise parameters $\sigma_{t,l}$ in the same manner, but in order to keep the number of plots manageable, we only display the class $l$ with the highest value of $3\sigma_{t,l}$ at $t=1,\!000$. This is reflected in the plots, since the values of $3\sigma_{t,l}$ for the chosen class $l$ are generally larger than $3\sigma_t$.

\paragraph{Walking through an example (Figure 2)} To explain the plots from the user's perspective, suppose the user trains an initial ensemble of $t_0=200$ classifiers with the `census income' data. (The following considerations will apply in the same way to the other datasets in Figures 3-7.)  At this stage, the user may  compute either of the estimators $\hat{\sigma}_{200,\textsc{o}}$ or  $\hat{\sigma}_{200,\textsc{h}}$. In turn, the user may follow the definitions~\eqref{extrapests} to plot the extrapolated estimators for all $t\geq 200$ at no additional cost. These curves will look like the purple or green curves in the left panel of Figure 2, up to a small amount of variation indicated by the error bars.

 If the user wants to select $t$ so that $3\sigma_t$ is at most, say 0.5\%, then the purple or green curves in the left panel of Figure~2 would tell the user that 200 classifiers are already sufficient, and no extra classifiers are needed  (which is correct in this particular example). Alternatively, if the user happens to be interested in the class-wise error rate for $l=1$, and if the user wants $3\sigma_{t,1}$ to be at most 0.5\%, then the curve for the \textsc{oob} estimator accurately predicts that approximately 600 total (i.e.~400 extra) classifiers are needed. By contrast, the hold-out method is conservative, and indicates that approximately 1,000 total (i.e.~800 extra) classifiers should be trained. So, in other words, the hold-out estimator would still provide the user with the desired outcome, but at a higher computational cost.


\vspace{0.3cm}
\begin{figure*}[h!]
\centering
{\includegraphics[angle=0,
  width=.45\linewidth,height=.45\linewidth]{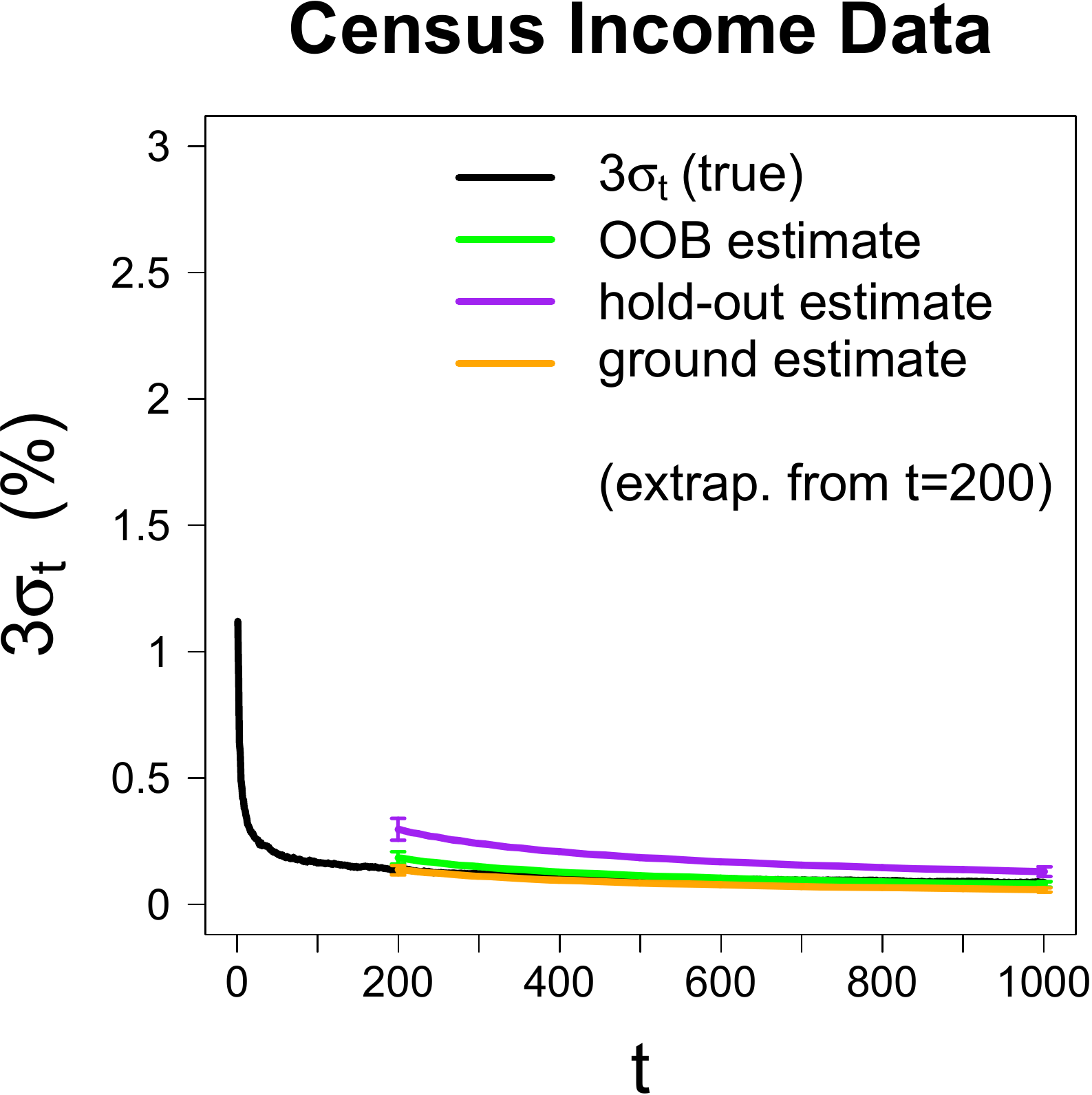}} \ \ \ 
{\includegraphics[angle=0,
  width=.47\linewidth,height=.47\linewidth]{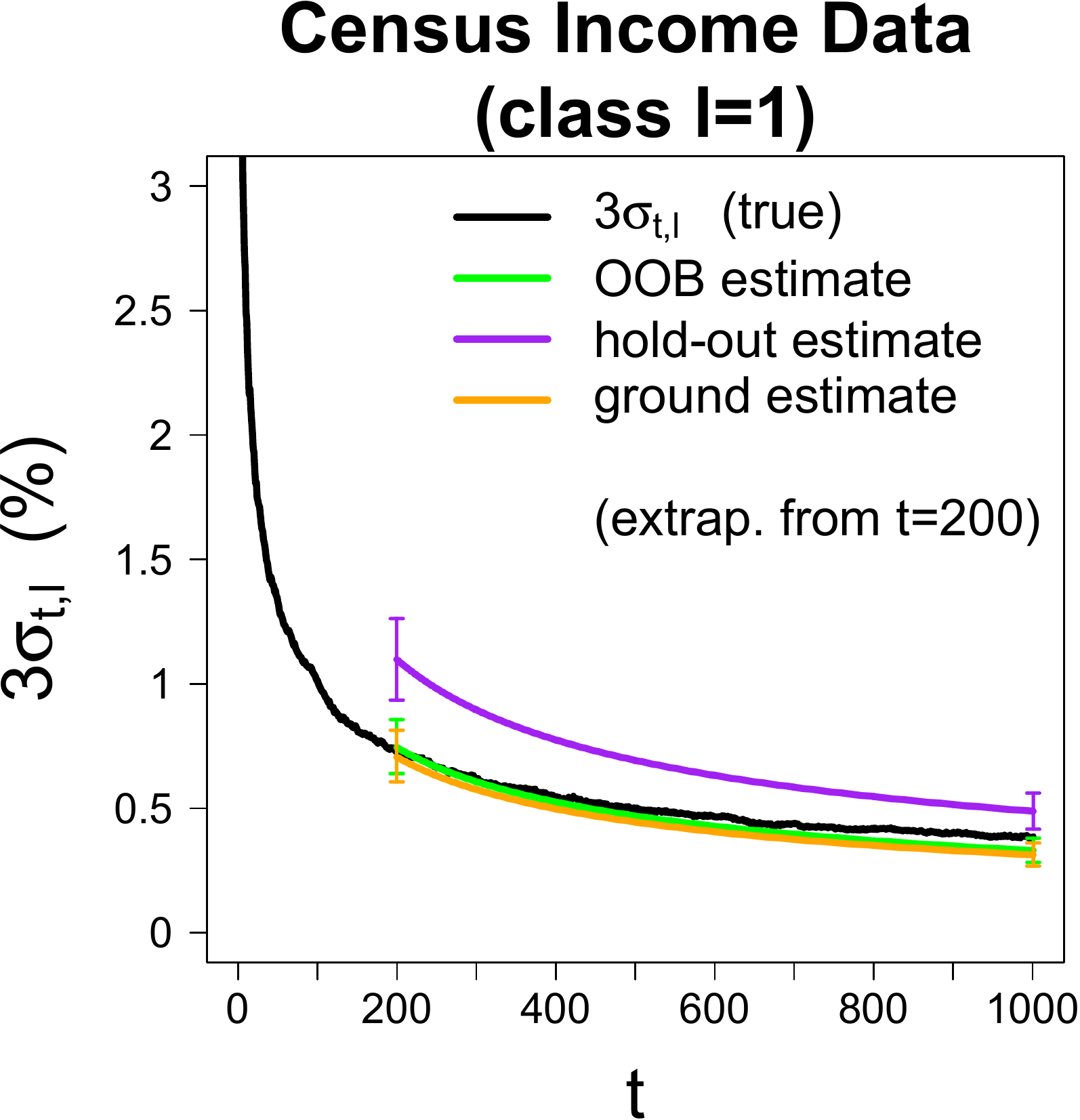}}

\caption{Results for census income data.}
\label{fig:abalone}
\end{figure*}

~\\
\vspace{.2cm}

\begin{figure*}[h!]
\centering
{\includegraphics[angle=0,
   width=.45\linewidth,height=.45\linewidth]{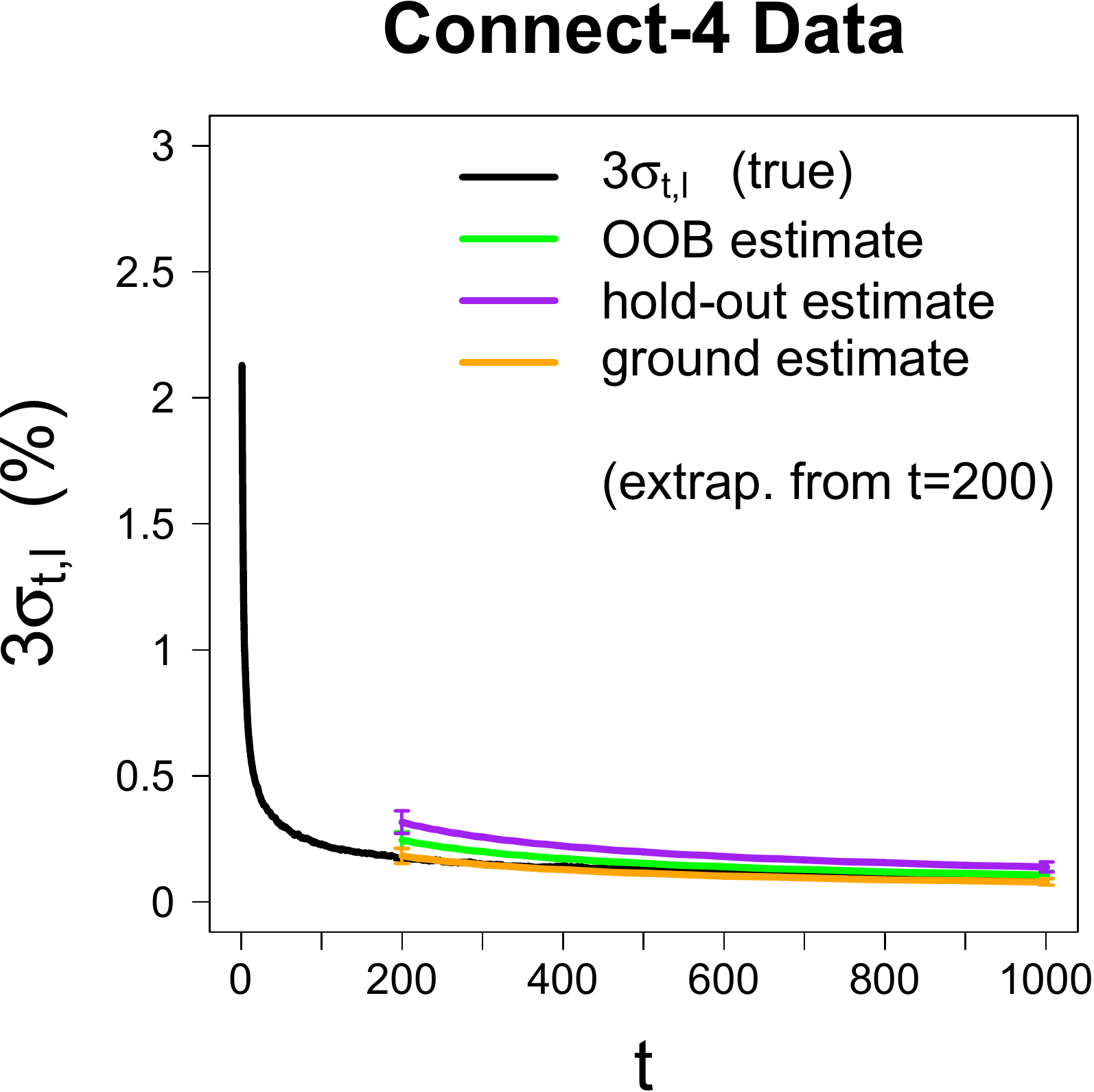}}  \ \ \ 
{\includegraphics[angle=0,
  width=.47\linewidth,height=.47\linewidth]{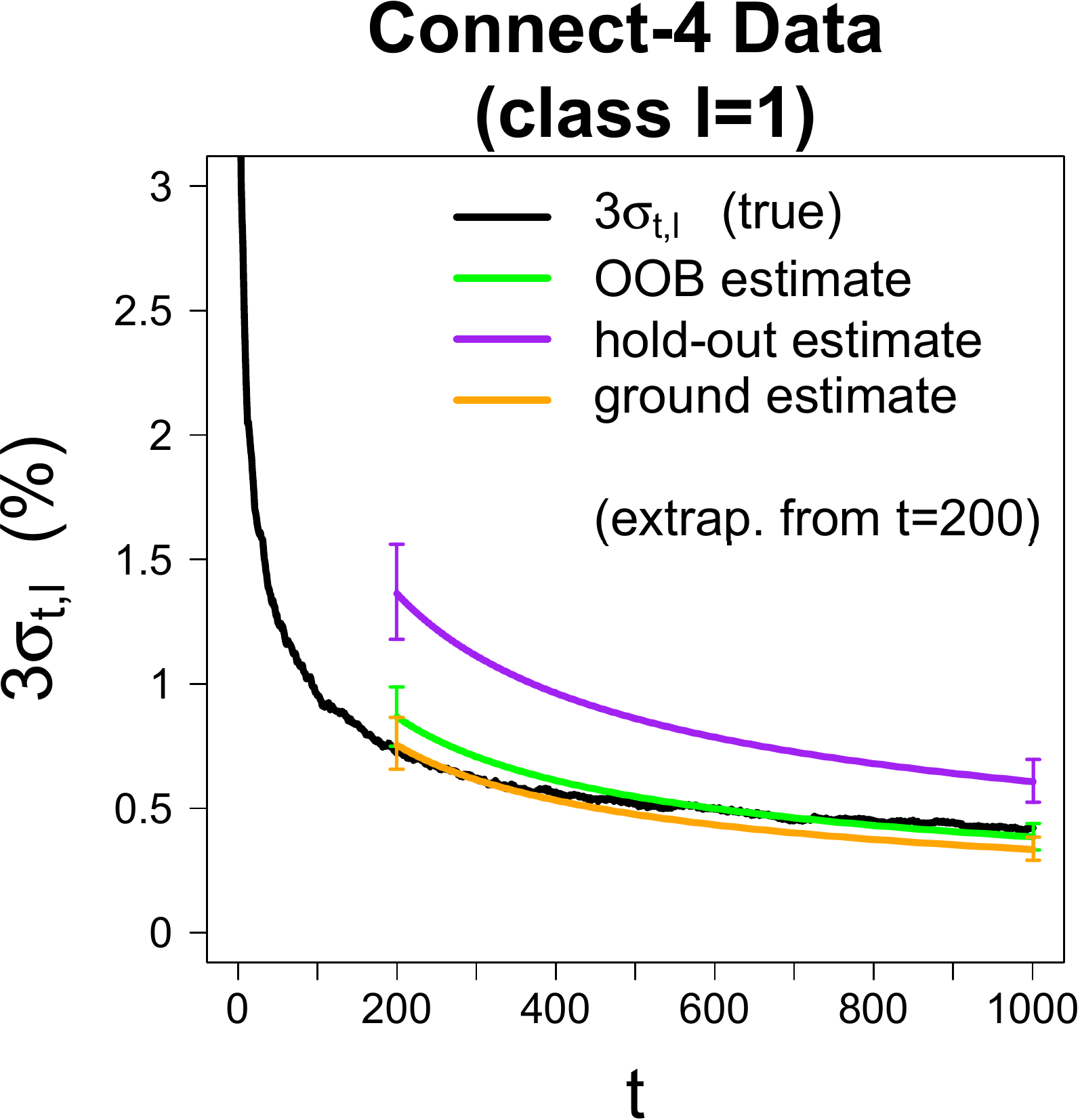}}

\caption{Results for connect-4 data.}
\label{fig:contraceptive}
\end{figure*}


%
%
%

%
%

\begin{figure*}[h!]
\centering
  \  {\includegraphics[angle=0,
  width=.48\linewidth,height=.48\linewidth]{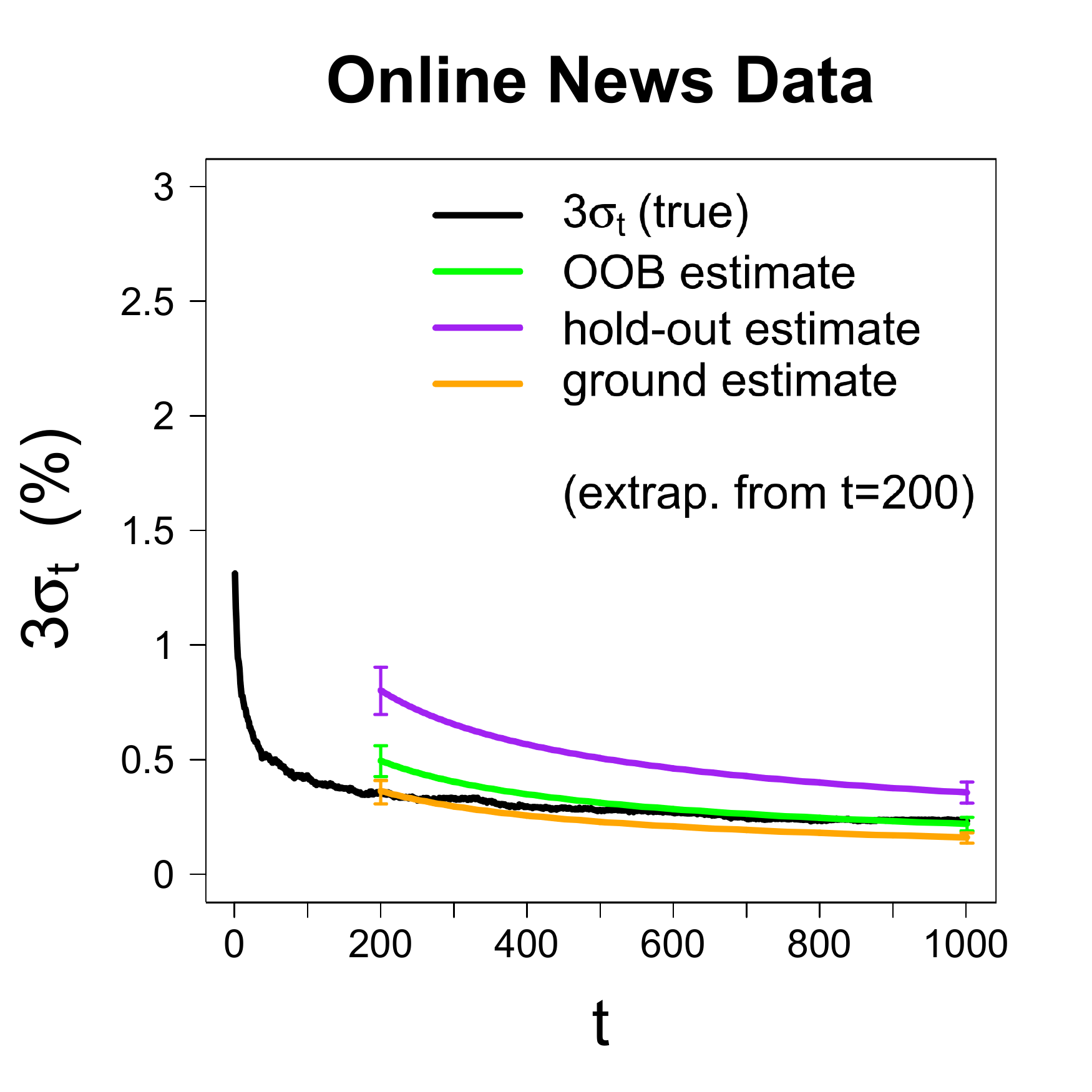}} \hspace{-.2cm} 
{\includegraphics[angle=0,
  width=.49\linewidth,height=.49\linewidth]{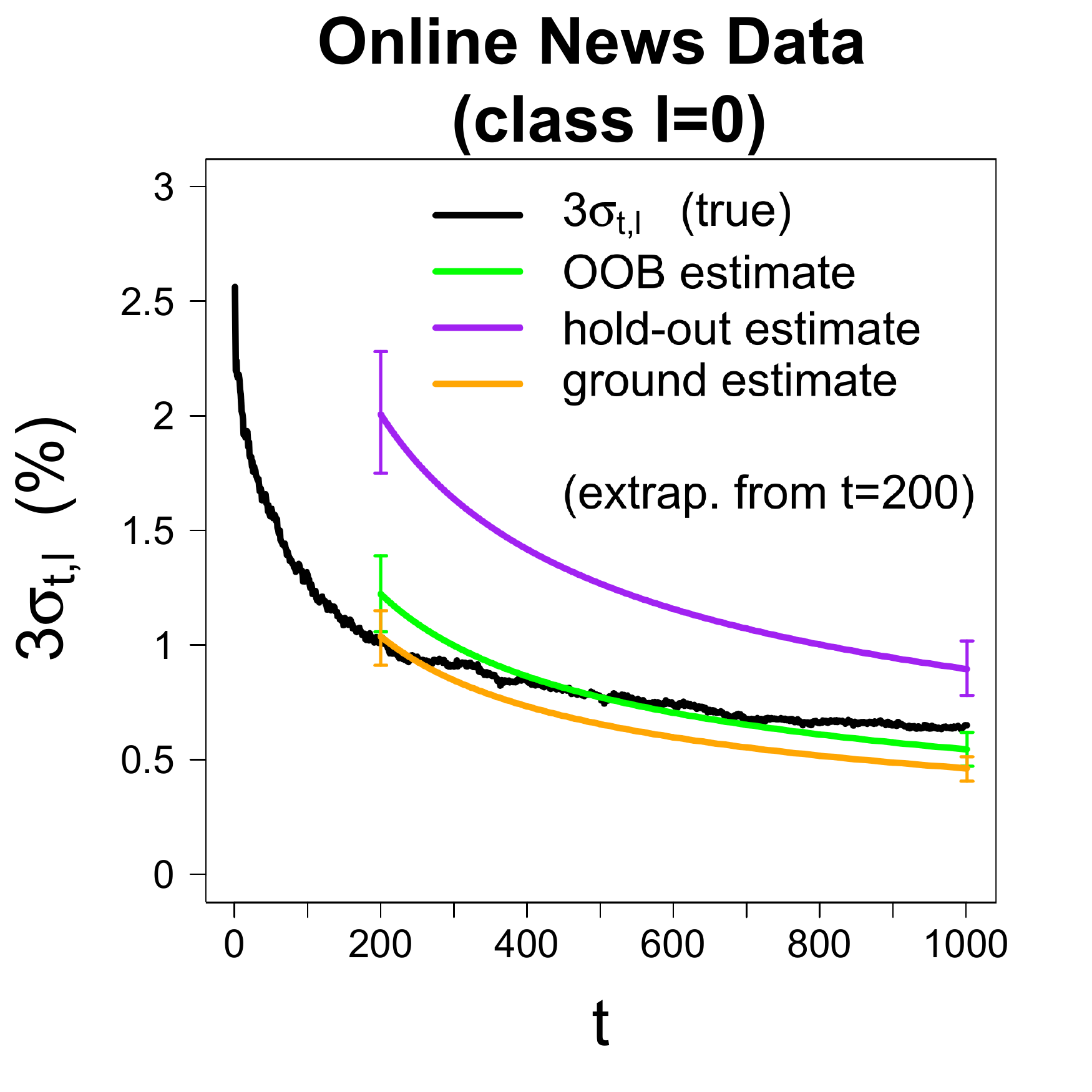}}
 \vspace{-.3cm} 
\caption{Results for online news data.}
\label{fig:contraceptive}
\end{figure*}

~\\
\begin{figure*}[h]
\centering
{\includegraphics[angle=0,
  width=.45\linewidth,height=.45\linewidth]{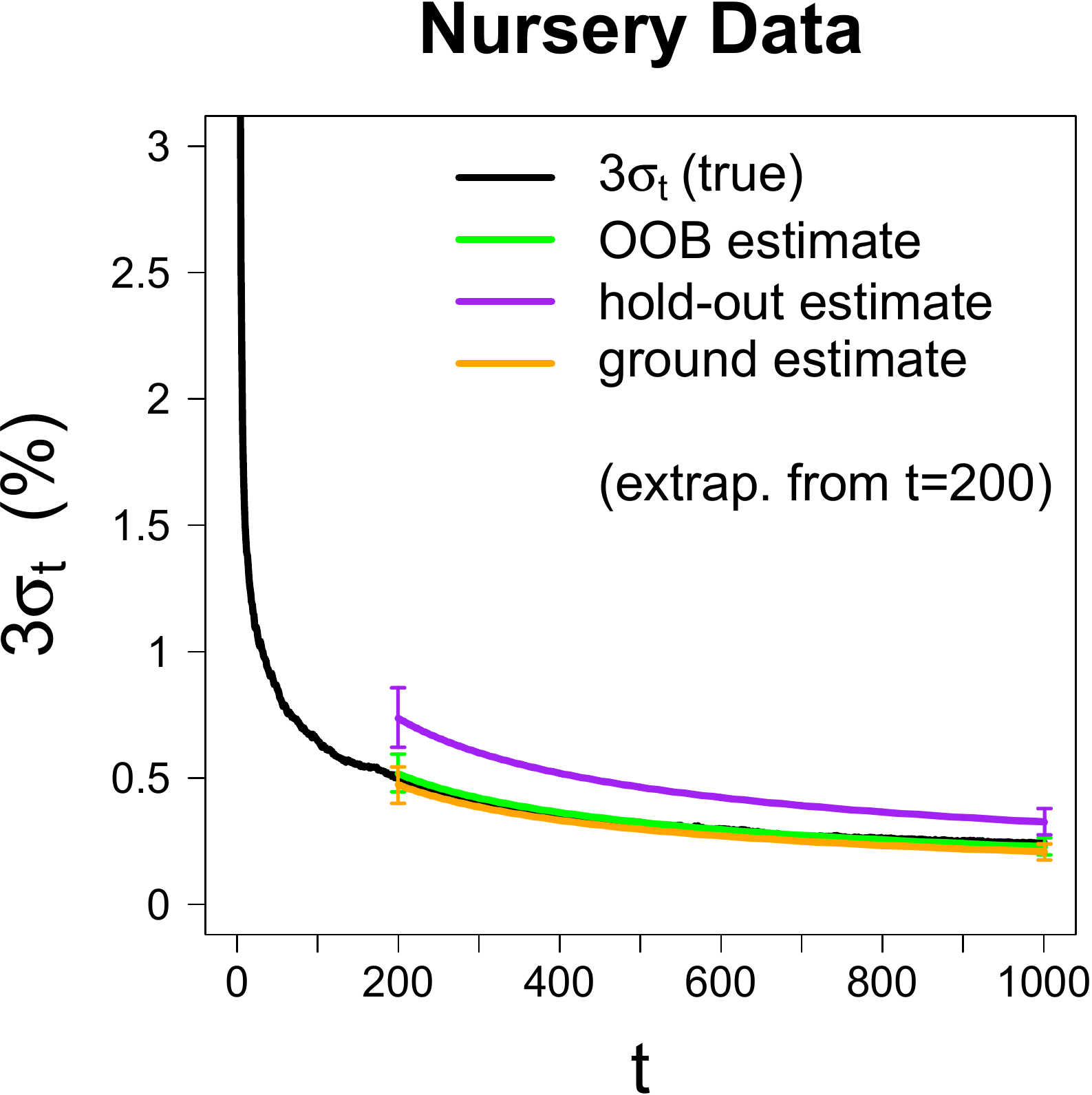}} \ \ \ 
{\includegraphics[angle=0,
  width=.47\linewidth,height=.47\linewidth]{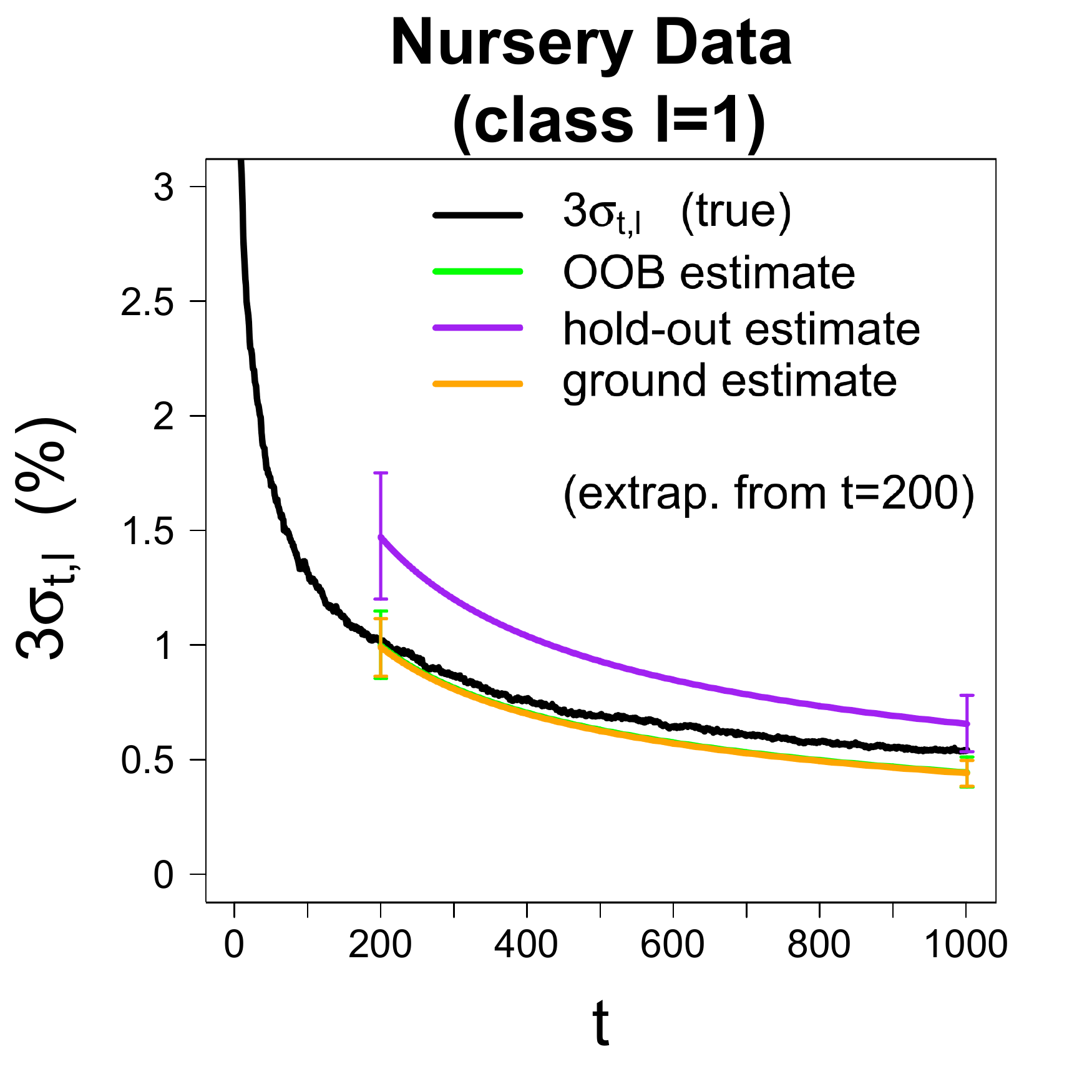}}

\caption{Results for nursery data.}
\label{fig:contraceptive}
\end{figure*}


\begin{figure*}[h!]
\centering
{\includegraphics[angle=0,
  width=.45\linewidth,height=0.45\linewidth]{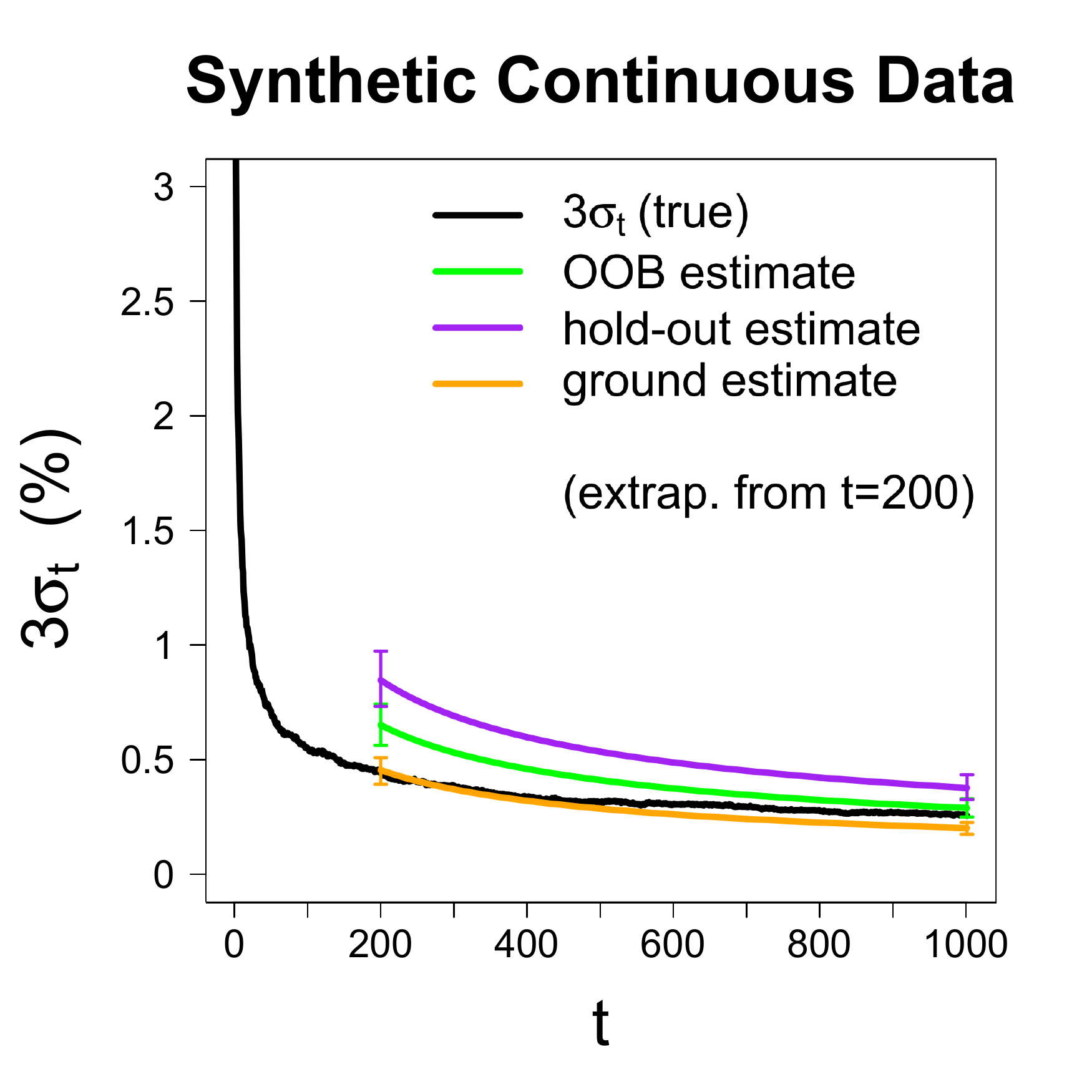}} \ 
{\includegraphics[angle=0,
  width=.45\linewidth,height=.45\linewidth]{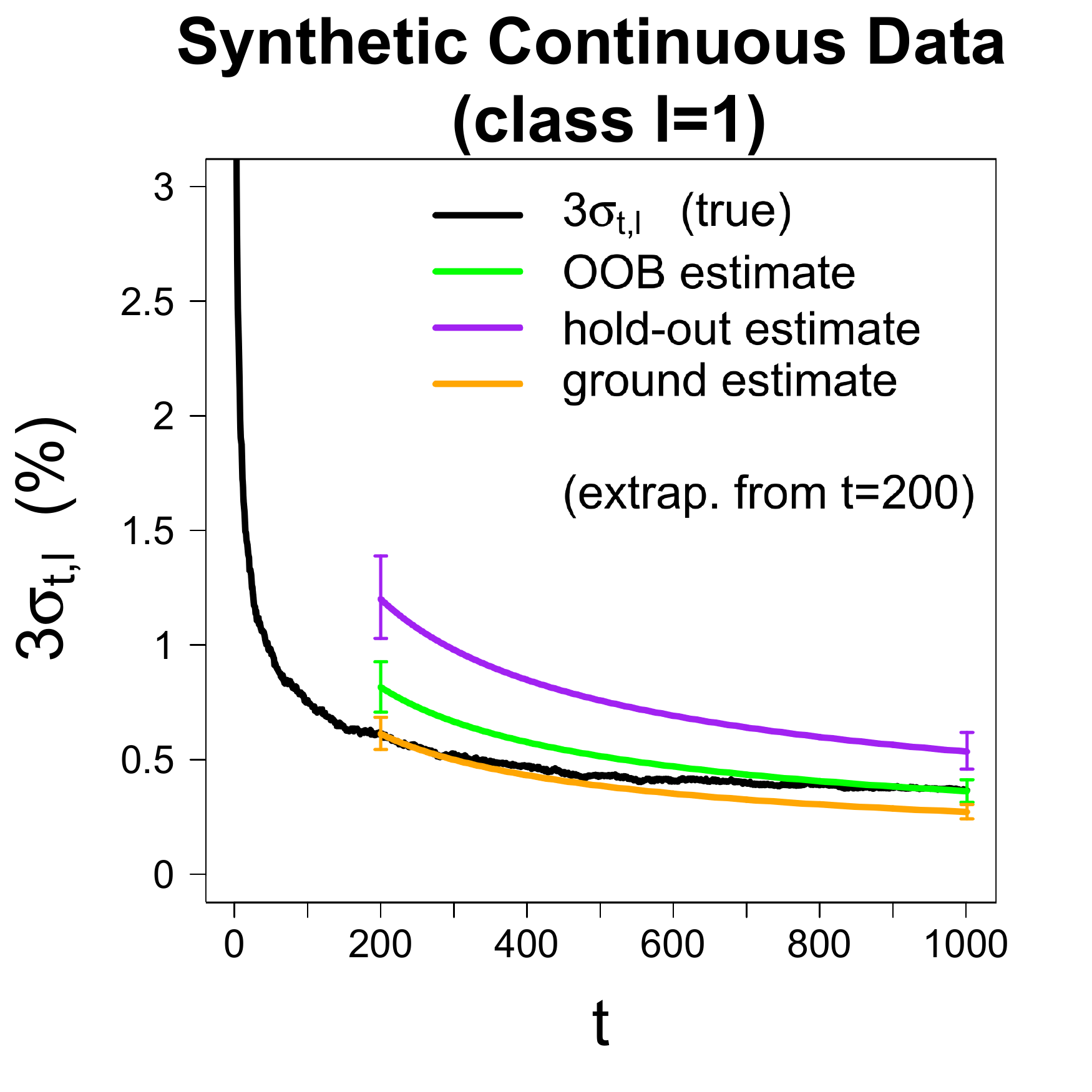}}

\caption{Results for synthetic continuous data.}
\label{fig:contraceptive}
\end{figure*}


\begin{figure*}[h!]
\centering
{\includegraphics[angle=0,
  width=.47\linewidth]{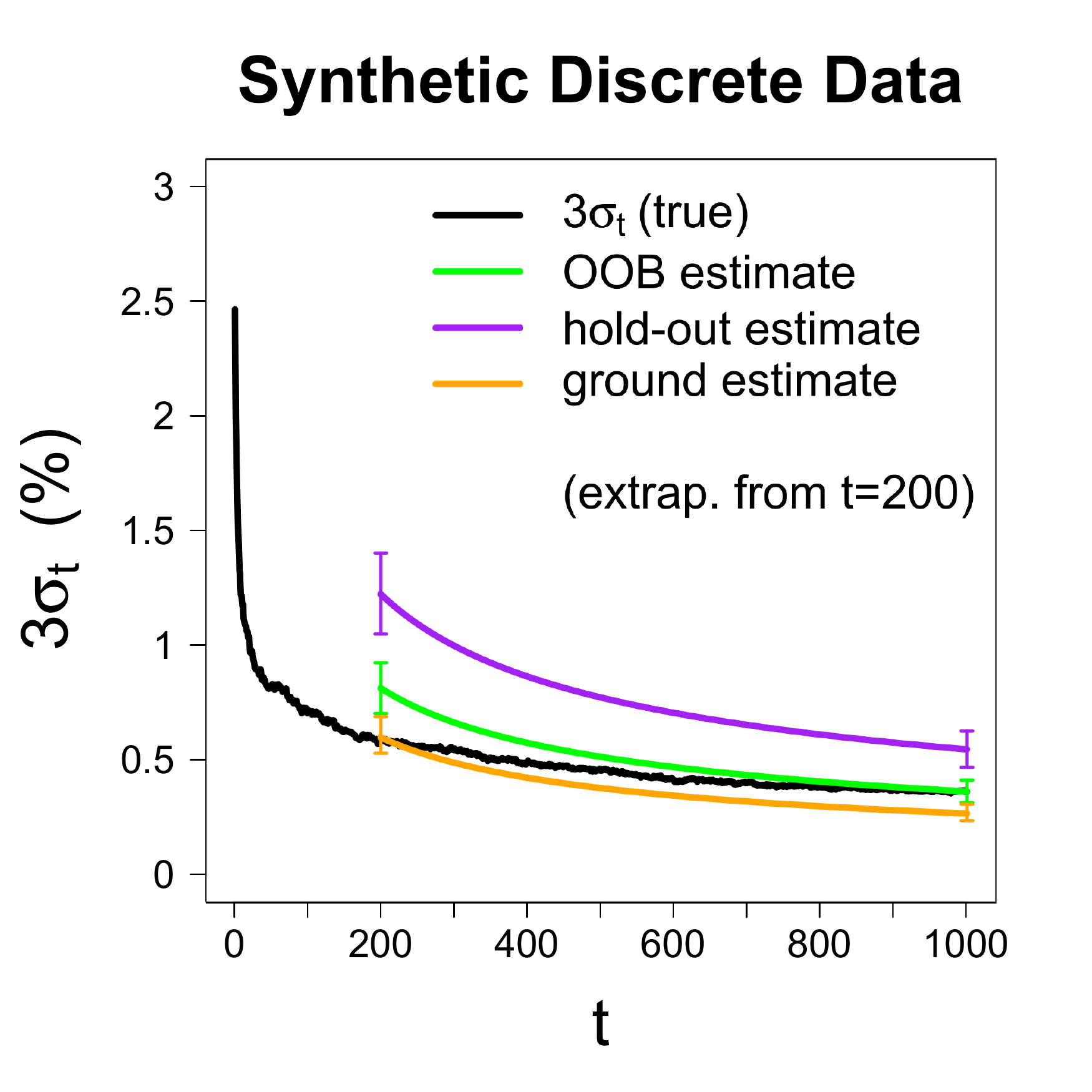}}
{\includegraphics[angle=0,
  width=.47\linewidth]{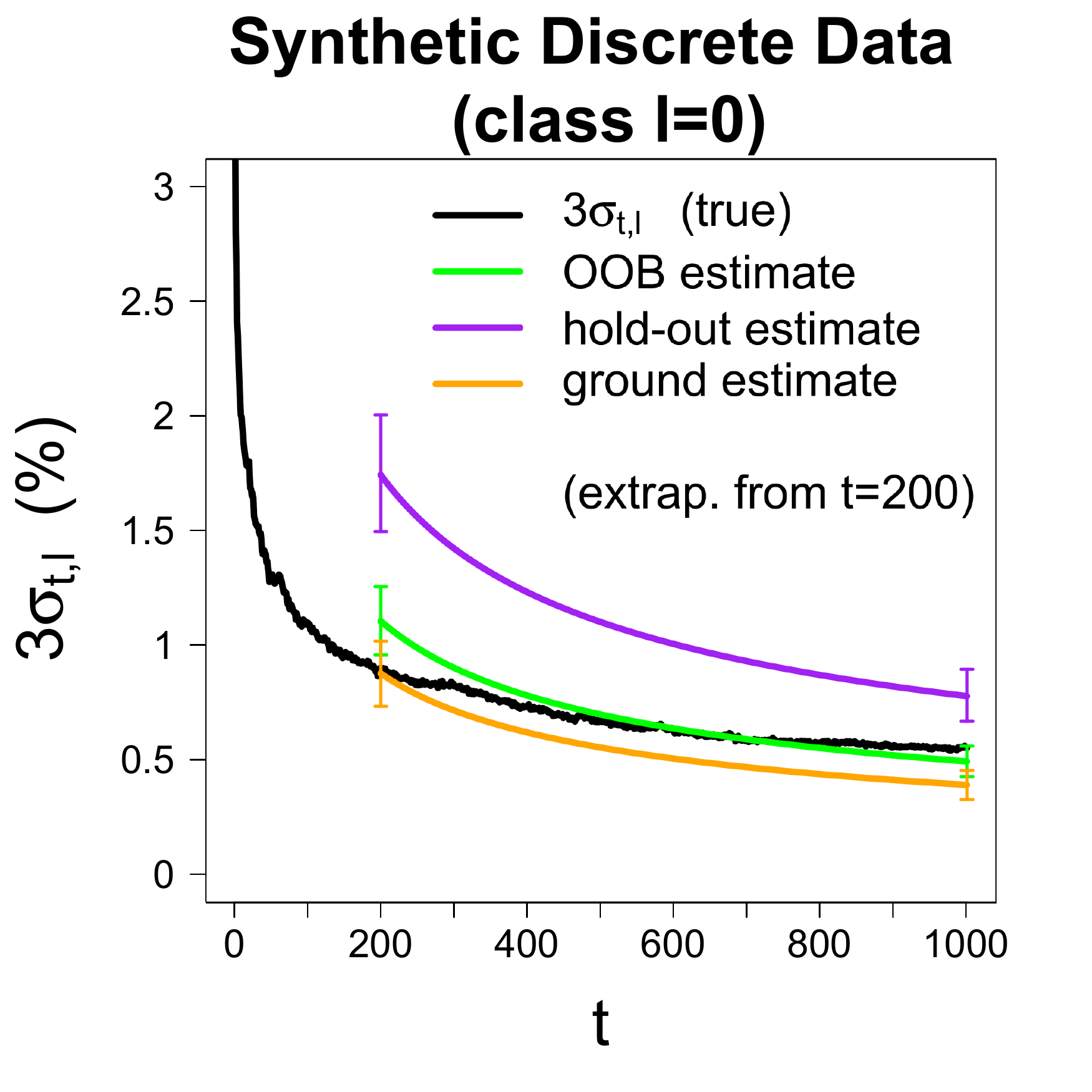}}

\caption{Results for synthetic discrete data.}
\label{fig:contraceptive}
\end{figure*}

\paragraph{Comments on numerical results} Considering all of the datasets collectively,  the plots show that the extrapolated \textsc{oob} and ground estimators are generally quite accurate. Meanwhile, the hold-out estimator tends to be conservative, due to an upward bias. Consequently, the \textsc{oob} method should be viewed as preferable, since it is both more accurate, and does not require data to be held out. Nevertheless, when considering the hold-out estimator, it is worth noting that the effect of the bias actually diminishes with extrapolation, and even if the initial value $3\hat{\sigma}_{t_0,\textsc{h},\text{extrap}}$  has noticeable bias at $t_0=200$, it is possible for the extrapolated value $3\hat{\sigma}_{t,\textsc{h},\text{extrap}}$ to have relatively small bias at $t=1,\!000$.

To understand where the bias of the hold-out estimator comes from, imagine that two ensembles have been trained on the same data, and suppose their accuracy is compared on a small hold-out set. In this situation, it is possible for their observed error rates on the hold-out set to noticeably differ --- even if the true error rates are very close. For this reason, the small size of $\D_{\text{hold}}$ leads to greater variation among the estimated values $\widehat{\Err}(\tilde{A}^*)$ generated in the hold-out version of Algorithm 2, which leads to an inflated estimate of $\sigma_t$. By contrast, the ground estimator suffers less from this bias because it relies on the much larger ground truth set $\D_{\text{ground}}$ in place of $\D_{\text{hold}}$. Similarly, the  \textsc{oob} estimate of $\sigma_t$ is less susceptible to this bias, because it will typically use every point in the larger training set $\D$ as an ``effective test point'', preventing excess variation among the values $\widehat{\Err}(A^*)$. (Even for small choices of $t_0$,  all training points are likely to be an \textsc{oob} point for at least one classifier.)

One last point to mention is that many of the datasets have discrete features, which may violate the theoretical conditions in Assumption~\ref{ASSUMPTION}. Nevertheless, the presence or absence of discrete features does not seem to substantially affect on the performance of the estimators. So, to this extent, the bootstrap does not seem to depend too heavily on Assumption~\ref{ASSUMPTION}. (See Appendix~\ref{app:empirical} for further empirical assessment of that assumption.)

\section{Conclusion}\label{sec:conc}

We have studied the notion of algorithmic variance $\sigma_t^2=\var(\Err_t|\D)$ as a criterion for deciding when a randomized ensemble will perform nearly as well as an infinite one (trained on the same data).  To estimate this parameter, we have developed a new bootstrap method, which allows the user to directly measure the convergence of randomized ensembles with a guarantee that has not previously been available.

With regard to practical considerations, we have shown that our bootstrap method can be enhanced in two ways. First, the use of a hold-out set can be avoided with the  \textsc{oob} version of Algorithm 2, and our numerical results show that the \textsc{oob} version is preferable when hold-out points are scarce. Second, the extrapolation technique substantially reduces the cost of bootstrapping. Furthermore,
we have analyzed the cost of the method in terms of floating point operations to show that it compares favorably with the cost of training a single ensemble via random forests.

From a theoretical standpoint, we have analyzed the proposed method within the framework of a first-order model for randomized ensembles. In particular, for a generic ensemble whose classifiers are conditionally i.i.d.~given $\D$, there is a corresponding first-order model that matches the generic ensemble with respect to its average error rate $\E[\Err_t|\D]$. Under this setup, our main result shows that the proposed method consistently approximates $\mathcal{L}(\sqrt{t}(\Err_t-\err_{\infty})|\D)$ as $t\to\infty$.  Some extensions of this result could include generalizations to other models of randomized ensembles (e.g., along the lines of~\citet{biau2008,biau2012,scornetconsistency,scornetIEEE}), as well as corresponding results in the context of regression ensembles, which we hope to pursue in future work. More generally, due to the versatility of bootstrap methods, our approach may also be relevant to measuring the convergence of other randomized algorithms, as in~\cite{lopes_matrix,lopes_LS}.

\section*{Acknowledgements} The author thanks Peter Bickel,  Philip Kegelmeyer, and Debashis Paul for helpful discussions. In addition, the author thanks the editors and referees for their valuable feedback, which significantly improved the paper.

\bibliographystyle{imsart-nameyear}

\bibliography{rfBoot_classif_14aug2016.bib,MajVote2016bib.bib}

~\newpage

\begin{center}

{\bf{\large{Supplementary Material
}}}

\end{center}

\section*{Outline of appendices and the proof of Theorem~\ref{THM:CLT}} Here we explain how the the main components of the proof of Theorem~\ref{THM:CLT} fit together. First, in Appendix~\ref{sec:lift}, we show that under a first-order model, there is an explicit functional $\phi$ such that
\begin{equation}\label{eq:basicrelation}\tag{1}
\Err_t-\err_{\infty} = \phi(\F_t)-\phi(\textup{id}),
\end{equation}
where   $\F_t(u)$ is the empirical c.d.f. of independent~Uniform[0,1] random variables $\textup{\bf{U}}_t=(U_1,\dots,U_t)$, and $\id$ is the identity function on $[0,1]$. Also note that under the first-order model, the set of variables $\textup{\bf{U}}_t$ plays the role of the randomizing parameters $\boldsymbol \xi_t$. It follows that bootstrapping the classifier functions $Q_1,\dots,Q_t$ corresponds to bootstrapping the variables $\textup{\bf{U}}_t$. More specifically, if we let\begin{equation*}
\F_t^*(u) :=\ts\frac{1}{t}\sum_{i=1}^t 1\{U_i^*\leq u\},
\end{equation*}
where $U_1^*,\dots,U_t^*$ are drawn with replacement from  $\textup{\bf{U}}_t$, then the bootstrap counterpart of~\eqref{eq:basicrelation} is given by
\begin{equation}\label{eq:bootbasicrelation}\tag{2}
 \Err_t^*-\Err_t= \phi(\F_t^*)-\phi(\F_t).
\end{equation}
Due to the relations~\eqref{eq:basicrelation} and~\eqref{eq:bootbasicrelation}, the remainder of the proof deals with converting ``bootstrap consistency for $\F_t$'' into ``bootstrap consistency for $\Err_{t}$'', and this is handled with the functional delta method. In this way, the proof boils down to establishing the Hadamard differentiability of $\phi$, which is done in Theorem~\ref{THM:HADAMARD} of Appendix~\ref{sec:hadamard}, and this should be viewed as the core technical result of the paper. In turn, the details of applying the functional delta method, as well as the conclusion of the proof, are given in Appendix~\ref{sec:delta}.

The remainder of the appendices are organized as follows. Appendices~\ref{sec:highlevel},~\ref{app:variation}, and~\ref{app:smooth} contain the arguments for proving Theorem~\ref{THM:HADAMARD} on Hadamard differentiability. Throughout these arguments, we will refer to various technical lemmas that are stated and proved in Appendix~\ref{lemmas}. Also, we henceforth assume that the first-order model holds, so that $Q_i=T_i$ for $i\geq 1$, and the sets $\textup{\bf{U}}_t$ and $\boldsymbol \xi_t$ are synonymous. Lastly, Appendix~\ref{app:assumption} discusses the theoretical and empirical assessment of Assumption~\ref{ASSUMPTION}.

\paragraph{Notation}  If $\mathcal{\mathcal{A}}$ is a generic subset of Euclidean space, we write $\mathcal{A}^{\circ}$ for the interior and $\partial \mathcal{A}$ for the boundary. The identity function on $\mathcal{A}$ is denoted $\text{id}_{\mathcal{A}}$, or just $\id$ when there is no ambiguity. The spaces $C[0,1]$, $D[0,1]$, $\mathcal{F}[0,1]$,  and $\ell^{\infty}[0,1]$ respectively denote the following sets of  real-valued functions on [0,1], all equipped with the supremum norm:  continuous functions, c\`adl\`ag functions, distribution functions with no point mass at 0, and bounded functions. Hence, if we write $h_t\to h$ as $t\to\infty$ for functions $h_t$ and $h$ in these spaces, it is implicit that the convergence is uniform. (The analogous definitions apply when the open interval $(0,1)$ replaces $[0,1]$.) Lastly, we use $\ell^{\infty}[\mathcal{A};\R^{k-1}]$ to denote the space of coordinate-wise bounded functions from $\mathcal{A}$ to $\R^{k-1}$, with norm
\mbox{$\|g\|=\max_{1\leq l\leq k-1}\sup_{\theta\in\mathcal{A}}|g_l(\theta)|$}, where $g=(g_1,\dots,g_{k-1})$

\appendix

\section{Proof of Theorem~{\ref{THM:CLT}}}\label{app:mainthm}

\subsection{Representing $\Err_{t}$ as a functional of the empirical process}\label{sec:lift}

Working under the first-order model, our aim in this subsection is to construct an explicit functional $\phi$ such that
\begin{equation*}
\Err_t = \phi(\F_t).
\end{equation*}
 (Later on, we will formally define $\err_{\infty}:=\phi(\id)$, which will turn out to be the limit of $\Err_t$ as $t\to\infty$.)
To construct this functional,  consider the basic relation $\Err_t=\sum_{l=0}^{k-1} \pi_l \Err_{t,l}$, where $\pi_l=\P(Y=\boldsymbol e_l)$ is the  $l$th class proportion, and $\Err_{t,l}$ is defined in equation~\eqref{classerrdef}. Hence, for each $l=0,\dots,k-1$, it is enough to construct a functional $\phi_l$ such that
\begin{equation}\label{eqn:functionalrelation}
\Err_{t,l}=\phi_l(\F_t),
\end{equation}
which implies that $\phi$ may be written as $\phi=\ts\sum_{l=0}^{k-1} \pi_l \phi_l$.

\paragraph{The lifting operator $\boldsymbol L$} As a first step in deriving the relation~\eqref{eqn:functionalrelation}, we will show that the stochastic process $\bar{T}_t:\mathcal{X}\to \Delta$ can be obtained from a linear operator that ``lifts'' the univariate function  $\F_t$ on $[0,1]$ to a multivariate  map on the simplex $\Delta$.
The relationship between $\bar{T}_t$ and $\F_t$ is unraveled by noting that for any $l\geq 2$, the definition of the interval $I_l$ gives
\begin{equation}\label{unravel}
\begin{split}
[\bar{T}_t(x)]_l
&=\ts\frac{1}{t}\sum_{i=1}^t \ [T_i(x)]_l\\[0.4cm]
&= \ts\frac{1}{t}\sum_{i=1}^t 1\big \{ U_i \in I_l(\vartheta(x))\big\},\\[0.4cm]
&= \F_t(\vartheta_1(x)+\cdots+\vartheta_l(x))-\F_t(\vartheta_1(x)+\cdots+\vartheta_{l-1}(x)),
\end{split}
\end{equation}
and also $[\bar T_t(x)]_1=\F_t(\vartheta_1(x))$.
Due to the formula above, it is natural to consider an operator, denoted $\boldsymbol L$, that lifts a generic function $h:\R\to \R$ to a new function $\boldsymbol L(h):\R^{k-1}\to \R^{k-1}$ whose $l$th coordinate is given by
\begin{equation}\label{liftdef}
[ \boldsymbol L(h)(\theta)]_l := h(\theta_1+\cdots+\theta_l) - h(\theta_1+\cdots+\theta_{l-1}), \text{ \ \ \ \ for any } \theta\in\R^{k-1},
\end{equation}
where $h(\theta_1+\dots+\theta_{l-1})$ is understood as 0 when $l=1$.
 In particular, the definition of $\boldsymbol L$, and the calculation in equation~\eqref{unravel}, show that $\bar{T}_t(x)$ can be expressed by lifting $\mathbb{F}_t$ and then evaluating at $\vartheta(x)$, 
\begin{equation}\label{TLrep}
\bar{T}_t(x)=\boldsymbol L(\F_t)(\vartheta(x)).
\end{equation}

\noindent Since many of our arguments will rely on special properties of $\boldsymbol L$, we briefly summarize these properties below.\\

\paragraph{Properties of the lifting operator}
If we let $f,g,$ and $F$ denote generic real-valued functions on $\R$, and let $a$ denote a scalar, then the following conditions hold.
\begin{itemize}
\small
\item Linearity.
$$\boldsymbol L(ag+h) = a\boldsymbol L(g)+\boldsymbol L(h).$$
\item Composition.
$$\boldsymbol L(g\circ h) = \boldsymbol L(g)\circ \boldsymbol L(h).$$
\item Identities.
$$\boldsymbol L(\text{id}_{[0,1]}) = \text{id}_{\Delta}.$$

\item Simplices.
$$\text{If $F$ is non-decreasing and maps $[0,1]$ to $[0,1]$, then $\boldsymbol L(F)$ maps $\Delta$ to $\Delta$.} $$
\item Inverses.
$$\text{If $F$ is invertible, then $\boldsymbol L(F)$ is invertible and $(\boldsymbol L(F))^{-1} = \boldsymbol L(F^{-1})$.}$$
\end{itemize}
The fact that $\boldsymbol L$ respects composition of functions is the only property that takes some care to verify, but we omit the calculation for brevity. 
In turn, the ``Inverses'' property follows by combining the ``Composition'' and ``Identities'' properties.

\paragraph{Deriving the relation $\Err_{t,l}=\phi_l(\F_t)$}
To finish the derivation of $\Err_t$ as a functional of $\F_t$, it is helpful to consider a certain subset of $\Delta$. Namely, observe that if a test point $X$ is drawn from class $l\in\{0,\dots,k-1\}$, then an error occurs for the plurality vote if and only if
\begin{equation}\label{Slfact}
\bar{T}_t(X)\in \mathcal{S}_l,
\end{equation}
where the set $\mathcal{S}_l\subset\Delta$ is defined as 
\begin{equation*}
\mathcal{S}_l:=\Big\{\theta \in\Delta \, \Big| \, \theta_l\leq  \theta_{l'} \text{ for some } l'\in \{0,\dots,k-1\}\setminus\{l\}\Big\},
\end{equation*}
where $\theta_0:=1-(\theta_1+\cdots+\theta_{k-1})$. Consequently, we have
\begin{equation*}
\Err_{t,l} = \P(\bar{T}_t(X) \in \mathcal{S}_l\big\bracevert \D, \textup{\bf{U}}_t, Y=\boldsymbol e_l).
\end{equation*}
 Next, the formula $\bar{T}_t(x)=\boldsymbol L(\F_t)(\vartheta(x))$ from equation~\eqref{TLrep} implies
$$\bar{T}_t(X)\in \mathcal{S}_l \ \ \ \text{ if and only if }  \ \ \ \  \vartheta(X)\in [\boldsymbol L(\F_t)]^{-1}(\mathcal{S}_l),$$
where the set $[\boldsymbol L(\F_t)]^{-1}(\mathcal{S}_l)$ is the pre-image of $\mathcal{S}_l$ under the map $\boldsymbol L(\F_t)$.
Finally, if we recall the definition of $\mu_l=\mathcal{L}(\vartheta(X)|\D,Y=\boldsymbol e_l)$, then $\Err_{t,l}$  is the $\mu_l$-probability of the pre-image,
\begin{equation}\label{fdelta0}
\Err_{t,l}  =\mu_l\Big([\boldsymbol L(\F_t)]^{-1}(\mathcal{S}_l)\Big),
\end{equation}
which provides us with the desired representation of $\Err_{t,l}$ as a functional of $\F_t$. More precisely, we obtain the relation \smash{$\Err_{t,l}=\phi_l(\F_t)$} by defining the   functional $\phi_l:\mathcal{F}[0,1]\to\R$ according to
\begin{equation}\label{phildef}
\phi_l(G) := \mu_l\Big([\boldsymbol L(G)]^{-1}(\mathcal{S}_l)\Big),  \ \ \ \ \text{ for any } \ \ \ G\in \mathcal{F}[0,1],
\end{equation}
where $l=0,\dots,k-1$, and we note that $\F_t$ lies in $\mathcal{F}[0,1]$ almost surely.

\subsection{Hadamard differentiability of $\phi_l$}\label{sec:hadamard}

Since the random distribution function $\F_t$ approaches $\id_{[0,1]}$ as $t\to\infty$, the asymptotic behavior of $\Err_{t}$ is closely linked to the smoothness of $\phi=\sum_{l=0}^{k-1}\pi_l \phi_l$ at the ``point'' $\id_{[0,1]}$. Specifically, we will use the standard notion of  Hadamard differentiability, as reviewed below~\citep{vaartWellner}. 
\begin{definition}[Hadamard differentiability\footnote{The added generality of letting $\textup{\bf{B}}$ be a normed space, rather than just $\R$, will be needed in Appendix~\ref{lemmas}.}] \emph{Let $\textup{\bf{B}}$ be a normed space. A map $\psi:\mathcal{F}[0,1]\to \textup{\bf{B}}$ is \emph{Hadamard differentiable} at $G_0\in\mathcal{F}[0,1]$ tangentially to $C[0,1]$, if there is a continuous linear map $\psi_{G_0}':C[0,1]\to \textup{\bf{B}}$
 such that as $t\to\infty$,
\begin{equation*}
\frac{\psi(G_0+\ve_t h_{t})-\psi(G_0)}{\ve_t} \xrightarrow{ \ \ \ } \psi_{G_0}'(h),
\end{equation*}
for all converging sequences of positive numbers $\ve_t\to 0$ and functions $h_t\to h$, such that $G_0+\ve_th_t\in \mathcal{F}[0,1]$ for every $t\geq 1$, and $h\in C[0,1]$. In particular, the linear map $\psi'_{G_0}$ is referred to as the Hadamard derivative of $\psi$ at $G_0$.
}
\end{definition}

The following result (Theorem~\ref{THM:HADAMARD}) is the core technical result of the paper. It provides conditions under which the functionals $\phi_l$ are Hadamard differentiable at $G_0=\id_{[0,1]}$, tangentially to $C[0,1]$. 
Likewise, the theorem implies that the same differentiability property holds for the linear combination $\phi=\sum_{j=1}^{k-1}\pi_l \phi_l$, for which  $\Err_t=\phi(\F_t)$.

Although each functional $\phi_l$ is defined in terms of the particular set $\mathcal{S}_l$ and the particular measure $\mu_l$, the Hadamard differentiability of $\phi_l$ is only mildly dependent on their structure.  For each measure $\mu_l$, the only property we need is  that it satisfies Assumption~\ref{ASSUMPTION}. Regarding the set $\mathcal{S}_l$, it is simple to check that its complement in $\Delta$ is a convex set with non-empty interior. In other words, the functional $1-\phi_l(G)$ may be written as $\mu_l\big([\boldsymbol L(G)]^{-1}(\mathcal{S})\big)$ for some convex set $\mathcal{S}\subset\Delta$ with non-empty interior. So, given that the Hadamard derivative of $1-\phi_l$ is that same as that of $\phi_l$, up to a sign, we state the result in terms of a generic functional $\psi$ that arises from such a set $\mathcal{S}$, and such a measure $\mu$.
\begin{thm}[Hadamard differentiability of $\phi_l$] \label{THM:HADAMARD} Let $\mathcal{S}\subset \Delta$ be a convex set with non-empty interior, and let $\mu$ be a distribution on $\Delta$
satisfying the conditions of Assumption~\ref{ASSUMPTION}.
Also let \mbox{$\psi: \mathcal{F}[0,1]\to \R$} be defined for any $G\in\mathcal{F}[0,1]$ according to
\begin{equation*}
\psi(G) := \mu\Big([\boldsymbol L(G)]^{-1}(\mathcal{S})\Big).
\end{equation*}
Then, the functional $\psi$ is Hadamard differentiable at $\textup{id}\in \mathcal{F}[0,1]$ tangentially to $C[0,1]$. 
Furthermore, the Hadamard derivative $\psi_{\textup{id}}':C[0,1]\to\R$ is given by
\begin{equation}\label{bigthm}
\psi_{\textup{id}}'(h) = \int_{\partial \mathcal{S}} \big\langle \!-\!\boldsymbol L(h)(\theta),\mathbf{n}(\theta)\big\rangle f(\theta) d\sigma(\theta),
\end{equation}
where $\mathbf{n}(\theta)$ is the outward normal to $\partial \mathcal{S}$ at the point $\theta$, 
and $d\sigma$ is Hausdorff measure on $\partial \mathcal{S}$.
\end{thm}
\paragraph{Remark} A high-level proof is given in Appendix~\ref{sec:highlevel}. In the next subsection, we apply this result in conjunction with the functional delta method to complete the proof of bootstrap consistency (Theorem~\ref{THM:CLT}).

\subsection{Functional delta method}\label{sec:delta}
Recall that the bootstrap works for $\F_t$ in the sense that the law $\mathcal{L}(\sqrt{t}(\F_t^*-\F_t)\big| \textbf{\textup{U}}_t)$ converges weakly in probability to the same limit as $\mathcal{L}(\sqrt{t}(\F_t-\id))$.
The following well-known fact shows that bootstrap consistency is preserved when $\F_t$ and $\F_t^*$ are composed with a Hadamard differentiable functional.
\begin{lemma}[Functional delta method~\citep{vaart}, Sec.~23.2.1]\label{thm:fdm}
Let $\psi:\mathcal{F}[0,1]\to\R$ be Hadamard differentiable at $\textup{\id}\in\mathcal{F}[0,1]$, tangentially to $C[0,1]$. Also, let $\mathbb{B}$ be a standard Brownian bridge on $[0,1]$. Then, as $t\to\infty$,
\begin{equation*}
\mathcal{L}\big(\sqrt{t}(\psi(\F_t)-\psi(\textup{id})\big)\xrightarrow{ \ \ w \ \ } \mathcal{L}\big( \psi'_{\textup{id}}(\mathbb{B})\big),
\end{equation*}
and
\begin{equation*}
\mathcal{L}\big(\sqrt{t}(\psi(\F_t^*)-\psi(\F_t))\big| \textup{\bf{U}}_t\big)\xrightarrow{ \ \ w \ \ } \mathcal{L}\big(\psi'_{\textup{id}}(\mathbb{B})\big), \ \text{ in } \  \P_{\textup{\bf{U}}}\text{-probability}.
\end{equation*}
\end{lemma}
\paragraph{Concluding the proof of Theorem~\ref{THM:CLT}} With this lemma in hand, Theorem~\ref{THM:CLT} on bootstrap consistency follows quickly from Theorem~\ref{THM:HADAMARD}. Specifically, if we consider $\phi=\sum_{l=0}^{k-1}\pi_l \phi_l$, with each $\phi_l$ as in equation~\eqref{phildef}, and define $\err_{\infty}:=\phi(\id)$, then the relations~\eqref{eq:basicrelation} and~\eqref{eq:bootbasicrelation} lead to
\begin{align}
\mathcal{L}\big(\sqrt{t}(\Err_t-\err_{\infty})|\D\big) &= \mathcal{L}\big(\sqrt{t}(\phi(\F_t)-\phi(\textup{id}))|\D\big),\label{top}\\[0.3cm]
\mathcal{L}\big(\sqrt{t}(\Err_t^*-\Err_t)\big| \D,\textup{\bf{U}}_t\big) &= \mathcal{L}\big(\sqrt{t}(\phi(\F_t^*)-\phi(\F_t))\big| \D,\textup{\bf{U}}_t\big),\label{bottom}
\end{align}
where we recall that $\boldsymbol \xi_t=\textup{\bf{U}}_t$ in the first-order model, as explained on p.\pageref{thetaiui} of the main text. Note also that $\phi$ implicitly depends on $\D$ through the function $\vartheta$.

It follows from Theorem~\ref{THM:HADAMARD} and Lemma~\ref{thm:fdm} that the right hand sides of~\eqref{top} and~\eqref{bottom} 
tend to the same weak limit, namely $\phi'_{\text{id}}(\mathbb{B})$. The only remaining detail to be proven in Theorem~\ref{THM:CLT} is that $\phi'_{\text{id}}(\mathbb{B})$ has a centered Gaussian distribution. This follows from the fact that $\mathbb{B}$ is a centered tight Gaussian process in $C[0,1]$, and the fact that $\phi'_{\text{id}}$ is a continuous linear functional on $C[0,1]$ \cite[Lemma~3.9.8]{vaartWellner}.


%
\section{A high-level proof of Theorem~\ref{THM:HADAMARD}}\label{sec:highlevel}
Here we give a proof of Theorem~\ref{THM:HADAMARD} that focuses on the key ideas and delegates the technical pieces to Appendices~\ref{app:variation},~\ref{app:smooth}, and~\ref{lemmas}. Consider a sequence of positive numbers $\ve_t\to 0$ and functions \smash{$h_t\to h$} such that $\id_{[0,1]}+\ve_th_t\in \mathcal{F}[0,1]$ for every $t\geq 1$, and $h\in C[0,1]$. Define the distribution function $F_t:[0,1]\to[0,1]$ by
\begin{equation*}
F_{t}:=\id_{[0,1]}+\ve_t h_{t},
\end{equation*} 
and define its lifted version $V_t:\Delta\to \Delta$ by
\begin{equation*}
V_{t}:=\boldsymbol L(F_t).
\end{equation*}
The fact that the range of $V_t$ is contained in $\Delta$ follows from the properties of $\boldsymbol L$ listed earlier.
Our aim is to evaluate the limit of the following difference as $t\to\infty$,
\begin{equation}\label{maindiff}
\ts\frac{1}{\ve_t}\big(\psi(F_t)-\psi(\id_{[0,1]})\big) = \ts\frac{1}{\ve_t}\Big(\mu(V_t^{-1}(\mathcal{S})) - \mu(\mathcal{S})\Big).
\end{equation}
Here, we have used the fact that $\psi(\id_{[0,1]})=\mu(\mathcal{S})$, which follows from $\boldsymbol L(\id_{[0,1]})=\id_{\Delta}$. Since $V_t$ approaches $\id_{\Delta}$ as $t\to\infty$, we may view the preimage $V_t^{-1}(\mathcal{S})$ as a perturbed version of the set $\mathcal{S}$. From this perspective, it is natural to interpret the right side of equation~\eqref{maindiff} through the lens of the first variation formula, introduced in Section~\ref{sec:contrib}, with $\mu$ playing the role of a volume, $\mathcal{S}$ playing the role of a manifold, and $V_t$ playing the role of the map $f_{\delta}$ with $\delta=\ve_t$.

In its classical form, the first variation formula deals with smooth maps on smooth manifolds. However, since the map $V_t$ need not be smooth, our proof proceeds by constructing a smoothed version of $V_t$. In order to do this, it is enough to smooth the univariate function $F_t$ and apply the linear operator $\boldsymbol L$.
The smoothing will be done using the linear \emph{Bernstein smoothing operator}, denoted $\mathcal{B}_s$, where $s\geq 1$ is an integer-valued smoothing parameter~\citep{lorentz,devoreconstructive}. 

For any function $G:[0,1]\to \R$, the Bernstein smoothing operator returns a new function $\mathcal{B}_s(G):[0,1]\to \R$ defined by
\begin{equation*}
\mathcal{B}_s(G)(u):= \sum_{j=0}^s G(j/s) \cdot b_{j}(u;s),
\end{equation*}
where $b_{j}(u;s) := \binom{s}{j}u^j(1-u)^{s-j}$ is the $j$th Bernstein basis polynomial where $u\in [0,1]$, and $j$ ranges over $\{0,1,\dots, s\}$. Below, we will use of some special properties of this operator. First, when $G$ is a cumulative distribution function, it turns out that $\mathcal{B}_s(G)$ is also a cumulative distribution function. Second, when $G$ is continuous, we have the uniform limit $\mathcal{B}_s(G)\to G$ as $s\to\infty$. The details of all the properties of $\mathcal{B}_s$ we will use are summarized in Lemma~\ref{bernprops} of Appendix~\ref{lemmas}.

When applying the operator $\mathcal{B}_s$, the smoothed version of $F_t$ will be denoted by
\begin{equation*}
F_{t,s}:=\mathcal{B}_s(F_t).
\end{equation*}
Likewise, the smoothed version of $V_t$ will be denoted by
\begin{equation}
V_{t,s}:=\boldsymbol L(F_{t,s}) = [\boldsymbol L\circ \mathcal{B}_s](F_t),
\end{equation}
which is a map from $\Delta$ to itself.
(It is not immediately obvious that $V_{t,s}$ takes values in $\Delta$, and this follows from  $F_{t,s}$ being a cumulative distribution function, by Lemma~\ref{bernprops}, as well as the ``Simplices'' property of $\boldsymbol L$.)

The remainder of the proof involves two essential parts. First, we prove a special version of the first variation formula for the smoothed maps $V_{t,s}$. 
Second, we show that this smoothing leads to negligible approximation error. 
To quantify the approximation error from smoothing, define the remainder $R_{t,s}$ according to the following equation
\small
\begin{equation}\label{rsn}
\frac{\mu(V_{t}^{-1}(\mathcal{S}))-\mu(\mathcal{S})}{\ve_t}= \frac{\mu(V_{t,s}^{-1}(\mathcal{S}))-\mu(\mathcal{S})}{\ve_t} +R_{t,s}.
\end{equation}
\normalsize
Due to the smoothness of $V_{t,s}$, the difference quotient on the right may be represented with a change of variable formula
\small
\begin{equation*}
 \frac{\mu(V_{t,s}^{-1}(\mathcal{S}))-\mu(\mathcal{S})}{\ve_t}
=\int_{\mathcal{S}^{\circ}} \frac{ f(V_{t,s}^{-1}(\theta)) |\det J(V_{t,s}^{-1})(\theta)| -f(\theta)}{\ve_t}d\theta,\label{changeofvar},
\end{equation*}
\normalsize
where $J(V_{t,s}^{-1})(\theta)$ is the Jacobian matrix of $V_{t,s}^{-1}$ at the point $\theta$.
This step is justified by Lemmas~\ref{bijective},~\ref{jacdensity}, and~\ref{meas0}, which also prove invertibility of $V_{t,s}$. From the above integral formula, Proposition~\ref{divprop} in Appendix~\ref{app:variation} provides the following limit
\small
\begin{equation*}\label{divthmpf}
\lim_{s\to\infty} \lim_{t\to\infty}\int_{\mathcal{S}^{\circ}} \frac{ f(V_{t,s}^{-1}(\theta)) |\det J(V_{t,s}^{-1})(\theta)| -f(\theta)}{\ve_t}d\theta
=\int_{\partial \mathcal{S}} \big\langle \!-\!\boldsymbol L(h)(\theta),\mathbf{n}(\theta)\big\rangle f(\theta) d\sigma(\theta).
\end{equation*}
\normalsize 
Next, Proposition~\ref{volumearg} in Appendix~\ref{app:smooth} shows that replacing $V_{t}^{-1}$ with its smoothed version $V_{t,s}^{-1}$ leads to negligible approximation error, i.e.
\begin{equation*}
\lim_{s\to\infty} \limsup_{t \to\infty} |R_{t,s}| = 0.
\end{equation*}
Consequently, by separately applying the operations $\liminf_{t\to\infty}$ and  $\limsup_{t\to\infty}$ to equation~\eqref{rsn}, and then taking $\lim_{s\to\infty}$ in each case, it follows that
\small
\begin{equation*}
\lim_{t\to\infty} \frac{\mu(V_{t}^{-1}(\mathcal{S}))-\mu(\mathcal{S})}{\ve_t} = \int_{\partial \mathcal{S}} \big\langle \!-\!\boldsymbol L(h)(\theta),\mathbf{n}(\theta)\big\rangle f(\theta) d\sigma(\theta).
\end{equation*}
\normalsize
Lastly, it is simple to check that the right side is a continuous linear functional of $h$, as required by the definition of Hadamard differentiability.
\qed

\section{A first variation formula}\label{app:variation}

\begin{proposition}\label{divprop}
Assume the conditions of Theorem~\ref{THM:HADAMARD}. Then, in the notation of Appendix~\ref{sec:highlevel}, the following limit holds,\\[0.2cm]
\small
\begin{equation}\label{divthm}
\lim_{s\to\infty} \lim_{t\to\infty} \int_{\mathcal{S}^{\circ}} \frac{ f(V_{t,s}^{-1}(\theta)) |\det J(V_{t,s}^{-1})(\theta)| -f(\theta)}{\ve_{t}}d\theta\\[0.3cm]
\ \ = \ \int_{\partial \mathcal{S}} \Big\langle \!-\!\boldsymbol L(h)(\theta),\mathbf{n}(\theta)\Big\rangle f(\theta) d\sigma(\theta).
\end{equation}
\normalsize
\end{proposition}
\proof  For any $\theta\in\Delta^{\circ}$, let
\begin{equation}\label{wdef}
W_{t,s}(\theta):=\ts\frac{1}{\ve_t}(V_{t,s}(\theta)-\theta),
\end{equation}
and
\begin{equation}\label{wtildedef}
\widetilde{W}_{t,s}(\theta):=\ts\frac{1}{\ve_t}(V_{t,s}^{-1}(\theta)-\theta).
\end{equation}
Using an expansion for $|\det J(V_{t,s}^{-1})(\theta)|$ given in Lemma~\ref{detexp} of Appendix~\ref{lemmas}, as well as the boundedness of $f$ on $\Delta$, the integral on the left side of equation~\eqref{divthm}
may be written as
\begin{equation}\label{writtenas}
\small
\int_{\mathcal{S}^{\circ}}\frac{ f(V_{t,s}^{-1}(\theta))  -f(\theta)}{\ve_{t}}d\theta \ 
-\int_{\mathcal{S}^{\circ}}\div \, W_{t,s}(V_{t,s}^{-1}(\theta))\, f(V_{t,s}^{-1}(\theta)) \, d\theta \
 +\mathcal{O}(\ve_t K_s),
\end{equation}
where $\div \, W_{t,s}(V_{t,s}^{-1}(\theta))$ is the divergence of $W_{t,s}$ evaluated at $V_{t,s}^{-1}(\theta)$, and $K_s\in [0,\infty)$ is a sequence of numbers not depending on $t$. (In addition, note that $V_{t,s}^{-1}(\theta)$ lies in $\Delta^{\circ}$ whenever $\theta\in\Delta^{\circ}$, which follows from  Lemma~\ref{bijective}, and the invariance of domain principle~\cite[Theorem 2B.3]{Hatcher}.) 

We now evaluate the limits of the two integrals in line~\eqref{writtenas} separately. %
Due to the multivariate mean value theorem, for each fixed $\theta\in\Delta^{\circ}$, there is a point $\zeta_{t,s}(\theta)\in \Delta^{\circ}$ on the line segment between $V_{t,s}^{-1}(\theta)$ and $\theta$ such that 
\begin{equation}\label{mvt}
\ts\frac{1}{\ve_t}\big(f(V_{t,s}^{-1}(\theta)) -f(\theta)\big)=\big\langle \nabla f(\zeta_{t,s}(\theta)),\widetilde{W}_{t,s}(\theta)\big\rangle.
\end{equation}
 Also, the points $\zeta_{t,s}(\theta)$ may be taken to satisfy $\zeta_{t,s}(\theta)\to \theta$ as $t\to\infty$, since $\lim_{t\to\infty} V_{t,s}^{-1}(\theta)= \theta$ for every fixed $\theta\in\Delta^{\circ}$ and fixed $s\geq 1$, which follows from Lemma~\ref{inversion}.
By the Cauchy-Schwarz inequality, the inner product in equation~\eqref{mvt} is dominated by a number depending only on $s$, since $\|\nabla f \|_2$ is bounded on $\Delta^{\circ}$ by assumption, and  $\sup_{t\geq 1} \sup_{\theta \in \Delta^{\circ}} \|\widetilde{W}_{t,s}(\theta)\|_2$ is finite by Lemma~\ref{inversion}. Furthermore, Lemma~\ref{inversion} gives the convergence $\lim_{t\to\infty}\widetilde{W}_{t,s}(\theta)\to -W_s(\theta)$ for every $\theta\in\Delta^{\circ}$, where we define
\begin{equation}W_s(\theta):=\boldsymbol L(\mathcal{B}_s(h))(\theta).\label{wsdef}
\end{equation} 
Hence, the continuity of $\nabla f$ and the dominated convergence theorem lead to 
\begin{equation*}
\lim_{t\to\infty} \int_{\mathcal{S}^{\circ}} \frac{ f(V_{t,s}^{-1}(\theta))  -f(\theta)}{\ve_{t}}d\theta \  = \ -\int_{\mathcal{S}^{\circ}} \big\langle \nabla f(\theta),W_s(\theta)\big\rangle d\theta.
\end{equation*}
\normalsize

 Turning our attention to the second integral in line~\eqref{writtenas}, Lemma~\ref{derivbound} ensures that the divergence 
\begin{equation*}\div \, W_{t,s}(\theta):=\sum_{l=1}^{k-1} \ts\frac{\partial}{\partial \theta_l}[W_{t,s}]_l(\theta)
\end{equation*} 
converges uniformly to $\div \, W_s(\theta)$ on $\Delta^{\circ}$ as $t\to\infty$. Combining this with the facts that $\lim_{t\to\infty}V_{t,s}^{-1}(\theta)=\theta$ for every $\theta\in\Delta^{\circ}$ (Lemma~\ref{inversion}) and that $\div(W_s)$ is continuous (since $W_s$ is smooth), it follows that
$$
\small
 \lim_{t\to\infty} \div \, W_{t,s}(V_{t,s}^{-1}(\theta)) \, f(V_{t,s}^{-1}(\theta)) \ = \ \div \, W_s(\theta)f(\theta).$$
Furthermore, it is simple to check that this pointwise limit is dominated by a constant (using Lemma~\ref{derivbound}), and so
\small
\begin{equation*}
\lim_{t\to\infty} \int_{\mathcal{S}^{\circ}} \div \,W_{t,s}(V_{t,s}^{-1}(\theta))\, f(V_{t,s}^{-1}(\theta))d\theta \ = \  \int_{\mathcal{S}^{\circ}}  \div \, W_s(\theta) f(\theta)d\theta.
\end{equation*}
\normalsize

The preceding calculations may now be combined using the basic differential identity
\begin{equation*}
\div\big(f\, W_s\big)(\theta) = \langle \nabla f(\theta),W_s(\theta)\rangle+ \div\big(W_s\big)(\theta) f(\theta),
\end{equation*}
which shows that the quantity in line~\eqref{writtenas} tends to $-\int_{\mathcal{S}^{\circ}} \div\big(f\, W_s\big)(\theta) d\theta$ as $t\to\infty$ with $s$ held fixed.
In turn, Stokes' theorem may be applied to this divergence integral. To be specific, an applicable version of the Stokes' theorem that holds for convex domains may be found in \citet[Theorem 9.64]{Leoni}. (When referencing that result, note that bounded convex domains have Lipschitz boundaries~\cite[Corollary 1.2.2.3]{Grisvard}, and also, that the regularity assumptions on $f$ imply that $f\,W_s$ is a Lipschitz vector field on $\Delta$.)
Hence, for every $s\geq 1$,
\small
\begin{equation}\label{stokes}
-\int_{\mathcal{S}^{\circ}} \div\big(f \,W_s\big)(\theta) d\theta = -\int_{\partial \mathcal{S}} \big\langle W_s(\theta),\mathbf{n}(\theta)\big\rangle f(\theta)d\sigma(\theta).
\end{equation}
\normalsize

Finally, since $h\in C[0,1]$, the uniform approximation property of Bernstein polynomials for continuous functions (Lemma~\ref{bernprops}) implies that 
$$W_s(\theta)=\boldsymbol L(\mathcal{B}_s(h))(\theta) \to \boldsymbol L(h)(\theta)$$
uniformly on $\Delta$ as $s\to\infty$. Hence, the right side of equation~\eqref{stokes} tends to $\int_{\partial \mathcal{S}} \langle - \boldsymbol L(h)(\theta),\mathbf{n}(\theta)\rangle f(\theta) d\sigma(\theta)$ by the dominated convergence theorem, which completes the proof.
\qed

\section{Smoothing error is negligible}\label{app:smooth}
\begin{proposition}\label{volumearg}
Let $R_{t,s}$ be as defined in equation~\eqref{rsn}, and suppose the conditions of Theorem~\ref{THM:HADAMARD} hold. 
Then,
\begin{equation*}
\lim_{s\to\infty} \limsup_{t\to\infty} |R_{t,s}|=0.
\end{equation*}
\end{proposition}
\vspace{-0.4cm}
\proof
It is simple to check that $R_{t,s}$ may be written as
\begin{equation}\label{diffs}
R_{t,s}=  \ts\frac{1}{\ve_t}\mu\big(V_t^{-1}(\mathcal{S})\setminus V_{t,s}^{-1}(\mathcal{S})\big) - \ts\frac{1}{\ve_t}\mu\big(V_{t,s}^{-1}(\mathcal{S})\setminus V_{t}^{-1}(\mathcal{S})\big).
\end{equation}
We will argue that both terms on the right side are small.
Note that the set $V_t^{-1}(\mathcal{S})\setminus V_{t,s}^{-1}(\mathcal{S})$ consists of the points $\theta\in \Delta$ such that $V_t(\theta)\in \mathcal{S}$ and $V_{t,s}(\theta)\not\in \mathcal{S}$. Now consider a particular point $\theta\in V_t^{-1}(\mathcal{S})\setminus V_{t,s}^{-1}(\mathcal{S})$, and note that if the Euclidean distance between the points $V_t(\theta)$ and $V_{t,s}(\theta)$ is written as $d_{t,s}(\theta)$, then both of the points $V_{t,s}(\theta)$ and $V_t(\theta)$ must be within a distance $d_{t,s}(\theta)$ from the boundary $\partial \mathcal{S}$. In other words, both of the points $V_t(\theta)$ and $V_{t,s}(\theta)$ lie within the tubular neighborhood of $\partial \mathcal{S}$ of radius $d_{t,s}(\theta)$. (The same reasoning applies to the other set $V_{t,s}^{-1}(\mathcal{S})\setminus V_{t}^{-1}(\mathcal{S})$.) 

To express the previous paragraph more formally, let $\mathcal{T}(\partial \mathcal{S}; r)\subset \R^{k-1}$ denote the tubular neighborhood of $\partial \mathcal{S}$ of radius $r\geq 0$, 
\begin{equation*}
\mathcal{T}(\partial \mathcal{S}; r):=\big\{\theta\in\R^{k-1} \, \big| \, d_2(\theta,\partial \mathcal{S})\leq r\big\},
\end{equation*}
where $d_2(\theta,\partial \mathcal{S})=\inf\{\|\theta-v\|_2 \, | \, v\in \partial \mathcal{S}\}$ is the Euclidean distance from a point $\theta$ to the boundary $\partial \mathcal{S}$.
Also, for any $\theta\in \R^{k-1}$, recall that
\begin{equation}\label{distance}
d_{t,s}(\theta) = \|V_{t,s}(\theta)-V_t(\theta)\|_2.
\end{equation}
Consequently, applying the reasoning above to both terms on the right side of equation~\eqref{diffs} gives,
\begin{equation}\label{tubeneighb}
|R_{t,s}| \ \leq \  2\cdot \ts\frac{1}{\ve_t}\,\mu\big\{\theta\in\Delta \ \big| \ V_{t,s}(\theta)\in \mathcal{T}\big(\partial \mathcal{S}; d_{t,s}(\theta)\big) \big\}.
\end{equation}
Furthermore, if $\Theta\in\Delta$ is a random vector distributed according to $\mu$, then this upper bound may be written as
\begin{equation*}
|R_{t,s}| \ \leq \  2\cdot \ts\frac{1}{\ve_t}\,\P\Big(V_{t,s}(\Theta)\in \mathcal{T}\big(\partial \mathcal{S}; d_{t,s}(\Theta)\big) \Big).
\end{equation*}
To modify the last bound somewhat, let
\begin{equation}\label{deltatsdef}
\delta_{t,s}:=\sup_{\theta\in\Delta}d_{t,s}(\theta),
\end{equation}
and put
\begin{equation}
r_{t,s}:=\max\{\delta_{t,s},\ve_t/s\}.\label{smallrtsdef}
\end{equation}
(This definition of $r_{t,s}$ will be used later on  for handling the formal possibility that $\delta_{t,s}$ may be zero. Note that $\ve_t/s$ is positive for all $s$ and $t$.)
Clearly, the definition of $r_{t,s}$ gives
\begin{equation*}
|R_{t,s}| \ \leq  \ 2\cdot \ts\frac{1}{\ve_t}\,\P\Big(V_{t,s}(\Theta)\in \mathcal{T}\big(\partial \mathcal{S}; r_{t,s}\big) \Big).
\end{equation*}
Next, we obtain an upper bound on this probability using a volume argument. Lemma~\ref{jacdensity} in Appendix~\ref{lemmas} shows that $V_{t,s}(\Theta)$ has a density, say $g_{t,s}$, with respect to $(k-1)$-dimensional Lebesgue measure. Also, Lemmas~\ref{bijective} and~\ref{meas0} imply that the random vector $V_{t,s}(\Theta)$ lies in $\Delta^{\circ}$ with probability 1 for all $t\geq 1$ and $s\geq 1$. Therefore,
\begin{equation}\label{mlbound}
|R_{t,s}| \ \leq \  2\cdot \ts\frac{1}{\ve_t} \Big(\displaystyle\sup_{\theta\in\Delta^{\circ}} g_{t,s}(\theta)\Big) \cdot \Big(\vol_{k-1}\!\mathcal{T}\big(\partial \mathcal{S}; r_{t,s}\big)\Big),
 \end{equation}
 where $\text{vol}_{k-1}$ denotes Lebesgue measure on $\R^{k-1}$.

  To control the volume of the right side of the bound~\eqref{mlbound}, it is convenient to consider the Hausdorff measure of $\partial \mathcal{S}$. Specifically, it is a fact from geometric measure theory that the $(k-2)$-dimensional Hausdorff measure of $\partial \mathcal{S}$, denoted $\mathcal{H}^{(k-2)}(\partial \mathcal{S})$, can be expressed as
\begin{equation}\label{hausdorff}
\mathcal{H}^{(k-2)}(\partial \mathcal{S}) = \lim_{r\to 0^+}\, \ts\frac{1}{2r}\vol_{k-1}\mathcal{T}(\partial \mathcal{S}; r).
\end{equation}
(See the beginning of Section 3.2.37 and  Theorem 3.2.39 in the book~\citep{federer} for additional details. In addition, that result depends on the fact that  bounded convex sets in $\R^{k-1}$ have $(k-2)$-rectifiable boundaries~\citep[p.743]{ambrosio}.) In order to make use of the expression above, we will rely on Lemma~\ref{distancelemma} in Appendix~\ref{lemmas},
 which states that for every $s\geq 1$, there is a  number $\kappa_s<\infty$ such that
\begin{equation*}
\lim_{t\to\infty}\ts\frac{1}{\ve_t}r_{t,s}= \kappa_s.
\end{equation*}
In particular,  $r_{t,s}$ is a sequence of positive numbers with $r_{t,s}\to0$ as $t\to\infty$, and so when $s$ is fixed, this means
\begin{equation}\label{minkowski}
\begin{split}
\lim_{t\to\infty} \ts\frac{1}{\ve_t} \vol_{k-1}\!\mathcal{T}\big(\partial \mathcal{S}; r_{t,s}\big) &=
\lim_{t\to\infty} \ts\frac{2 \, r_{t,s}}{\ve_t}\cdot \Big(\ts\frac{1}{2 \, r_{t,s}} \vol_{k-1}\!\mathcal{T}\big(\partial \mathcal{S}; r_{t,s}\big) \Big)\\[0.2cm]
&=2\kappa_s\cdot \mathcal{H}^{(k-2)}(\partial \mathcal{S}). 
\end{split}
\end{equation}
We now turn our attention to the factor $\sup_{\theta\in\Delta^{\circ}} g_{t,s}(\theta)$ in the bound~\eqref{mlbound}. Lemma~\ref{jacdensity} ensures there is a constant $c_0<\infty$, such that the following bound holds for every $s\geq 1$,
\begin{equation}\label{limsupdensity}
\limsup_{t\to\infty} \sup_{\theta\in\Delta^{\circ}} g_{t,s}(\theta) \leq c_0.
\end{equation}
Combining lines~\eqref{mlbound},~\eqref{minkowski}, and~\eqref{limsupdensity}, we conclude that for every $s\geq 1$,
\begin{equation*}
\limsup_{t\to\infty} |R_{t,s}| \leq 4c_0\cdot \kappa_s\cdot \mathcal{H}^{(k-2)}(\partial \mathcal{S}). 
\end{equation*}
Finally, the proof is completed using the fact that $\kappa_s\to 0$ as $s\to\infty$, which is shown in Lemma~\ref{distancelemma}.\qed

\section{Technical lemmas}\label{lemmas}

\paragraph{Remark} To simplify the statements of the lemmas in this section, the notation in the previous appendices will be generally assumed without comment.
\paragraph{Proof of equation~\eqref{sameerr}}Recall
$\bar{T}_t(\cdot):= \ts\frac{1}{t}\sum_{i=1}^t T_i(\cdot),$ and
that we may express $\Err_{t}$ and $\Err_{t}'$ as
\small
\begin{equation*}
\Err_{t} = \int_{\mathcal{X}\times \mathcal{Y}} 1\{{\tt{V}}(\bar{Q}_t(x))\neq y\}\, d\nu(x,y) \text{ \ \ \ and \ \ \ } \Err_{t}' =\int_{\mathcal{X}\times \mathcal{Y}} 1\{{\tt{V}}(\bar{T}_t(x))\neq y\}\, d\nu(x,y).
\end{equation*}
\normalsize
By Fubini's theorem,
\small
\begin{align*}
\E[\Err_{t}|\mathcal{D}] &= \int_{\mathcal{X}\times \mathcal{Y}} \P({\tt{V}}(\bar{Q}_t(x)) \neq y|\mathcal{D}) \,d\nu(x,y), \ \ \ \ \text{ and }\\[0.2cm]
\E[\Err_{t}'|\D] &=  \int_{\mathcal{X}\times\mathcal{Y}} \P({\tt{V}}(\bar{T}_t(x)) \neq y|\mathcal{D}) \,d\nu(x,y).
\end{align*}
\normalsize
Consequently, the claim $\E[\Err_{t}|\mathcal{D}] =\E[\Err_{t}'|\mathcal{D}] $ will follow if we can verify the equation
\begin{equation*}
\P({\tt{V}}(\bar{Q}_t(x)) \neq y \big| \mathcal{D}) = \P({\tt{V}}(\bar{T}_t(x)) \neq y \big| \mathcal{D}), 
\end{equation*}
for every fixed $(x,y)\in\mathcal{X}\times \mathcal{Y}$. This holds because when $x$ is fixed, we have $\mathcal{L}(\bar{Q}_t(x)|\mathcal{D}) = \mathcal{L}(\bar{T}_t(x)|\D)$, due to the first-order matching condition~\eqref{firsto}.\qed

\begin{lemma}[Properties of Bernstein polynomials]\label{bernprops} The Bernstein smoothing operator $\mathcal{B}_s$ satisfies the following properties.
	\begin{enumerate}
	\item For any function $h\in C[0,1]$, the following limit holds
	 \begin{equation*}
	 \sup_{u\in[0,1]}|\mathcal{B}_s(h)(u)-h(u)|\to 0 \text{ as } s\to\infty.
	 \end{equation*}
	 \item If $h:[0,1]\to\R$ is linear, then $\mathcal{B}_s(h)=h$ for all $s\geq 1$, and in particular  $\mathcal{B}_s(\textup{id}_{[0,1]})=\textup{id}_{[0,1]}$.
	 \item If $h:[0,1]\to\R$ is non-decreasing and non-constant, then for every $s\geq1$, the function $\mathcal{B}_s(h)$ is strictly increasing on $[0,1]$.
	\end{enumerate}
\end{lemma}
\proof We refer to the book~\citep{devoreconstructive} for general background on these properties. The  first property is given in Theorem 2.3 of \cite[Chapter 1]{devoreconstructive}. The  second property is proven after equation 1.7 of~\cite[Chapter 1]{devoreconstructive}. 

Regarding the third property, if we let $u\in (0,1)$ be arbitrary, it is enough to show that the derivative $\ts\frac{d}{du} \mathcal{B}_s(h)(u)$ is strictly positive. To this end, it is shown in equation 2.2 of~\cite[Chapter 10]{devoreconstructive} that $\mathcal{B}_s(h)$ satisfies
\begin{equation}\label{bernderiv}
\ts\frac{d}{du} \mathcal{B}_s(h)(u) = s \displaystyle\sum_{j=0}^{s-1} \underline{h}(\ts\frac{j}{s})\cdot  b_{j}(u;s-1),
\end{equation}
where we put $\underline{h}(u):=h(u+1/s)-h(u)$. If $h$ is non-decreasing and non-constant then
\begin{equation*}
0< h(1)-h(0)=\sum_{j=0}^{s-1}\underline{h}(\ts\frac{j}{s}).
\end{equation*}
So, because all the terms $\underline{h}(\ts\frac{j}{s})$ are non-negative, at least one of them must be positive. In turn, equation~\eqref{bernderiv} implies that $\ts\frac{d}{du} \mathcal{B}_s(h)(u)$ must be positive, since all the values $ b_{j}(u;s-1)$ are positive for all $u\in (0,1)$. \qed

\paragraph{Remark} Let $\boldsymbol D$ denote the differentiation operator on univariate functions. If this operator is applied to a function $g$ with domain $[0,1]$, then the domain of the derivative $\boldsymbol D(g)$ is taken to be $(0,1)$.

\begin{lemma}\label{bridgediff} Let $1_{(0,1)}$ be the indicator function of $(0,1)$. Then, for any fixed $s\geq 1$, we have the identity,
\begin{equation}\label{identity}
\ts\frac{1}{\ve_t}(F_{t,s}'-1_{(0,1)})=[\boldsymbol D\circ \mathcal{B}_s](\ts\frac{1}{\ve_t}(F_t - \textup{id}_{[0,1]})),
\end{equation}
and the uniform limit,
\begin{equation}\label{bridgedifflimit}
\lim_{t\to\infty} \ts\frac{1}{\ve_t}(F_{t,s}'-1_{(0,1)})\, = \, [\boldsymbol D\circ \mathcal{B}_s](h) \  \text{ \ in \ } \ell^{\infty}(0,1).
 \end{equation}

\end{lemma}
\proof The identity~\eqref{identity} is a direct consequence of $\mathcal{B}_s(\id_{[0,1]})=\id_{[0,1]}$ from Lemma~\ref{bernprops}.
To prove the limit, first note that because the functions $b_j(\cdot;s)$ are polynomials on $[0,1]$, the supremum \mbox{$C_s:=\displaystyle\max_{0\leq j\leq s} \sup_{u\in (0,1)} |b_{j}'(u;s)|$} is finite. Hence, for fixed $s$, 
\footnotesize
\begin{equation*}
\begin{split}
\sup_{u\in (0,1)} \Bigg|[\boldsymbol D\circ \mathcal{B}_s]&\Big(\ts\frac{1}{\ve_t}(F_t- \textup{id}_{[0,1]}) -h\Big)(u)\Bigg|\\[0.2cm]
&=\sup_{u\in (0,1)}\Bigg|\sum_{j=0}^s \Big( h_t(j/s)-h(j/s)\Big) b'_j(u;s)\Bigg|\\[0.2cm]
&\leq C_s \sum_{j=0}^s \big| h_t(j/s)-h(j/s)\big|\\[0.2cm]
&= o(1) \ \text{ as }  \ t\to\infty. 
\end{split}
\end{equation*}
\normalsize
\qed
\begin{lemma}\label{derivbound}
For any fixed $s\geq 1$, there is a number $K_s\in [0,\infty)$ not depending on $t$ such that the inequality 
 \begin{equation}\label{supderivbound}
 \sup_{t\geq1} \sup_{\theta\in \Delta^{\circ}} \big|\ts\frac{\partial}{\partial\theta_l}[ W_{t,s}]_l(\theta)\big|\leq K_s
 \end{equation}
 holds for all $l=1,\dots,k-1$.
 Furthermore, as $t\to\infty$, we have the uniform limit
 \begin{equation}\label{divunif}
 \text{\emph{div}} (W_{t,s})\xrightarrow{} \text{\emph{div}}(W_s) \text{ \ \ \ in \ \ \  } \ell^{\infty}[\Delta^{\!\circ}; \R],
\end{equation}
where $W_s=\boldsymbol L(\mathcal{B}_s(h))$.
\end{lemma}
\noindent \emph{Proof.} From the definition of $W_{t,s}$ in equation~\eqref{wdef}, we have
\begin{equation*}
\small
\ts\frac{\partial [W_{t,s}]_l(\theta)}{\partial \theta_l} = 
\ts\frac{1}{\ve_t}\Big(\ts\frac{\partial [V_{t,s}]_l(\theta)}{\partial \theta_l} - 1\Big)=\ts\frac{1}{\ve_t}\Big( F_{t,s}'(\theta_1+\dots+\theta_l) - 1\Big).
\end{equation*}
\normalsize
The last expression is bounded in absolute value for every $l=1,\dots,k-1$, and every $t \geq 1$, by 
\begin{equation*}
K_s:=\sup_{t\geq 1}\sup_{u\in(0,1)}\ts\frac{1}{\ve_t}\big| F_{t,s}'(u)-1\big|,
\end{equation*}
which is finite by the uniform limit in Lemma~\eqref{bridgediff}. (Hence, the bound~\eqref{supderivbound} is proved.)
Lastly, the limit~\eqref{divunif} follows from Lemma~\ref{bridgediff} and the definition $W_s=\boldsymbol L(\mathcal{B}_s(h))$. \qed

\begin{lemma}\label{bijective}
For any $t\geq 1$ and $s\geq 1$, the following three statements are true:
\begin{itemize}
\item[--] The map $V_{t,s}:\Delta\to\Delta$ is bijective and continuous on $\Delta$, and is also $C^1$ on $\Delta^{\circ}$. 
\item[--] The Jacobian matrix $J(V_{t,s})(\theta)$ is non-singular for all $\theta\in\Delta^{\!\circ}$.
\item[--]  The inverse map $V_{t,s}^{-1}:\Delta\to\Delta$ is $C^1$ on $\Delta^{\circ}$.
\end{itemize}
\end{lemma}

\proof It is simple to verify that $V_{t,s}$ is continuous on $\Delta$, and is $C^1$ on $\Delta^{\circ}$, due to the smoothness of $F_{t,s}$. To show that $V_{t,s}$ is bijective, it is enough to show that $F_{t,s}$ is bijective (due to the ``inverses'' property of $\boldsymbol L$ stated in Appendix~\ref{sec:lift}).
In turn, the fact that $F_{t,s}=\mathcal{B}_s(F_t)$ is bijective follows from part (c) of Lemma~\ref{bernprops}.

To see that the Jacobian matrix $J(V_{t,s})(\theta)$ is non-singular for $\theta\in\Delta^{\circ}$, it is enough to check that $V_{t,s}^{-1}$ is $C^1$ on the set $\Delta^{\circ}$, because we may  differentiate the identities 
\begin{equation*}
V_{t,s}\circ V_{t,s}^{-1} = \id_{\Delta} \text{ \ and \ } V_{t,s}^{-1}\circ V_{t,s} = \id_{\Delta},
\end{equation*}
to establish the inverse of $J(V_{t,s})(\theta)$ via the chain rule. Using the ``Inverses'' property of $\boldsymbol L$, we have $V_{t,s}^{-1}=\boldsymbol L(F_{t,s}^{-1})$, and it follows that $V_{t,s}^{-1}$ will be $C^1$ as long as $F_{t,s}^{-1}$ is. Finally, the fact that $F_{t,s}^{-1}$ is  $C^1$ follows from the strict monotonicity of $F_{t,s}$ and the univariate inverse function theorem.

\qed

\paragraph{Remark} 
The next lemma gives a uniform expansion for the determinant of $J(V_{t,s}^{-1})$ on $ \Delta^{\circ}$.
\begin{lemma}\label{detexp}
 For each \smash{$s\geq 1$} and $t\geq 1$, define the remainder function $\rho_{t,s}:\Delta^{\!\circ}\to \R$ to satisfy
\begin{equation*}
|\det J(V_{t,s}^{-1})(\theta)| = 1-\ve_t\, \text{\emph{\div}}\,W_{t,s}(V_{t,s}^{-1}(\theta))+\rho_{t,s}(\theta).
\end{equation*}
Then, for any $s\geq 1$, there is a number $K_s\in [0,\infty)$ not depending on $t$, such that the following bound holds for all large $t$,
\begin{equation*}
\sup_{\theta\in\Delta^{\circ}} |\rho_{t,s}(\theta)|\leq \ve_t^2\cdot \textup{poly}(K_s,\ve_t),
\end{equation*}
where $\textup{poly}(\cdot,\cdot)$ is a bivariate polynomial whose degree and coefficients do not depend on $t$ or $s$.
\end{lemma}

\proof 
It is simple to check that the Jacobian matrix $J(V_{t,s})(\theta)$ is lower-triangular for all $\theta\in\Delta^{\circ}$, and so the determinant of $J(V_{t,s})(\theta)$ is the product of the diagonal entries. Consequently, using the invertibility of $J(V_{t,s})(\theta)$ shown in Lemma~\ref{bijective}, we obtain the following expression for all $\theta\in\Delta^{\circ}$,
\begin{equation}\label{inversedet}
\det J(V_{t,s}^{-1})(\theta)= 1\big/\det \big(J(V_{t,s})(\theta_{t,s}')\big)= 1\Big/\ts \prod_{l=1}^{k-1}\ts\frac{\partial[V_{t,s}]_l(\theta_{t,s}')}{\partial \theta_l},
\end{equation}
where we put $\theta_{t,s}':=V_{t,s}^{-1}(\theta)$, which lies in $\Delta^{\circ}$ when $\theta$ does.
By the definition of $W_{t,s}$ in equation~\eqref{wdef}, we have for each $l=1,\dots,k-1$,
\begin{equation}\label{wreln}
\ts\frac{\partial[V_{t,s}]_l(\theta_{t,s}')}{\partial \theta_l} = 1+\ve_t\frac{\partial [W_{t,s}]_l(\theta_{t,s}')}{\partial \theta_l}.
\end{equation}
Due to Lemma~\ref{derivbound}, there is a number $K_s\in[0,\infty)$ not depending on $t$ such that the bound
\begin{equation}\label{rembound}
\sup_{t\geq 1}\,\sup_{\theta\in \Delta^{\circ}} \ts \big|\frac{\partial [W_{t,s}]_l(\theta)}{\partial \theta_l}\big| \leq K_s
\end{equation}
 holds for all $l$ simultaneously.  Next, define the function $\tilde{\rho}_{t,s}:\Delta^{\circ}\to\R$ so that the equation
\begin{equation}\label{bigproduct}
\begin{split}
\prod_{l=1}^{k-1} \Big(1+\ve_t\ts\frac{\partial [W_{t,s}]_l(\theta_{t,s}')}{\partial \theta_l}\Big) &=1+\ve_t\, \div\, W_{t,s}(\theta_{t,s}') +\tilde{\rho}_{t,s}(\theta_{t,s}')
\end{split}
\end{equation}
holds. Then, the bound~\eqref{rembound} implies
\begin{equation}\label{poly}
|\tilde{\rho}_{t,s}(\theta_{t,s}')| \leq \ve_t^2\cdot \text{poly}(\ve_t, K_s),
\end{equation}
where $\text{poly}(\cdot,\cdot)$ is a bivariate polynomial function whose degree and coefficients do not depend on $t$ or $s$. In particular, this upper bound does not depend on the point $\theta_{t,s}'$.
The proof is completed by combining lines~\eqref{inversedet} and~\eqref{bigproduct} with the elementary bound
\begin{equation*}
\big|\ts\frac{1}{1+\eta}-(1-\eta)\big|\leq 2\eta^2,
\end{equation*}
which holds for real numbers $|\eta|<\ts\frac{1}{2}$. Namely, we view $\eta$ as playing the role of $\ve_t \,\div \, W_{t,s}(\theta_{t,s}')+\tilde{\rho}_{t,s}(\theta_{t,s}')$. (Also note that taking the absolute value of $\det J(V_{t,s}^{-1})(\theta)$ does not matter asymptotically, because it approaches 1.)
\qed
\begin{lemma}\label{jacdensity} Assume the conditions of Theorem~\ref{THM:HADAMARD} hold, and let $\Theta$ be a random vector distributed according to $\mu$.  Then, the random vector $V_{t,s}(\Theta)$ has a density $g_{t,s}(\theta)$ with respect to Lebesgue measure on $\Delta^{\circ}$, given by
\begin{equation}\label{density}
g_{t,s}(\theta) = f(V_{t,s}^{-1}(\theta))|\det J(V_{t,s}^{-1})(\theta)|.
\end{equation}
Furthermore, the density $g_{t,s}$ is asymptotically bounded, in the sense that for each $s\geq 1$, we have
\begin{equation}\label{asympdensitybound}
 \limsup_{t\to\infty} \sup_{\theta\in\Delta^{\circ}} g_{t,s}(\theta)\leq \|f\|_{\infty},
 \end{equation}
 where $\|f\|_{\infty}:=\sup_{\theta\in\Delta} f(\theta)$.
\end{lemma}

\proof Lemma~\ref{bijective} and the standard change of variables formula (\cite[Theorem 2.47]{folland}) give the stated expression for $g_{t,s}(\theta)$. The boundedness condition~\eqref{asympdensitybound} follows from Lemmas~\ref{derivbound} and~\ref{detexp}, as well as the boundedness of $f$ on $\Delta$.\qed
\begin{lemma}\label{meas0}
Suppose the conditions of Theorem~\ref{THM:HADAMARD} hold, and let $A\subset \Delta$ be a convex set. Then, 
\begin{equation*}
\mu(\partial A)=0,
\end{equation*}
and for all $s,t\geq 1$,
\begin{equation*}
\mu(V_{t,s}^{-1}(\partial A))=0.
\end{equation*}

\end{lemma}
\proof 
The first equation follows from the fact that the boundary of any convex set in $\R^{k-1}$ has Lebesgue measure zero~\cite[Lemma 2.4.3]{dudleyclt}, and $\mu$ is assumed to have a density with respect to Lebesgue measure on $\Delta$.  The second equation follows similarly using Lemma~\ref{jacdensity}, and the fact that $V_{t,s}^{-1}$ maps $\partial \Delta$ to $\partial \Delta$. \qed
\begin{lemma}\label{distancelemma}
Let $r_{t,s}$ be as defined in equation~\eqref{smallrtsdef}. Then, there is a sequence of numbers $\kappa_s\in[0,\infty)$ such that
\begin{equation*}
\lim_{t\to\infty}\ts\frac{1}{\ve_t}r_{t,s} = \kappa_{s},
\end{equation*}
and furthermore
 \begin{equation*}
 \lim_{s\to\infty} \kappa_s = 0.
 \end{equation*}
\end{lemma}
\proof

Suppose we can show there are numbers $\tilde{\kappa}_s\in[0,\infty)$ such that $\frac{1}{\ve_t} \delta_{t,s}\to\tilde{\kappa}_s$ as $t \to\infty$, where $\delta_{t,s}$ is defined in equation~\eqref{deltatsdef}. Then, the definition of $r_{t,s}$ gives
\begin{equation*}
\begin{split}
\lim_{t\to\infty}\ts\frac{1}{\ve_t} r_{t,s} &= \lim_{t\to\infty} \max \{ \ts\frac{1}{\ve_t} \delta_{t,s}, 1/s\}
= \max\{ \tilde{\kappa}_s, 1/s\}.
\end{split}
\end{equation*}
Hence, it is enough to show that such numbers $\tilde{\kappa}_s$ exist, and that $\tilde{\kappa}_s\to 0$ as $s\to\infty$.

To proceed, recall the identities
\begin{align*}
V_{t,s} - \id_{\Delta} &= \boldsymbol L\circ \mathcal{B}_s(F_t-\id_{[0,1]})\\[0.2cm]
V_t-\id_{\Delta} &= \boldsymbol L(F_t-\id_{[0,1]}).
\end{align*}
Letting $\boldsymbol I$ denote the identity operator on $\ell^{\infty}[0,1]$, it is straightforward to check that the linear operator $\boldsymbol L(\mathcal{B}_s-\boldsymbol I)$ is a continuous map from $\ell^{\infty}[0,1]$ to $\ell^{\infty}[\Delta;\R^{k-1}]$.
Consequently, holding $s$ fixed, we have the following uniform limit in $\ell^{\infty}[\Delta; \R^{k-1}]$,
\begin{equation*}
\begin{split}
\ts\frac{1}{\ve_t}\big(V_{t,s}-V_t \big)
&= \ \ \, \ts\frac{1}{\ve_t}\big(V_{t,s} - \id_{\Delta} \big)-\ts\frac{1}{\ve_t}\big(V_t-\id_{\Delta}\big)\\[0.3cm]
&= \ \ \, \boldsymbol L(\mathcal{B}_{s}-\boldsymbol I)\big[\ts\frac{1}{\ve_t}\big(F_t-\id_{[0,1]}\big)\big]\\[0.3cm]
&  \xrightarrow{ \ t\to\infty   } \ \  \boldsymbol L(\mathcal{B}_{s}-\boldsymbol I)h.
\end{split}
\end{equation*}
It follows that as $t\to\infty$,
\begin{equation}\label{Mlimit}
\sup_{\theta\in \Delta} \ts\frac{1}{\ve_t}\big\|V_{t,s}(\theta) - V_t(\theta)\big\|_2 \ \xrightarrow{ \ \ \  \ }\  \displaystyle\sup_{\theta\in\Delta}\big\|\big[\boldsymbol L(\mathcal{B}_{s}-\boldsymbol I)h\big](\theta)\big\|_2=: \tilde{\kappa}_s, 
\end{equation}
which proves the desired numbers $\tilde{\kappa}_s$ exist. To show that $\tilde{\kappa}_s$ tends to 0 as $s\to\infty$, note that since $h$ lies in $C[0,1]$, the uniform approximation property of Bernstein polynomials on $C[0,1]$ in Lemma~\ref{bernprops} ensures that as $s\to\infty$,
\begin{equation*}
(\mathcal{B}_s-\boldsymbol I)h \xrightarrow{} 0 \ \text{ in } \ C[0,1].
\end{equation*}
Consequently, by using the continuity of the operator $\boldsymbol L$, we have $\tilde{\kappa}_s\to 0.$\qed
~\\

\begin{lemma}\label{inversion}

Let $s\geq 1$ be fixed. Then, the following limits hold respectively in the spaces $\ell^{\infty}[\Delta;\R^{k-1}]$ and $\ell^{\infty}[\Delta^{\circ};\R^{k-1}]$ as $t\to\infty$:
\begin{equation}\label{firstmap}
W_{t,s}=\ts\frac{1}{\ve_t}(V_{t,s}-\text{\emph{id}}_{\Delta})\xrightarrow{ \ \ } \boldsymbol L\circ \mathcal{B}_s(h),
\end{equation}
and 
\begin{equation}\label{secondmap}
\widetilde{W}_{t,s}=\ts\frac{1}{\ve_t}(V_{t,s}^{-1}-\text{\emph{id}}_{\Delta^{\circ}})\xrightarrow{ \ \ } -\boldsymbol L\circ \mathcal{B}_s(h).
\end{equation}
\end{lemma}
\proof 
The limit~\eqref{firstmap} follows from the identity 
$$\ts\frac{1}{\ve_t}\big(V_{t,s} - \id_{\Delta}\big) = \boldsymbol L\circ \mathcal{B}_s(\ts\frac{1}{\ve_t}\big(F_t-\id_{[0,1]})\big),$$
 and the continuity of $\boldsymbol L\circ \mathcal{B}_s$.

 To prove the second limit~\eqref{secondmap}, let  $\Psi$ denote the map from $\mathcal{F}[0,1]$ to $\ell^{\infty}(0,1)$, that sends a distribution function $G$ to its generalized inverse $G^{-1}$, defined by
 $G^{-1}(y):=\inf\{x \, | \, G(x)\geq y\}$.
 Also note that for every $s,t\geq 1$, the function $F_{t,s}$ lies in $\mathcal{F}[0,1]$ by Lemma~\ref{bernprops}, and so it makes sense to evaluate $\Psi$ on $F_{t,s}$. 
 
 Due to Lemma 3.9.23 in the book~\citep{vaartWellner},\footnote{Note that the space $\mathcal{F}[0,1]$ is denoted $\mathbb{D}_2$ in Lemma 3.9.23 of the book \citet{vaartWellner}} the map $\Psi$ is Hadamard differentiable at $\id_{[0,1]}\in \mathcal{F}[0,1]$ tangentially to $C[0,1]$, and the derivative acts on $C[0,1]$ as the negative of the identity operator, i.e.  $\Psi_{\id_{[0,1]}}'=-\boldsymbol I$. Hence, using the linearity of $\mathcal{B}_s$, and the property \mbox{$\mathcal{B}_s(\id_{[0,1]})=\id_{[0,1]}$} from Lemma~\ref{bernprops}, as well as the chain rule for Hadamard differentiation, we obtain
\begin{equation*}
\begin{split}
\ts\frac{1}{\ve_t}(F_{t,s}^{-1}-\id_{(0,1)})& \, \ = \ \ts\frac{1}{\ve_t}\Psi(F_{t,s})-\Psi(\id_{[0,1]}))\\[0.2cm]
& \, \ = \ \ts\frac{1}{\ve_t}\big(\Psi(\mathcal{B}_s(F_t))-\Psi(\mathcal{B}_s(\id_{[0,1]}))\big)\\[0.2cm]
& \xrightarrow{ t\to\infty } \ \Psi_{\id_{[0,1]}}'(\mathcal{B}_s(h))\\[0.2cm]
& \, \ =  \ -\mathcal{B}_s(h),
\end{split}
\end{equation*}
where the limit is with respect to the sup-norm on $\ell^{\infty}(0,1)$, which is the codomain of $\Psi$. Applying the operator $\boldsymbol L$ to the previous limit, and using the ``Inverses'' property $V_{t,s}^{-1}= [\boldsymbol L(F_{t,s})]^{-1}=\boldsymbol L(F_{t,s}^{-1})$, we have
\begin{equation*}
\ts\frac{1}{\ve_t}(V_{t,s}^{-1}-\id_{\Delta^{\circ}}) = \ts\frac{1}{\ve_t}\boldsymbol L\big(F_{t,s}^{-1}-\id_{(0,1)}\big)\xrightarrow{ \ t\to\infty \ } -\boldsymbol L(\mathcal{B}_s(h)),
\end{equation*}
which completes the proof.
\qed

\section{Assessment of Assumption~\ref{ASSUMPTION}}\label{app:assumption}

In the first portion of this section, we provide theoretical support for Assumption~\ref{ASSUMPTION} in the context of two types of ensemble methods: the voting Gibbs classifier, and bagged decision stumps. Later on, we also provide empirical justification in the context of random forests. 

\subsection{Theoretical assessment of Assumption~\ref{ASSUMPTION}}
 Before dealing with examples of specific ensemble methods, we first give a general result concerning the existence of the density $f_l$ in Assumption~\ref{ASSUMPTION}. In essence, the following proposition shows that $f_l$ exists when the function $\vartheta$ is sufficiently smooth. 

\begin{proposition}\label{prop:assumption} Let $\mathcal{X}=\R^p$ with $p\geq k$, and suppose  the test point distribution \smash{$\mathcal{L}(X|Y=\boldsymbol e_l)$} has a density $\rho_l$ with respect to Lebesgue measure on $\R^p$. In addition, suppose the function \smash{$\vartheta(x)=\E[Q_1(x)|\D]$} is $C^{p-k+2}$ on $\mathcal{X}$. Then, the distribution $\mathcal{L}(\vartheta(X)|\mathcal{D},Y=\boldsymbol e_l)$ has a density $f_l$ with respect to Lebesgue measure on $\Delta$, and for almost every $\theta\in\Delta$, the density is given by 
\begin{equation}\label{coarea}
f_l(\theta)=\int_{\vartheta^{-1}(\theta)} \ts\frac{\rho_l(x)}{\sqrt{\det(J(\vartheta)(x)J(\vartheta)(x)\!\ttop)}} d\textup{vol}_{\vartheta^{-1}(\theta)}(x),
\end{equation}
where $J(\vartheta)(x)$ is the Jacobian matrix of $\vartheta$ evaluated at $x$, the region of integration is the pre-image $\vartheta^{-1}(\theta)=\{x\in \mathcal{X}: \vartheta(x)=\theta\}$,  and $d\textup{vol}_{\vartheta^{-1}(\theta)}(x)$ refers to $(p-k+1)$-dimensional Hausdorff measure on $\vartheta^{-1}(\theta)$. 
\end{proposition}

\proof The result is a consequence of the co-area formula and Sard's Theorem. The details may be found by combining Theorem 10.4 and line 10.6 in the book~\citep{simon}.\qed

\paragraph{Remarks}Beyond the existence of $f_l$, Assumption~\ref{ASSUMPTION} also requires the gradient of $f_l$ to bounded and continuous. However, given that the general formula~\eqref{coarea} for $f_l$ is quite complex, the analysis of the gradient of $f_l$ seems to be prohibitive. For this reason, we focus primarily on the existence of $f_l$ in the examples below --- by analyzing the smoothness of $\vartheta$. Indeed, even verifying the smoothness of $\vartheta$ is non-trivial in general.

\subsubsection{Voting Gibbs classifier}
The voting Gibbs classifier is a Bayesian ensemble method for binary classification, $k=2$~\citep{ngjordan}. 
Suppose that $\mathcal{X}=\R^p$, and that for each $l\in\{0,1\}$, the class-wise test point distribution $\mathcal{L}(X|Y=\boldsymbol  e_l)$ has a parametric density $\rho_{l,\beta}(x)$ with respect to Lebesgue measure, where $\beta$ denotes a vector of parameters. Also, let $p(\beta|\D)$ denote a posterior distribution for $\beta$, and define the function 
\begin{equation}\label{psidef}
\eta_{\beta}(x):=\P(Y=1 |X=x)=\frac{\pi_1 \rho_{1,\beta}(x)}{\pi_0\rho_{0,\beta}(x)+\pi_1\rho_{1,\beta}(x)},
\end{equation}
where $\pi_l=\P(Y=\boldsymbol e_l)$ is the $l$th class proportion.
The voting Gibbs classifier works by drawing i.i.d.~samples $\beta_1,\dots,\beta_t$ from the posterior $p(\beta|\D)$, as well as an independent set of i.i.d.~variables $U_1,\dots,U_t$ from Uniform[0,1]. Next, each pair $(\beta_i,U_i)$ gives rise to a classifier defined by
$$Q_i(x) := 1\big\{U_i\leq \eta_{\beta_i}(x)\big\},$$
which is to say that $Q_i$ randomly labels $x$ as 1 with probability $\eta_{\beta_i}(x)$. The classifiers $Q_1,\dots,Q_t$ are then aggregated via majority voting.

In this context, averaging out the algorithmic randomness corresponds to integrating over the posterior, as well as the uniform variables $U_1,\dots,U_t$. Consequently, the function $\vartheta(x)=\E[Q_1(x)|\D]$ can be represented as

$$\vartheta(x)=\int \eta_{\beta}(x) p(\beta|\D)d\beta.$$
In turn, the smoothness of $\vartheta(x)$ will be inherited from the smoothness of $\eta_{\beta}(x)$ via equation~\eqref{psidef}. For example, in the case of Bayesian logistic regression, we have 
$$\eta_{\beta}(x)=\frac{1}{1+\exp(-w\ttop x-b)},$$
where the parameter vector is written as $\beta=(w,b)\in\R^{p+1}$. Likewise, all of the mixed partial derivatives of $\vartheta(x)$ with respect to $x=(x_1,\dots,x_p)$ will exist, provided that these derivatives can be applied to the smooth function $1/(1+\exp(-w\ttop x-b))$ under the integral
$$\vartheta(x)=\int \ts\frac{1}{1+\exp(-w\ttop x-b)}p(\beta|\D)d\beta,$$
and it can be checked that this is permitted (for instance) when $p(\beta|\D)$ is continuous in $\beta$, and is supported on a compact rectangular domain.

\subsubsection{Bagged decision stumps} 
Recall that in the context of bagging, the classifiers $Q_1,\dots,Q_t$ are randomized, conditionally on $\D$, by training them on bootstrapped datasets $\D_1^*,\dots,\D_t^*$. Since the number of possible bootstrapped versions of $\D$ is finite, it follows that the function $\vartheta(x)=\E[Q_1(x)|\D]$ represents a finite average. Hence, if each $Q_i(x)$ is a non-smooth function of $x$, then the finite average $\vartheta(x)$ will generally be non-smooth in a strict sense.

Nevertheless, if we write of $\vartheta_n(x)\equiv \vartheta(x)$ to reflect the dependence of $\vartheta$ on the sample size $n=|\D|$, then it is possible to argue that $\vartheta_n(x)$  can become smooth as $n\to\infty$.
For example, in the seminal paper \citep{buhlmannyu}, the authors consider  a classifier of the form \smash{$1\{\hat{d}_n\leq x\}$}, where $\mathcal{X}=\R$, and $\hat{d}_n$ is a ``split point'' that is estimated from $\D$. (Such a classifier can be viewed as a one-level decision tree, and is known as a ``decision stump''.) Hence, if the values $\hat{d}_{n,1}^*,\dots,\hat{d}_{n,t}^*$ are obtained from bootstrapped datasets $\D_{1}^*,\dots,\D_{t}^*$, then the associated ensemble of bagged classifiers is given by
 $$Q_i(x):=1\{\hat{d}_{n,i}^*\leq x\},$$
 and $\vartheta_n(x)$ represents an average over all bootstrap samples, with
  \begin{equation}\label{thetanrep}
 \vartheta_n(x)=\P(\hat{d}_{n,1}^*\leq x|\D).
 \end{equation}
In this situation, the statement below formalizes the asymptotic smoothness of $\vartheta_n(x)$, and is a slight reformulation of Proposition 2.1 in the paper~\citep{buhlmannyu}.  The significance of this fact is that it allows the asymptotic smoothness of $\vartheta_n(x)$ to be understood in terms of the limit of the standardized bootstrap distribution,  $\mathcal{L}(\sqrt{n}(\hat{d}_n^*-\hat{d}_n)|\D)$, which can be derived analytically in special cases.

\begin{proposition}[{{\cite[Proposition 2.1]{buhlmannyu}}}]  Under the conditions described above, suppose there is a distribution function \smash{$G:\R\to [0,1]$,} and an increasing sequence of positive numbers $\{b_n\}$, such that as $n\to\infty$,
 \begin{equation}\label{BYassumption}
  \sup_{v\in\R}\bigg|\P\big(b_n(\hat{d}_n^*-\hat{d}_n)\leq v\big|\D\big) - G(v)\bigg| =o_{\P_{\D}}(1).
 \end{equation}
Then, as $n\to\infty$, 
 \begin{equation*}
\sup_{x\in\R} \bigg|\vartheta_n(x)-G\big(b_n(x-\hat{d}_n)\big)\bigg|=o_{\P_{\D}}(1).
 \end{equation*}
 \end{proposition}

 \paragraph{Remarks} In essence, the proposition states that if the condition~\eqref{BYassumption} holds with a smooth limiting distribution function $G$,  then $\vartheta_n(x)$ will also be asymptotically smooth, which leads to the question of verifying~\eqref{BYassumption}. As noted in the paper~\cite[p.933]{buhlmannyu}, when the estimate $\hat{d}_n$ depends on the training data in a smooth way, and when the observations $\D=\{(X_j,Y_j)\}_{j=1}^n$ are i.i.d., one often expects that $\mathcal{L}(\sqrt{n}(\hat{d}_{n,1}^*-\hat{d}_n)|\D)$ will have a smooth (Gaussian) limiting distribution function. Similarly, when $\mathcal{X}=\R^p$, this reasoning extends to classification based on discriminant functions, such as $1\{\hat{w}_n\ttop x+\hat{b}_n\leq 0\}$, with estimated coefficients $(\hat{w}_n,\hat{b}_n)\in\R^{p+1}$, and in fact, an analogue of assumption~\eqref{BYassumption} can be verified with a Gaussian limit in the case of logistic regression~\citep{lee_logistic}. However, in more general situations, the condition~\eqref{BYassumption} can be difficult to verify --- for instance, when the bootstrapped law $\mathcal{L}(b_n(\hat{d}_{n,1}^*-\hat{d}_{n,1})|\D)$ does not asymptotically agree with $\mathcal{L}(b_n(\hat{d}_{n,1}-\E[\hat{d}_{n,1}]))$. Indeed, in this type of situation, the question of verifying the asymptotic smoothness of $\vartheta_n(x)$ was left as an open problem by B\"uhlmann and Yu~\cite[p.941]{buhlmannyu}, and they offer the following remark about the ability of bagging to work as a smoothing operation:\\

 ``It is worth noting that [bootstrap consistency] is not necessary for bagging to work as long as the resulting bagged estimator is sensible itself. Conditional on the original sample, $\hat d_n^*$ spreads around [its population counterpart] by taking one of the discrete values between original sample points. The resulting bagged stump estimator is a weighted average of the stump estimators with split points between the original sample values. Thus, bagging is still a smoothing operation, similar to the assertion in Proposition 2.1, although exact analysis seems difficult and we leave it as an open research problem.''\\

\subsection{Empirical assessment of Assumption~\ref{ASSUMPTION}}\label{app:empirical}

Here, we empirically assess Assumption~\ref{ASSUMPTION} by seeing how well \smash{$\mu_l=\mathcal{L}(\vartheta(X)|\mathcal{D},Y=\boldsymbol e_l)$} can be approximated by a distribution with a smooth density function. For convenience, we only consider the situation of binary classification, because in this case,  $\mu_l$ is a univariate distribution on [0,1], which simplifies the assessment of goodness-of-fit.

A natural class of smooth distributions on [0,1] is the Beta$(\alpha,\beta)$ family, parameterized by $\alpha,\beta>0$.
For any fixed $\tau\in[0,1]$, the densities in this family are given by
$$g(\tau; \alpha,\beta)=\frac{1}{B(\alpha,\beta)} \tau^{\alpha-1}(1-\tau)^{1-\beta},$$
where $B(\alpha,\beta)$ is the Beta function.

 The main idea of these experiments is to generate approximate samples from $\mu_l$, and then see how well these samples can be fit by a member of the Beta$(\alpha,\beta)$ family. Noting that $\mu_l$ depends on a particular training set $\mathcal{D}$, we will consider three instances of $\mu_l$ arising form the datasets `census income', `synthetic discrete', and `synthetic continuous' from the main text. \

For each of the three datasets, we prepared the training set $\mathcal{D}$ and the test set $\mathcal{D}_{\text{ground}}$ as described in Section~\ref{sec:expt}. To generate approximate samples from $\mu_l$, we first approximated the function $\vartheta$ using the sample average $\bar{Q}=\ts\frac{1}{t}\sum_{i=1}^tQ_i$, obtained from an ensemble of size $t=1,000$, trained on $\D$, via the package~{\tt{randomForest}} with default settings~\citep{randomForestsCitation}. Next, letting $X_{1,l}',\dots,X_{r_l,l}'$ denote the samples from class $l$ in the test set $\mathcal{D}_{\text{ground}}$, we used the values $\bar{Q}_t(X'_{1,l}),\dots,\bar{Q}_t(X_{r_l,l}')$ as approximate samples from $\mu_l$. In turn, these approximate samples were used to estimate $\alpha$ and $\beta$ via the method of moments, using the `mme' option in the package {\tt{fitdistrplus}}~\citep{fitdistrplus}. Below, we write $\hat{\alpha}_l$ and $\hat{\beta}_l$ to refer to the estimates associated with $\mu_l$.

To assess the quality of fit, we constructed quantile-quantile (QQ) plots by sorting the values $\bar{Q}_t(X_{1,l}'),\dots,\bar{Q}_t(X_{r_l,l}')$ and plotting them against a corresponding set of quantiles from the fitted distribution Beta($\hat{\alpha}_l,\hat{\beta}_l$), with the results shown below.
Overall, the plots indicate a good fit, with strong conformity to the diagonal line $y=x$.

~\\
~\\

\begin{figure*}[h!]
\centering
{\includegraphics[angle=0,
  width=.46\linewidth]{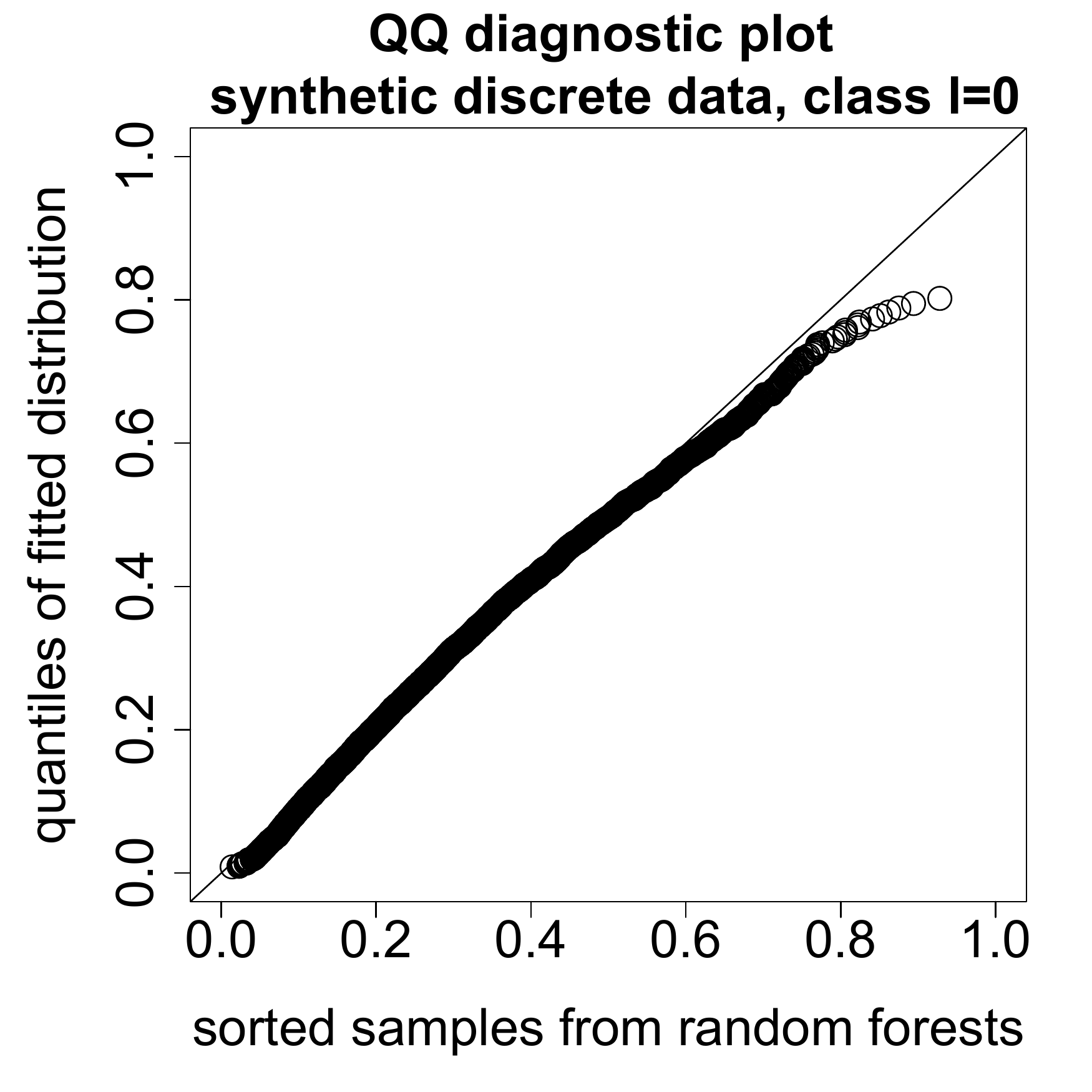}} \ \ \ \ 
{\includegraphics[angle=0,
  width=.46\linewidth]{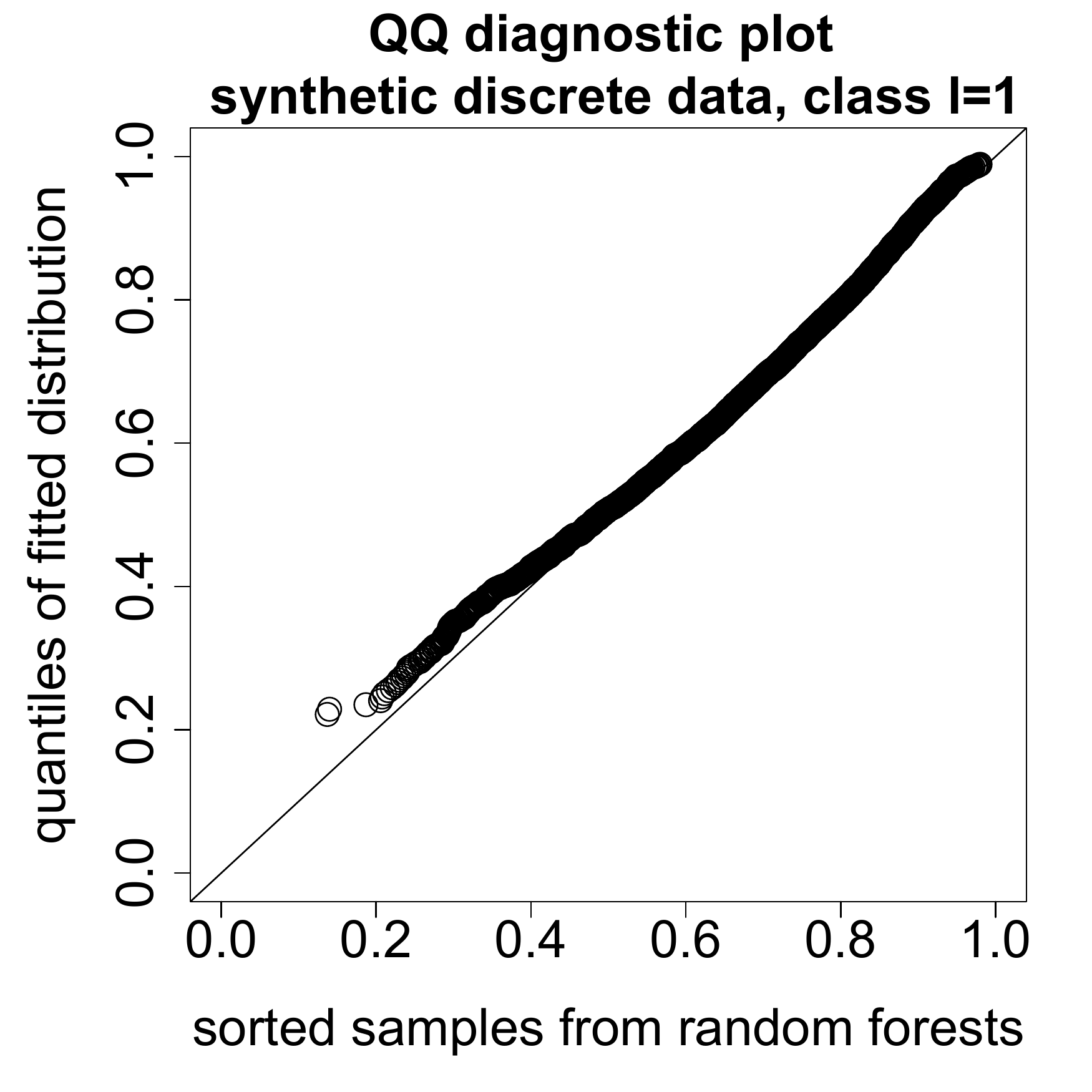}}

\caption{QQ plots for synthetic discrete data.}
\label{fig:contraceptive}
~\\
~\\
\end{figure*}

~\\
~\\
~\\
~\\
~\\
~\\
\begin{figure*}[h!]
\centering
{\includegraphics[angle=0,
  width=.46\linewidth]{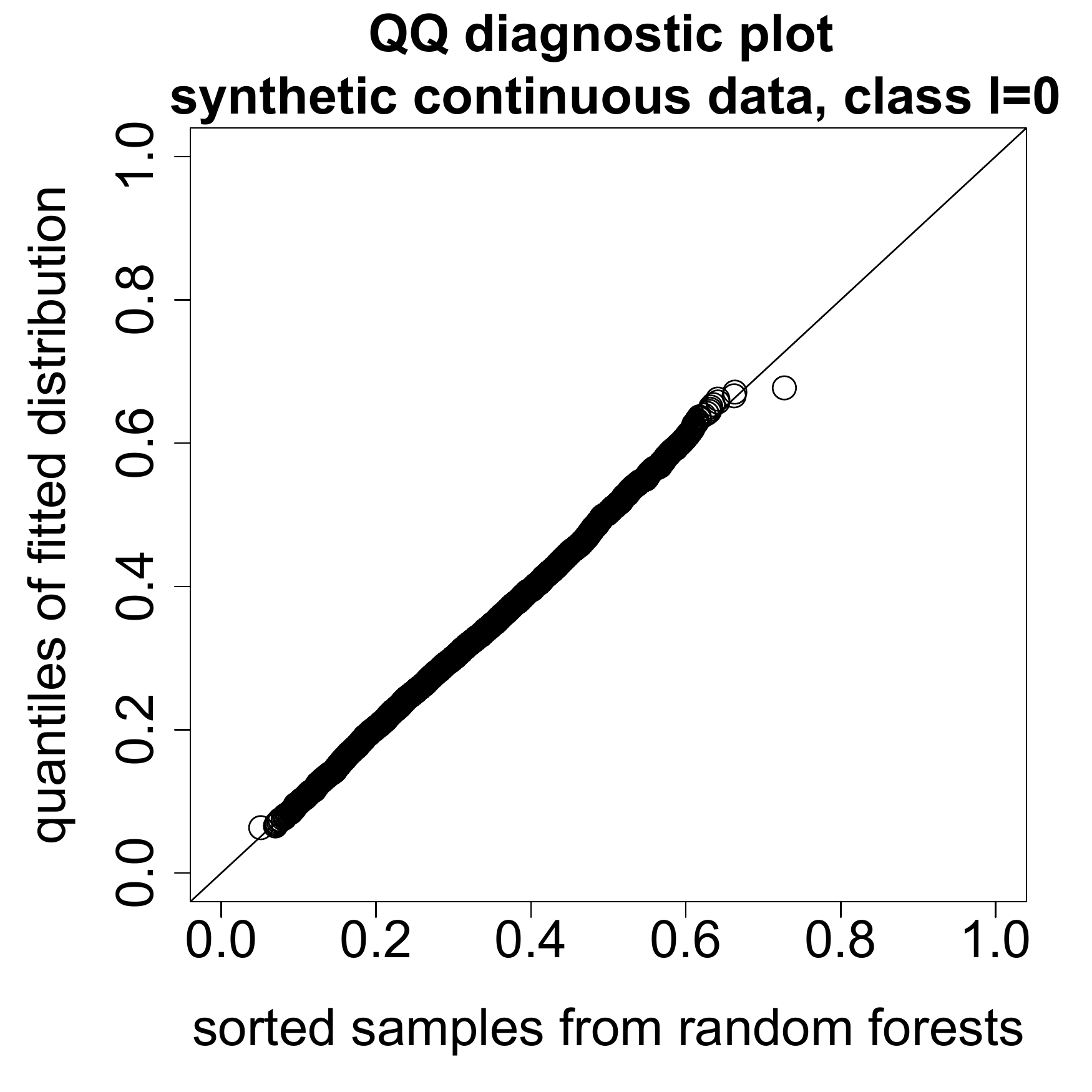}} \ \ \ \ 
{\includegraphics[angle=0,
  width=.46\linewidth]{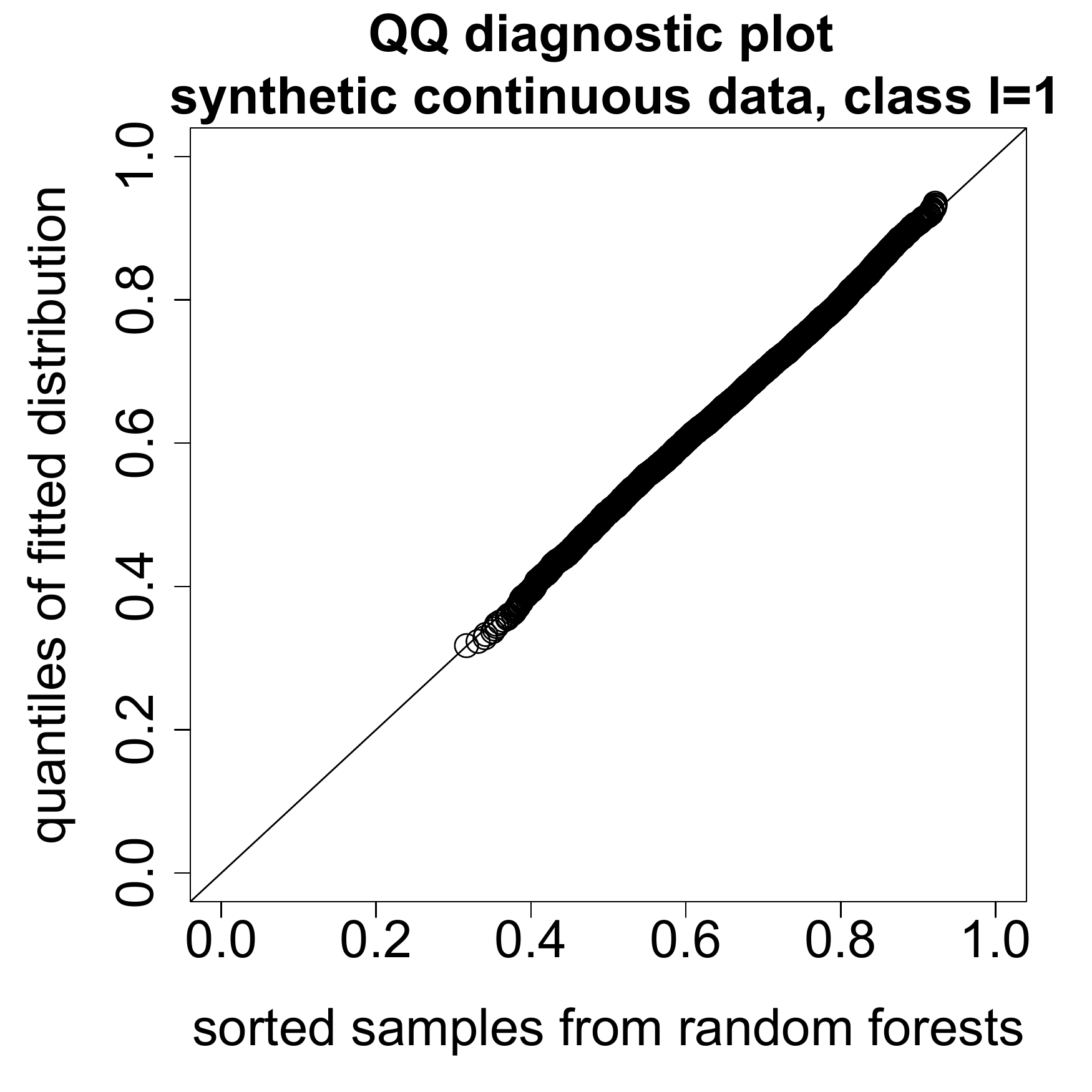}}

\caption{QQ plots for synthetic continuous data.}
\label{fig:contraceptive}
~\\
~\\
\end{figure*}


~\\
~\\
~\\

\begin{figure*}[h!]
\centering
{\includegraphics[angle=0,
  width=.43\linewidth]{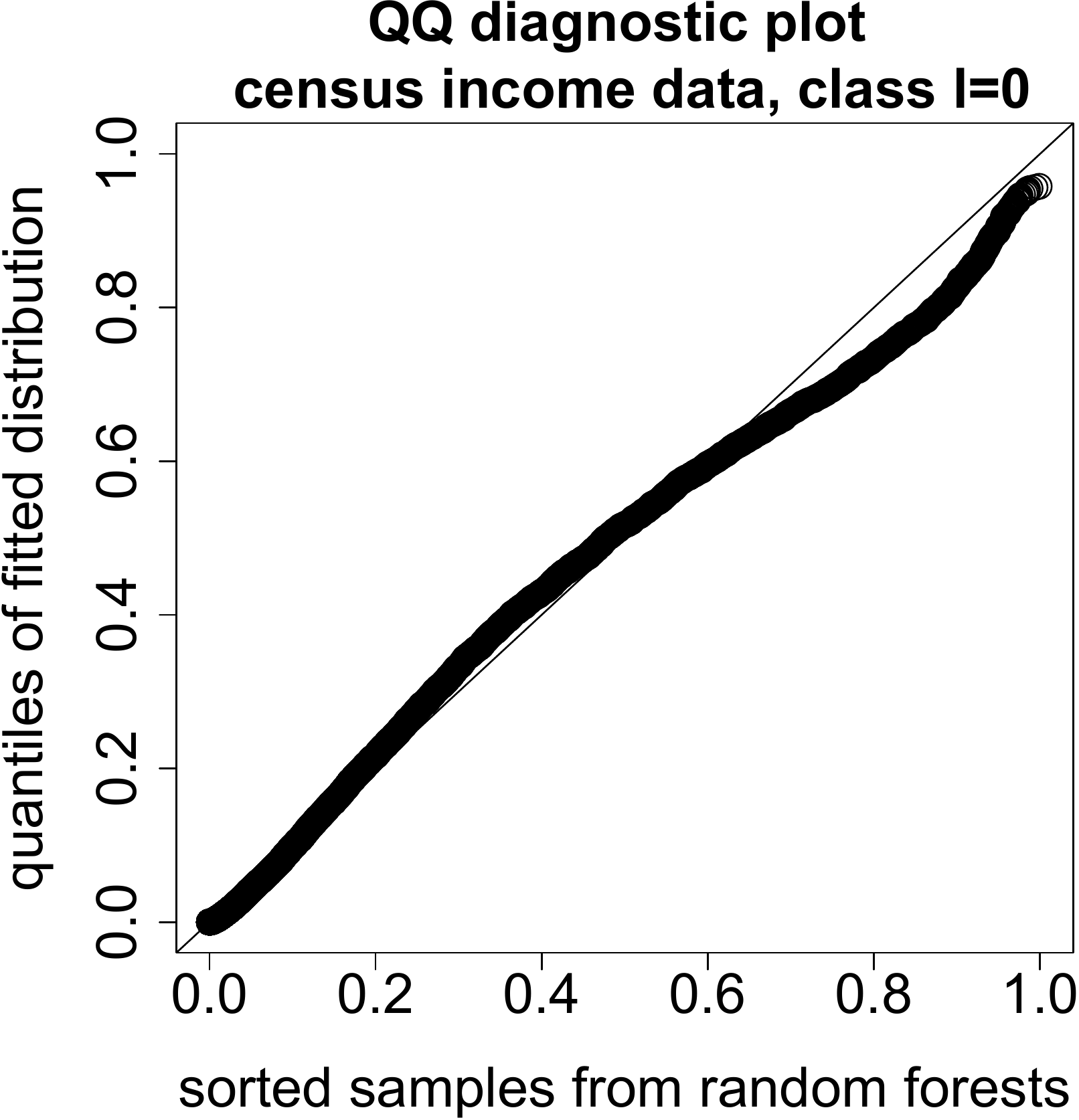}} \ \ \ \ \ \ \
{\includegraphics[angle=0,
  width=.43\linewidth]{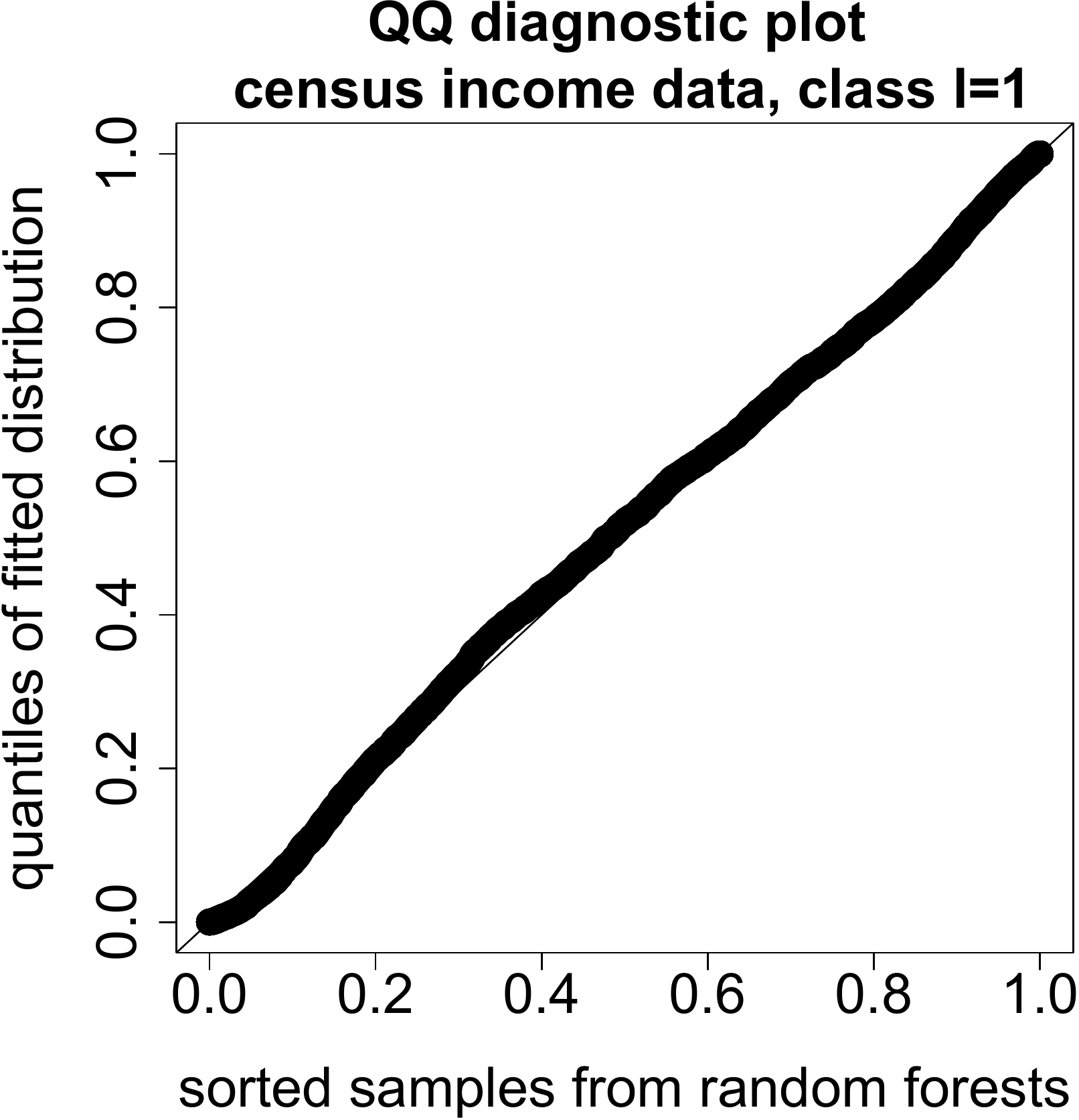}}

\caption{QQ plots for census income data.}
\label{fig:contraceptive}
\end{figure*}


\end{document}